\documentclass[a4paper,10pt]{article} 
\usepackage{amsmath, amssymb, esint, comment}
\usepackage[a4paper,
        left=2.6cm,
        right=2.6cm,
        top=2.1cm,
        bottom=2.2cm]{geometry}
\usepackage{amscd, accents}
\usepackage{shadow, epsfig, color}
\usepackage{graphicx}
\usepackage{url}
\usepackage{mathrsfs}
\usepackage[final]{changes_0}
\usepackage{xcolor,cancel}
\usepackage{authblk}
\usepackage{fancyhdr}
\usepackage{pb-diagram, environ}

\usepackage[dvipdfm,
colorlinks=true,
linkcolor=black,
filecolor=black,
urlcolor=black,
citecolor=black
]{hyperref}

\newcommand{\N}{\mathbb N}

\newcommand{\Q}{\mathbb Q}
\newcommand{\R}{\mathbb R}

\newcommand{\e}{\varepsilon}
\newcommand{\1}{\mathbf 1}
\newcommand{\qed}{\ \hfill \fbox{} \bigskip}
\newcommand{\proof}{{\it Proof.} }
\newcommand{\dis}{\displaystyle}

\newcommand{\E}{\mathcal E}

\newcommand{\F}{\mathcal F}
\newcommand{\x}{\overline{x}}

\newcommand{\EN}{\overline{\N}}

\newcommand{\red}{\textcolor{red}}

\newcommand{\CC}{\accentset{\circ}}
\newcommand{\OM}{\overline{\mathcal M(N,K,D)}}

\renewcommand{\b}{\mathbf{b}}

\numberwithin{equation}{section}

\newtheorem{thm}{Theorem}[section]
\newtheorem{defn}[thm]{Definition}
\newtheorem{lem}[thm]{Lemma}
\newtheorem{rem}[thm]{Remark}
\newtheorem{cor}[thm]{Corollary}
\newtheorem{exa}[thm]{Example}
\newtheorem{prop}[thm]{Proposition}

\renewenvironment{abstract}
 {
  \begin{center}
  \bfseries \abstractname\vspace{-.5em}\vspace{0pt}
  \end{center}
  \list{}{%
    \setlength{\leftmargin}{7mm}
    \setlength{\rightmargin}{\leftmargin}%
  }%
  \item\relax}
 {\endlist}

\makeatletter
\newcommand{\subjclass}[2][2010]{%
  \let\@oldtitle\@title%
  \gdef\@title{\@oldtitle\footnotetext{#1 \emph{Mathematics Subject Classification.} #2}}%
}
\newcommand{\keywords}[1]{%
  \let\@@oldtitle\@title%
  \gdef\@title{\@@oldtitle\footnotetext{\emph{Key Words and Phrases.} #1.}}%
}
\makeatother

\renewcommand{\tilde}{\widetilde}

\title{Regularity and Stability of Invariant Measures for Diffusion Processes under Synthetic Lower Ricci Curvature Bounds
}
  \author{Kohei Suzuki \thanks{Scuola Normale Superiore di Pisa, Piazza dei Cavalieri 7,  56126 Pisa PI, Italy. E-mail: kohei.suzuki@sns.it}
}
            
\subjclass{Primary 53C23; Secondary 35J15.}
\keywords{Invariant Measure, Synthetic Lower Ricci Curvature Bound, Regularity of Elliptic PDE}

\begin{document}
\date{}
\maketitle
\begin{abstract}
The Sobolev regularity of invariant measures for diffusion processes is proved on non-smooth metric measure spaces with synthetic lower Ricci curvature bounds. 
As an application, the symmetrizability of semigroups is characterized, and 
the stability of invariant measures is proved under perturbations of drifts and the underlying spaces in the sense of the measured Gromov convergence.   
\end{abstract}
\tableofcontents

\section{Introduction}
Invariant measures describe equilibrium states of (random) dynamical systems, which play a fundamental role for understanding the long-time behaviour of their evolutions.   
For a given measure space $(X, m)$ and a strongly continuous semigroup $\{T_t\}$ on $L^p(X,m)$, a measure $\mu=\rho m$ with $\rho \in L^q(X,m)$ $(1/p+1/q=1)$ is said to be {\it an invariant measure for $\{T_t\}$} if 
\begin{align} \label{eq: invT5} 
\int_{X}T_t f d\mu=\int_X f d\mu, \quad \forall f \in L^p(X,m), \ \forall t \ge 0.
\end{align}
Invariant measures can also be characterized in terms of the infinitesimal generator.  Let $(L, \mathcal D(L))$ be the infinitesimal generator on $L^p(X,m)$ associated with $\{T_t\}.$ Then $\mu$ is an invariant measure for $\{T_t\}$ if, and only if 
\begin{align} \label{eq: invL0}
\int_{X}L\phi d\mu=0, \quad \forall \phi \in \mathcal D(L).
\end{align}
In this case, $\mu$ is also called an invariant measure for $(L, \mathcal D(L))$.
Equation \eqref{eq: invL0} comes down to a solution $\rho$ solving the following elliptic PDE:
\begin{align} \label{eq: invL1}
L^* \rho=0,
\end{align}
whereby $(L^*, \mathcal D(L^*))$ denotes the adjoint operator of $(L, \mathcal D(L))$ on $L^p(X,m)$. 

The existence, uniqueness, and regularity of invariant measures are fundamental questions both for the theory of Markov processes and for the theory of elliptic PDEs. There is a long history concerning these issues from both aspects: Metafune--Pallara--Rhandi \cite{MPR05} on the Euclidean space, Bogachev, Krylov, R\"ockner and their coauthors \cite{BKR96}--\cite{BRS00} for Riemannian manifolds and some related infinite-dimensional spaces; Shigekawa \cite{Shi87} and Hino \cite{H98-2} for the Wiener space; Bogachev--R\"ockner--Zhang \cite{BRZ00} for locally convex topological vector spaces; Hino \cite{H00} and Beznea--Cimpean--R\"ockner \cite{BCR15} for general state spaces.  

The aim of this paper is to study invariant measures and the related elliptic PDE \eqref{eq: invL1} on non-smooth metric measure spaces satisfying synthetic lower Ricci curvature bounds (RCD condition for brevity). The notion of the RCD condition on metric measure spaces is a generalization of the notion of lower Ricci curvature bounds in the framework of metric measure spaces, which are stable under geometric convergences such as the measured Gromov--Hausdorff (GH) convergence. They therefore contain various finite- and infinite-dimensional singular spaces such as Ricci limit spaces (Sturm \cite{Sturm06, Sturm06-2}, Lott--Villani \cite{LV09}), Alexandrov spaces (Petrunin \cite{Pet11}, Zhang--Zhu \cite{ZZ10}), warped products and cones (Ketterer \cite{Ket14, Ket14a}), quotient spaces (Galaz-Garc\'ia--Kell--Mondino--Sosa \cite{GKMS17}), stratified spaces (Bertrand--Ketterer--Mondello--Richard \cite{BKMR18}) and infinite-dimensional spaces such as Hilbert spaces with log-concave measures (Ambrosio--Savar\'e--Zambotti \cite{ASZ09}), which are related to various stochastic partial differential equations. In spite of the singularities of spaces, various analysis and geometries have been well-developed from the synthetic viewpoint of lower Ricci curvature bounds. 
 
We first establish the construction of diffusion processes and the existence and uniqueness of invariant measures. 
The infinitesimal generator considered in this paper can be written formally in the following form: 
\begin{align} \label{L: EO}
L=\frac{1}{2}\Delta+\b.
\end{align}
Here $\Delta$ denotes the Laplacian, and $\b$ denotes a derivation operator (also called {\it drift} in the context of probability theory), which is the first-order differential operator. In Section \ref{sec: CSD}, we construct diffusion processes and the associated semigroups by using the theory of Dirichlet forms and the Girsanov transformation (cf.\ Hino \cite{H98-2} in the Wiener space; Fitzsimmons \cite{F07} in general spaces).  In order to develop the Girsanov transformation in the generality of RCD spaces, we construct an isometry $I: L^2(TX) \to \accentset{\circ}{\mathcal M}$ whereby $L^2(TX)$ denotes the $L^2$-tangent module (Gigli \cite{G16}) and $\accentset{\circ}{\mathcal M}$ denotes the space of square-integrable martingales with finite energy associated with the Cheeger energy.  Then, for $\b^* \in L^2(TX)$, we define an exponential martingale:
\begin{align} \label{def: exp}
M_t=\exp\biggl\{N_t-\frac{1}{2}\langle N\rangle_t\biggr\},
\end{align}
 whereby $N=I(\b^*)$ and $\langle N \rangle$ denotes the quadratic variation process of $N$, and define the semigroup $T_t f(x)=\mathbb E_x[M_tf(B_t)]$, which is the Girsanov transformation of the Brownian motion $(B_t, \mathbb P_x)$ associated with $(1/2)\Delta$. By aid of the log-Sobolev inequality under the RCD condition, and the exponential integrability assumption on $\b$, the semigroup $\{T_t\}$ can be extended to a strongly continuous $L^p$-semigroup with $1<p<\infty$ and, when $p=2$,  $\{T_t\}$ coincides with the $L^2$-semigroup associated with the Dirichlet form constructed in Proposition \ref{prop: form}. In Section \ref{sec: EU}, by making use of the construction of semigroups in Section \ref{sec: CSD}, we show the existence and the uniqueness of invariant measures under several different assumptions.

As the main result of this paper, the Sobolev regularity of the density $\rho$ and the gradient estimate of $\rho$ are obtained in Section \ref{sec: NSDF} whenever $\rho$ satisfies the following equation:
 \begin{align} \label{eq: invL1-3}
\int_X \Bigl(\frac{1}{2}\Delta+\b\Bigr)(\phi) d\mu =0, \quad \forall \phi \in {\rm TestF}(X).
\end{align}
Note that \eqref{eq: invL1-3} is weaker than \eqref{eq: invL0} since ${\rm TestF}(X)$ is not necessarily dense in the domain of $(L, \mathcal D(L))$ with respect to the graph norm, but we also call $\mu=\rho m$ {invariant measure}.
Since RCD spaces are quite singular, the standard strategy based on the Friedrichs mollifier does not work for the Sobolev regularity of $\rho$. However, the heat semigroup $\{H_t\}$ possibly works as a good mollifier owing to the RCD condition, which enables us to show the Sobolev regularity of $\rho$. {Compared to the proof in the Euclidean case, we need to handle the integrability issue of the mollifier more carefully (see Theorem \ref{thm: mthm1-2} and Remark \ref{rem: 2}).} 
Note that, in the Wiener space, Shigekawa \cite{Shi87} proved the Sobolev regularity by using the finite-dimensional approximation of the Wiener space. Our proof relies only on regularity of the heat semigroup based on the Bakry--\'Emery gradient estimate,  which gives a different proof in the Wiener space.   

As an application of the Sobolev regularity, we can answer the question on the symmetrizability of the semigroup $\{T_t\}$ on $L^2(X,m)$. The question is based on the gap between the symmetrizing measures and the invariant measures. We say that $\{T_t\}$ is symmetrizable if there is a Borel measure $\mu=\rho m$ with $\rho \in L^1(X,m)$ so that 
\begin{align} \label{def: symm}
\int_X (T_t f) g d\mu=\int_X f T_t g d\mu, \quad f, g \in \mathcal B_b(X),
\end{align}
whereby $\mathcal B_b(X)$ denotes the set of real-valued bounded Borel-measurable functions on $X$. 
If $\{T_t\}$ is stochastically complete (i.e., $T_t\1=\1$ for any $t \ge 0$), it is easy to see that the symmetry \eqref{def: symm} with respect to $\mu$  implies the invariance \eqref{eq: invT5} with respect to $\mu$. The converse implication is, however, false in general. A natural question is, therefore, to ask 
\begin{itemize}
\item[(Q1)] 
when is $\{T_t\}$ symmetrizable?
\end{itemize} 
This question was considered by Kolmogorov \cite{K37}, originally in the context of probability theory, in order to investigate when time-reversed diffusion processes have the same generators as the original diffusion processes. In compact smooth Riemannian manifolds with smooth drifts, $\{T_t\}$ is symmetrizable if, and only if $\b$ is a gradient type, i.e., $\b(\cdot)=\langle \nabla f, \nabla \cdot \rangle$ for some smooth function $f$ (Kolmogorov \cite{K37}, and see also Ikeda--Watanabe \cite[Section 4 in Chapter V]{IW89}). In the case of infinite-dimensional spaces, (Q1) was considered by Shigekawa \cite{Shi87} in the Wiener space, and by Albeverio--R\"ockner \cite{AR95} in locally convex topological vector spaces (see also the references provided therein). 
In Theorem \ref{thm: mthm3}, under the RCD condition, we show that $\{T_t\}$ is symmetrizable if, and only if, there is a Lipschitz continuous function $f$ with $f \in L^2(m)$ so that $\b(\cdot)=\langle \nabla f, \nabla \cdot \rangle$. 
This result can be seen as a generalization of Kolmogorov \cite{K37} and Shigekawa \cite{Shi87} to non-smooth settings.
Another application of the Sobolev estimate is the stability of invariant measures under perturbations of the derivation operator $\b$ and the underlying space $X$.
Equality \eqref{eq: invL0} tells us that invariant measures are determined by the Laplacian $\Delta$ and the drift $\b$, whereby the Laplacian is determined by the underlying metric measure structure. 
A natural question regarding the stability is the following: 
\begin{itemize}
\item[(Q2)] 
If the underlying space $X_n$ converges to $X_\infty$ in a geometric sense, and the drift $\b_n$ converges to $\b_\infty$ in an operator sense, then does the invariant measure $\mu_n$ converge to $\mu_\infty$ weakly? 
\end{itemize} 
In this paper, the pointed measured Gromov (pmG) convergence and the $L^2$-convergence are adopted as a geometric convergence and a convergence of drifts in (Q2) respectively. In this case, it should be noted that the pmG convergence cannot a priori see any differential structure of $X_n$ (like the uniform convergence of functions), but invariant measures depend on the second-order differential structure (or the first-order differential structure in the weak form) by \eqref{eq: invL0}. Therefore, (Q2) is unreasonable without some uniform control of the underlying differential structure, for which the RCD condition comes into play.  

In Theorem \ref{thm: mthm4}, we prove the stability of $\mu_n$ and $\rho_n$ under the RCD condition. To be more precise,  we prove that $\mu_n$ is stable with respect to the weak convergence of measures, and that $\rho_n$ is stable with respect to the convergence of the Sobolev norm under the pmG convergence of $X_n$ and the $L^2$-convergence of $\b_n$.  
If the underlying spaces are fixed, i.e., $X_n=X_\infty$ for all $n \in \N$, then Theorem \ref{thm: mthm4} implies that $\mu_n$ converges to $\mu_\infty$ in the total variation distance.   
The key steps for the proof of Theorem \ref{thm: mthm4} are the regularity of $\rho$ obtained in Theorem \ref{thm: mthm1} and \ref{thm: mthm1-2},  and showing the convergence of non-symmetric Dirichlet forms under varying underlying spaces (cf.\ Hino \cite{H98}, T\"olle \cite{T06, T10}, and \cite{S18}). 

\section*{Acknowledgment}
The author appreciates Professor Shouhei Honda for valuable comments on the regularity of elliptic PDEs. 
He thanks Lorenzo Dello Schiavo for fruitful discussions on metric measure geometries and Proposition \ref{lem: Dom1}. He is indebted to Anna Katharina Suzuki-Klasen for her careful reading of his manuscript. 
This work was supported by World Premier International Research Center Initiative (WPI), MEXT, Japan, and JSPS Grant-in-Aid for Scientific Research on Innovative Areas ``Discrete Geometric Analysis for Materials Design": Grant Number 17H06465, and by JSPS Overseas Research Fellowships Grant number 290142.

\section{Notation and Preliminary Results} \label{sec: Pre}
\subsection{Preliminary from Metric Measure Geometry} \label{subsec: Pre}
Let $\N=\{0,1,2,...\}$ and $\overline{\N}=\N \cup \{\infty\}$ be the set of natural numbers and the set of extended natural numbers respectively.  
Let $(X,d)$ be a complete separable metric space. The open ball centered at $x \in X$ with radius $r>0$ is denoted by $B_r(x)=\{y \in X: d(x,y)<r\}$. The family of all Borel sets in $(X,d)$ is denoted by $\mathscr B(X)$. The set of real-valued bounded Borel-measurable functions on $X$ is written as $\mathcal B_b(X)$.
Let $C(X)$ denote the set of real-valued continuous functions on $X$, while $C_b(X), C_0(X)$ and $C_{bs}(X)$ denote the subsets of $C(X)$ consisting of bounded functions, functions with compact support, and bounded functions with bounded support respectively. {Let ${\rm Lip}(X)$ denote the set of real-valued Lipschitz continuous functions on $X$. The set of continuous functions on $[0,\infty)$ valued in $X$ is denoted by $C([0,\infty); X)$.  Let ${\rm Lip}_b(X)$ and $\mathrm{Lip}_{bs}(X)$ denote the subsets of ${\rm Lip}(X)$ consisting of bounded functions, and bounded functions with bounded supports, respectively.  For $f \in {\rm Lip}_{bs}(X)$, {\it the global Lipschitz constant ${\rm Lip}(f)$} is defined as the infimum of $h>0$ satisfying $|f(x)-f(y)| \le hd(x,y)$ for any $x,y \in X$.} 
The set $\mathcal P(X)$ denotes all Borel probability measures on $X$. A sequence of probability measures $\mu_n \in \mathcal P(X)$ is said to converge to $\mu$ if, for any $f \in C_b(X)$, it holds that $\mu_n(f) \to \mu(f)$ as $n \to \infty$. Here $\mu(f)=\int_X f d\mu.$

For a given Borel measure $m$, the support of $m$ is denoted by $\mathrm{supp}[m]=\{x \in X: m(B_r(x))>0, \ \forall r>0\}$.
Let $(Y,d_Y)$ be another complete separable metric space. For a Borel measurable map $f: X \to Y$, let $f_\#m$ denote the push-forward measure on $Y$:
$${f}_\#m(B) = m(f^{-1}(B)) \quad \text{for any Borel set} \quad  B \in \mathscr B(Y).$$
For a measurable set $A\subset X$, the indicator function is denoted by $\1_A$, which is equal to $1$ for $x \in A$ and $0$ otherwise. For  two functions $f, g: X \to \R$, their maximum and minimum are denoted by $f\vee g=\max\{f,g\}$ and $f \wedge g=\min\{f,g\}.$
Let $L^0(X,m)$ denote the set of all $m$-equivalence classes of Borel measurable functions $f: X \to \R$. Let $1 \le p \le \infty$.
A function $f \in L^0(X,m)$ is said to belong to $L^p(X,m)$ if $\|f\|_{p}^p=\int_X|f|^pdm<\infty$ for  $1 \le p <\infty$, and if $\|f\|_{\infty}={\text {ess-sup}}_{x \in X}|f(x)|<\infty$ for $p=\infty$. Let $L^p_+(X,m)$ (resp.\ $L^p_{>0}(X,m)$) denote the set of non-negative (resp.\ strictly positive) $L^p$-functions.
For $f,g \in L^2(X,m)$, the round bracket $(f,g)$ denotes the inner product $\int_Xfg dm.$ 

Throughout the paper, when we write $(X,d,m)$, we always assume that 
\begin{align} \label{cond: fund} 
&\text{ $(X,d)$ is a complete separable geodesic metric space}  \notag
\\
&\text{and $m$ is a Borel measure with $m(X)<\infty$ and ${\rm supp}[m]=X$.}
\end{align}
We simply call $(X,d,m)$ {\it a metric measure space} if \eqref{cond: fund} is satisfied.

A curve $\gamma:[0,1] \to X$ is said to be {\it absolutely continuous} if there exists a function $f \in L^1(0,1)$ so that 
\begin{align} \label{ineq: ABC}
d(\gamma_t, \gamma_s) \le \int_t^sf(r)dr, \quad \forall t,s \in [0,1], \quad t<s.
\end{align}
{\it The metric speed} $t \mapsto |\dot{\gamma}|_t \in L^1((0,1), dt)$ is defined as the essential infimum among all the functions $f$
satisfying \eqref{ineq: ABC}.
A Borel probability measure $\pi$ on $C([0,1];X)$ is said to be {\it a test plan} if there exists a non-negative constant $C(\pi)$ so that 
$$(e_t)_{\#}\pi \le C(\pi)m, \quad \forall t \in [0,1], \quad \text{with} \quad \int\int_0^1|\dot{\gamma}_t|^2dt d\pi(\gamma)<\infty.$$
Here $e_t(\gamma)=\gamma(t) \in X$ denotes the evaluation map. 
{\it The set of Sobolev functions} $S^2(X,d,m)$ (or, simply $S^2(X)$) is defined to be the space of all functions in $L^0(X,m)$ so that there exists a non-negative $G \in L^2(m)$ for which it holds
$$\int|f(\gamma_1)-f(\gamma_0)|d\pi(\gamma) \le \int\int_0^1G(\gamma_t)|\dot{\gamma}_t|dt d\pi(\gamma), \quad \text{for any test plan } \pi.$$ 
It turns out (see \cite{AGS14}), that for $f \in S^2(X)$ there exists a minimal $G$ in the $m$-a.e.\ sense so that the above inequality holds, which is denoted by $|\nabla f|$ and called {\it minimal weak upper gradient}. {\it The Sobolev space} is defined by $W^{1,2}(X,d,m)=S^2(X,d,m) \cap L^2(X,m)$ (or, simply $W^{1,2}(X)$).
{\it The Cheeger energy} ${\sf Ch}: W^{1,2}(X,d,m) \to \R$ is defined in the following manner: 
\begin{align*}
{\sf Ch}(f)&=\frac{1}{2}\int_X|\nabla f|^2dm, \quad f \in W^{1,2}(X,d,m).
\end{align*}
For $f \in {\rm Lip}(X)$, the local Lipschitz constant is defined as follows:
\begin{align*}|{\rm lip}(f)|(x)&=
\begin{cases} \dis
\limsup_{y \to x}\frac{|f(y)-f(x)|}{d(y,x)} \quad & \text{if $x$ is not isolated},
\\ \dis
0\quad &\text{otherwise}.
\end{cases}
\end{align*}
According to \cite{AGS14},
 the Cheeger energy can be represented as follows:
\begin{align} \label{limit Ch}
{\sf Ch}(f)&=\frac{1}{2}\inf\Bigl\{ \liminf_{n \to \infty} \int |{\rm lip} f_n|^2 dm: f_n \in \mathrm{Lip}_{b}(X), \int_X|f_n -f|^2 dm \to 0 \Bigr\},
\\
f &\in W^{1,2}(X,d,m). \notag
\end{align}
Note that ${\rm Lip}_b(X) \subset L^p(X, m)$ for all $1 \le p \le \infty$ since we always assume $m(X)<\infty$ in \eqref{cond: fund}.
%
The Cheeger energy ${\sf Ch}: L^2(X,m) \to [0,+\infty]$ is a lower semi-continuous and a convex functional, but not necessarily a quadratic form.
The metric measure space $(X,d,m)$ is said to be {\it infinitesimal Hilbertian (IH)} if ${\sf Ch}$ is a quadratic form, i.e.,
\begin{align}
	2{\sf Ch}(u)+2{\sf Ch}(v)={\sf Ch}(u+v) + {\sf Ch}(u-v), \label{defn: Cheeger}
	\end{align}
	for any $u,v \in W^{1,2}(X,d,m).$
Under Condition (IH), the point-wise scalar product is defined in the following manner:
\begin{align} \label{def: PWS}
\langle \nabla f, \nabla g \rangle=\lim_{\e\to 0}\frac{|\nabla(f+\e g)|^2-|\nabla f|^2}{2\e}, \quad f,g \in W^{1,2}(X,d,m),
\end{align}
whereby  the limit is taken with respect to $L^1(m)$. If the Cheeger energy ${\sf Ch}$ is quadratic, the point-wise inner product is an $L^1(m)$-valued bilinear form (see \cite[Definition 4.12]{AGS14b}, and \cite[Theorem 2.7]{AST16}). Define the bilinear form ${\sf Ch}(f,g)=1/2\int_X \langle \nabla f, \nabla g \rangle dm$ for $f, g \in W^{1,2}(X,d,m)$.

\subsection{$L^p$-Wasserstein Space} \label{subsec: L2}
Let $(X_i,d_i)$ be a complete separable metric space for $i=1,2$.
For $\mu_i \in \mathcal P(X_i)$, a probability measure $\pi \in \mathcal P(X_1 \times X_2)$ is called {\it a coupling of $\mu_1$ and $\mu_2$} if 
$${{\rm proj}_1}_\#\pi = \mu_1 \quad \text{and} \quad {{\rm proj}_2}_\#\pi=\mu_2,$$ 
whereby ${\rm proj}_i$ denotes the projection ${\rm proj}_i: X_1 \times X_2 \to X_i$ as $(x_1,x_2) \mapsto x_i$. The set of all couplings of $\mu$ and $\nu$ is denoted by $\Pi(\mu, \nu)$.
Let $(X,d)$ be a complete separable metric space. Let $1 \le p \le \infty$. The set of all Borel probability measures $\mu$ with finite $p$-th moment is denoted by $\mathcal P_p(X)$:
$$\mathcal P_p(X)=\Bigl\{\mu \in \mathcal P(X): \int_Xd^p(x,\x)d\mu(x)<\infty \quad \text{for some (and thus any) $\x \in X.$}\Bigr\}.$$
The transportation distance $W_p$ is defined on $\mathcal P_p(X)$, which is called {\it $L^p$-Wasserstein distance}:
\begin{align} \label{eq: Was}
W_p(\mu,\nu)=\Bigl(\inf_{\pi \in \Pi(\mu,\nu)}\int_{X\times X}d^p(x,y)d\pi(x,y)\Bigr)^{1/p}.
\end{align}
It turns out that $(\mathcal P_p(X), W_p)$ is a complete separable geodesic metric space (e.g., \cite[Theorem 6.18]{V09}). A coupling $\pi \in \Pi(\mu,\nu)$ is called {\it an optimal coupling} if $\pi$ attains the infimum in the equality \eqref{eq: Was}. For any $\mu,\nu$, there always exists an optimal coupling $\pi$ of $\mu$ and $\nu$ (e.g., \cite[\S 4]{V09}).
\subsection{Pointed Measured Gromov Convergence} \label{subsec: D}
The pmG convergence has been introduced by Gigli-Mondino-Savar\'e \cite{GMS13}. A metric measure space $(X,d,m)$ with a fixed point $\x \in X$ is called {\it a pointed metric measure space} and denoted shortly as p.m.m.\ space. 
\begin{defn}\normalfont \cite[Def.\ 3.9]{GMS13} \label{prop: Dconv} 
A sequence of p.m.m.\ spaces $(X_n, d_n, m_n, \x_n)$ is said to {\it converge to $\mathcal (X_\infty, d_\infty, m_\infty, \x_\infty)$ in the pointed measured Gromov (pmG) sense} if there exist a complete separable metric space $(X,d)$ and an isometric embedding $\iota_n: X_n \to X\  (n \in \EN)$
satisfying the following:
\begin{align}\label{eq: VGC}
\iota_n(\x_n) \to \iota_\infty(\x_\infty) \in {X_\infty}, \quad \text{and}\quad \int_{X}f \ d({\iota_n}_\#m_n) \to \int_{X}f \ d({\iota_\infty}_\#m_\infty),
\end{align}
for any $f \in C_{bs}(X)$. If the underlying space $X_n$ is compact for any $n \in \EN$, then the fixed point $\x_n$ is dropped and just called {\it measured Gromov (mG)} convergence.  
\end{defn}
\begin{rem} \normalfont \label{rem: GHMG}
Several remarks for Definition \ref{prop: Dconv} are presented.
\begin{enumerate}
\item[(i)] The pmG convergence is a notion of a convergence for {\it isomorphism classes} of pointed metric measure spaces while the pointed measured Gromov-Hausdorff (pmGH) convergence is defined for each pointed metric measure spaces, not for equivalence classes (\cite[Remark 3.25]{GMS13}). The pmG convergence is weaker than the pmGH convergence (\cite[\replaced{Theorem}{Proposition} 3.30, Example 3.31]{GMS13}). If \deleted{${\rm supp}[m_\infty]=X_\infty$ and }$\{\mathcal X_n\}_{n \in \N}$ satisfies a uniform doubling condition, then pmG and pmGH coincide \cite[\replaced{Theorem}{Proposition} 3.33]{GMS13}.
\item[(ii)] The  pmG convergence is metrizable by a certain distance $p\mathbb G_W$ on the collection $\mathbb X$ of all isomorphism classes of p.m.m.\ spaces (\cite[Definition 3.13]{GMS13}). The space $(\mathbb X, p\mathbb G_W)$ \replaced{is}{becomes} a complete and separable metric space (\cite[Theorem 3.17]{GMS13}). 
\end{enumerate}
\end{rem}

\subsection{$L^p(m)$-normed Module} \label{subsec: L2T}
In this subsection, the notion of the $L^p$-normed module is recalled by following \cite[\S 1.2]{G16}. 
\begin{defn}\normalfont {\bf ($L^\infty(m)$-premodule)}  \cite[Def.\ 1.2.1]{G16}
{\it An $L^\infty(m)$-premodule} is a Banach space $(\mathcal M, \|\cdot\|_\mathcal M)$ equipped with a bilinear map
$L^\infty(m) \times \mathcal M \ni (f,v) \mapsto f\cdot v \in \mathcal M$ satisfying the following:
$$(fg)\cdot v = f \cdot(g\cdot v), \quad \1\cdot v = v, \quad \|f\cdot v\|_\mathcal M \le \|f\|_{L^\infty(m)}\|v\|_\mathcal M,$$
for any $v \in \mathcal M$ and $f,g \in L^\infty(m)$.
\end{defn}
\begin{defn}\normalfont {\bf ($L^\infty(m)$-module/Hilbert Module)}  \cite[Def.\ 1.2.1, 1.2.20]{G16}
{\it An $L^\infty(m)$-module} is an $L^\infty(m)$-premodule $\mathcal M$ satisfying the following two conditions:
\begin{description}
\item[(i)] ({\bf Locality}) For any $v \in \mathcal M$ and $A_n \in \mathscr B(X)$ ($n \in \N$), it holds that 
$$\1_{A_n}\cdot v=0, \quad \forall n \in \N \quad \text{implies} \quad \1_{\cup_{n \in \N}A_n}\cdot v=0.$$
\item[(ii)] ({\bf Gluing}) For any sequence $\{v_n\}_{n \in \N} \subset \mathcal M$ and $\{A_n\}_{n \in \N} \subset \mathscr B(X)$ so that 
$$\1_{A_i \cap A_j}\cdot v_i=\1_{A_i \cap A_j}\cdot v_j, \quad \forall i,j \in \N, \quad \text{and} \quad \limsup_{n \to \infty} \|\sum_{i=1}^n\1_{A_i}\cdot v_i\|_\mathcal M <\infty,$$
there exists $v \in \mathcal M$ so that 
$$\1_{A_i}\cdot v=\1_{A_i}\cdot v_i, \quad \forall i \in \N, \quad \text{and} \quad \|v\|_\mathcal M \le \liminf_{n \to \infty}\|\sum_{i=1}^n\1_{A_i}\cdot v_i\|_\mathcal M.$$
\end{description}
If, furthermore, $(\mathcal M, \|\cdot\|_{\mathcal M})$ is a Hilbert space, then $\mathcal M$ is called {\it a Hilbert module}.
\end{defn}

For two $L^\infty(m)$-modules $\mathcal M_1$ and $\mathcal M_2$, a map $T: \mathcal M_1 \to \mathcal M_2$ is called {\it a module morphism} provided that it is a bounded linear map from $\mathcal M_1$ to $\mathcal M_2$ as a map between Banach spaces and satisfies the following:
\begin{align} \label{prop: mod}
T(f \cdot v)=f\cdot T(v), \quad \forall v \in \mathcal M_1, f \in L^\infty(m).
\end{align}
The set of all module morphisms is denoted by $\text{\sc Hom}(\mathcal M_1, \mathcal M_2)$. The set $\text{\sc Hom}(\mathcal M_1, \mathcal M_2)$ has a canonical $L^\infty(m)$-module structure. 
\begin{defn}\normalfont {\bf (Dual Module)}  \cite[Def.\ 1.2.6]{G16} 
For an $L^\infty(m)$-module $\mathcal M$, {\it the dual module $\mathcal M^*$} is defined as $\text{\sc Hom}(\mathcal M, L^1(m))$.
\end{defn}
\begin{defn} \normalfont {\bf ($L^p(m)$-normed module and its dual)} \label{defn: Lpnormed} \cite[Def.\ 1.2.10]{G16}
 Let $p,q \in [1,\infty]$ and $1/p+1/q=1$. 
 \begin{itemize}
\item[(i)]An $L^p(m)$-normed module is an $L^\infty(m)$-module $\mathcal M$ endowed with a map $|\cdot|: \mathcal M \to L^p(m)$ with non-negative values so that 
$$\||v|\|_{p}=\|v\|_\mathcal M, \quad |f\cdot v|=|f||v|, \quad \text{$m$-a.e,}$$
for every $v \in \mathcal M$ and $f \in L^\infty(m)$. The map $|\cdot|$ is called {\it point-wise norm}.
\item[(ii)]  Let $\mathcal M$ be an $L^p(m)$-normed module. Then the dual $\mathcal M^*$ is  an $L^q(m)$-normed module with its point-wise norm $|\cdot|_*$ defined as follows (\cite[Prop.\  1.2.14 (i)]{G16}):
$$|L|_*:={\rm esssup}_{v \in \mathcal M, |v| \le 1 \ \text{$m$-a.e.}}|L(v)|.$$
\end{itemize}
\end{defn}

\subsection{Tangent Module} \label{subsec: Tan}
In this subsection, by following \cite[\S 2]{G16}, the notion of the tangent module $L^p(TX)$ on a metric measure space $(X,d,m)$ is recalled.  The tangent module $L^p(TX)$ is an $L^p(m)$-normed module in the sense of Definition \ref{defn: Lpnormed} and  a generalized notion of the space of $L^p$-sections of the tangent bundle on a smooth manifold. 

\begin{defn} \normalfont {\bf (Pre-cotangent module)} \cite[p.62]{G16}
The set {\rm Pcm} defined as follows is called {\it Pre-cotangent module}:
\begin{align*}
{\rm Pcm}=\Bigl\{\{(f_i, A_i)\}_{i \in \N}: &\{A_i\}_{i \in \N} \subset \mathscr B(X)\ \text{is a partition of X}
\\
&f_i \in S^2(X), \forall i \in \N,\ and \sum_{i \in \N}\int_{A_i}|\nabla f|^2dm<\infty \Bigr\}.
\end{align*}
An equivalence relation between two elements in Pcm $\{(f_i, A_i)\}_{i \in \N} \sim \{(g_j, B_j)\}_{j \in \N}$ is defined as follows:
$$|\nabla(f_i-g_j)|=0, \quad m\text{-a.e. on}\ A_i \cap B_j, \quad \forall i,j \in \N.$$
\end{defn}

The vector space structure can be endowed with the quotient space Pcm$/\sim$ by defining the sum and the scalar multiplication as follows:
$$[(f_i, A_i)_i]+[(g_j, B_j)_j]=[(f_i+g_j, A_i\cap B_j)_{i,j}], \quad \lambda[(f_i, A_i)_i]=[(\lambda f_i, A_i)_i] \quad (\forall \lambda \in \R).$$

The product operation $\cdot: {\rm Sf}(m) \times {\rm Pcm}/\sim \to {\rm Pcm}/\sim$ can be defined in the following manner: Let ${\rm Sf}(m) \subset L^\infty(m)$ denote the set of all simple functions $f$, which means that $f$ attains only a finite set of values. Given $[(f_i, A_i)_i] \in {\rm Pmc}/\sim$ and $h=\sum_{j}a_j\1_{B_j} \in {\rm Sf}(m)$ with $\{B_j\}_{j \in \N}$ being a partition of $X$, the product $h \cdot [(f_i, A_i)_i]$ is defined as follows:
$$h \cdot [(f_i, A_i)_i]=[(a_jf_i, A_i\cap B_j)_{i,j}].$$

The point-wise norm $|\cdot|_*$ is now recalled (note that the notation $|\cdot|_*$ is used as a point-wise norm for the sake of consistency with the definition of tangent modules given later): Define $|\cdot|_*$ on ${\rm Pcm}/\sim \to L^2(X,m)$ as
$$\bigl| [(f_i, A_i)_i]\bigr|_*=|\nabla f_i|, \quad \text{$m$-a.e. on $A_i$ for all $i \in \N.$}$$
Then the map $\|\cdot\|_{L^2(T^*X)}: {\rm Pcm}/\sim \to [0,\infty)$ is defined as follows:
$$\|[(f_i, A_i)_i]\|_{L^2(T^*X)}^2=\int\bigl| [(f_i, A_i)_i]\bigr|^2_*dm=\sum_{i \in \N}\int_{A_i}|\nabla f_i|^2dm.$$
It can be seen that $\|\cdot\|_{L^2(T^*X)}$ is a norm on Pcm$/\sim.$
\begin{defn}\normalfont {\bf (Cotangent Module)} \label{defn: cotan} \cite[Def.\ 2.2.1]{G16}
The cotangent module $(L^2(T^*X), \|\cdot\|_{L^2(T^*X)})$ is defined as the completion of $({\rm Pcm}/\sim, \|\cdot\|_{L^2(T^*X)})$.
\end{defn}
The cotangent module $(L^2(T^*X), \|\cdot\|_{L^2(T^*X)})$ is an $L^2$-normed module with the product $\cdot$ (which can be extended to the map $\cdot: L^\infty(m) \times L^2(T^*X) \to L^2(T^*X)$), and the point-wise norm $|\cdot|_*$ (see \cite[\S 2.2]{G16} for more details).
\begin{defn} \normalfont {\bf (Tangent Module)} \label{defn: tan} \cite[Def.\ 2.3.1]{G16}
The tangent module $(L^2(TX), \|\cdot\|_{L^2(TX)})$ is defined as the dual module of $(L^2(T^*X), \|\cdot\|_{L^2(T^*X)})$.
The point-wise norm associated with the dual of $|\cdot|_*$ is written as $|\cdot|$.
\end{defn}
Under Condition (IH), the tangent module $(L^2(TX), \|\cdot\|_{L^2(TX)})$ is a Hilbert module and the point-wise norm $|\cdot|$ satisfies the parallelogram identity. Hence, the point-wise inner product can be defined by $\langle \cdot, \cdot \rangle$. The notion of $L^p(TX)$ can be also defined for $1 \le p <\infty$. See \cite[Section 1.3, 3.2]{G16}. 

The notions of differential and gradient for a function in Sobolev class are recalled.
\begin{defn}  \normalfont {\bf (Differential)} \label{defn: dif} \cite[Def.\ 2.2.2]{G16}
Let $f \in S^2(X)$. The differential $df \in L^2(T^*X)$ is defined as follows:
$$df=[(f,X)] \in {\rm Pcm}/\sim \subset L^2(T^*X).$$
Here $[(f,X)]$ means $[(f_i,A_i)_{i \in \N}]$ for $f_0=f$, $A_0=X$ and $A_i=\emptyset$ for $i \ge 1$.
\end{defn}
By definition, $|df|_*=|\nabla f|$. The notion of gradient of a Sobolev function is defined through duality with the notion of the differential. 
\begin{defn} \normalfont {\bf (Gradient)} \label{defn: grad} \cite[Def.\ 2.3.4]{G16}
Let $f \in S^2(X)$. A vector field $X \in L^2(TX)$ is said to be {\it a gradient of $f$} if 
$$df(X)=|X|^2=|df|^2_*.$$
The set of all gradients of $f$ is denoted by ${\rm Grad}(f)$.
\end{defn}
Under Condition (IH), the set ${\rm Grad}(f)$ has a unique element, which is denoted by $\nabla f$. In this case, the gradient $\nabla f$ satisfies the following linearity (\cite[Proposition 2.3.17]{G16}):
$$\nabla(f+g)=\nabla f+\nabla g, \quad m\text{-a.e.}, \quad f,g \in S^2(X).$$
Assume (IH) and let $\langle \cdot, \cdot \rangle: L^2(TX) \to L^1(X,m)$ be the point-wise inner product, which is induced by the structure of $L^2$-normed module with the point-wise norm $|\cdot|$ in $L^2(TX)$. Then, $\langle \nabla f, \nabla g \rangle$ can be identified in the $m$-a.e.\ sense with the same expression defined in \eqref{def: PWS} in Subsection \ref{subsec: Pre}.

\subsection{Derivation} \label{subsec: Der}
In this subsection, the notion of the derivation is briefly explained by following \cite[\S 2.3.1]{G16}. Let $(X,d,m)$ be a metric measure space. 
\begin{defn} {\bf (Derivation)} \normalfont \cite[Def.\ 2.3.2]{G16}
A linear map $\b:S^2(X,d,m) \to L^1(m)$ is said to be {\it a derivation}  if there exists $h \in L^2(m)$
 so that 
 $$|\b(f)| \le h|\nabla f|,\quad m\text{-a.e., for all $f \in S^2(X,d,m).$}$$
 The $m$-a.e.\ smallest function $h$ satisfying the above inequality is denoted by $|\b|$. 
 The space of all derivations with $|\b| \in L^p(X,m)$ is denoted by ${\rm Der}^p(X,d,m)$ for $1 \le p \le \infty$. 
 \end{defn}
 
Derivation operators satisfy the local property: for any $f,g \in S^2(X,d,m)$,
$$|\b(f-g)| \le h|\nabla (f-g)|=0, \quad m\text{-a.e. on }\ \{f=g\}.$$
By the local property, the chain rule holds:
$$\b(\phi(f))=(\phi' \circ f)\b(f), \quad \phi \in {\rm Lip}(\R), \quad m\text{-a.e.},$$
and the Leibniz rule also holds:
$$\b(fg)=\b(f)g+f\b(g), \quad m\text{-a.e.}$$

For any vector field $V \in L^2(TX)$, the map $V \circ d: S^2(X,d,m) \to L^1(X,m)$ is a derivation. Here $d$ has been defined in Definition \ref{defn: grad}.
Conversely,  by \cite[Theorem 2.3.3]{G16}, for a given derivation $\b$, there exists a unique vector field $V \in L^2(TX)$ so that the following diagram commutes:
$$
\begin{diagram}
\node{S^2(X,d,m)}
\arrow{e,t}{d}
\arrow{se,b}{\b}
\node{L^2(T^*X)}
\arrow{s,r}{V} \\
\node{}
\node{L^1(X,m)} 
\end{diagram}
$$
The vector field $V$ is denoted by $\b^*$.  

Now the notion of {\it divergence of derivations} is recalled. 
\begin{defn}{\bf (Divergence)} \normalfont \cite[Def.\ 2.3.11]{G16}
A derivation $\b \in {\rm Der}^2(X,d,m)$ is said to {\it have a divergence in $L^2(X, m)$} if there exists $g \in L^2(X,m)$ so that 
$$-\int_{X}\mathbf b(f)dm=\int_X fg dm, \quad \forall f \in W^{1,2}(X,d,m).$$  
Such a $g$ is uniquely determined if it exists, and denoted by ${\rm div}\mathbf b.$ The existence of such $g$ is not necessarily ensured for general $\b$, but when the notation ${\rm div}\b$ is used in this paper, the existence of such $g$ is implicitly assumed. Let ${\rm Div}^p(X,d,m)=\{\b \in {\rm Der}^p(X,d,m): {\rm div} \b \in L^p(X,d,m)\}$ for $1 \le p \le \infty$.
\end{defn}
By using the Leibniz rule, the following holds: 
\begin{align} \label{eq: DIV}
\int_X\mathbf b(f)\phi dm=-\int_X \mathbf b(\phi)fdm-\int_{X}f\phi {\rm div}\mathbf b dm,\quad \forall f, \phi \in {\rm Lip}_{b}(X).
\end{align}

\subsection{Dirichlet Forms} \label{subsec: NSD}
In this subsection, basic notions of Dirichlet forms are presented.  See \cite{MR92} and \cite{O13} for more details. Let $(X,d,m)$ be a metric measure space. 
Let $\mathcal F \subset L^2(X,m)$ be a dense linear subspace and $\mathcal E$ be a bilinear form on $\mathcal F$. Let $\E_{\alpha}(f,g)=\E(f,g)+\alpha(f,g)$ and $\E_{\alpha}(f)=\E_\alpha(f,f)=\E(f,f)+\alpha\|f\|_2^2$ for $\alpha \in [0,\infty)$. {\it The symmetric part} is defined by ${\E}^s(f,g)=(1/2)(\E(f,g)+\E(g,f))$ and {\it the anti-symmetric part} by ${\E}^a(f,g)=(1/2)(\E(f,g)-\E(g,f))$. The bilinear form $(\E, \F)$ is said to be {\it a coercive closed form} if there is a constant $\lambda \ge 0$ so that
\begin{description}
\item[($\mathcal E.1$)] {$\mathcal E_{\lambda}$ is non-negatively definite}: $\E_{\lambda}(f) \ge 0$ for all $f \in \mathcal F$;
\item[($\mathcal E.2$)] {$\mathcal E$ satisfies the weak sector condition}: there exists a constant $C \ge 1$ so that 
$$|\E_{\lambda+1}(f,g)| \le C \E_{\lambda+1}(f)\E_{\lambda+1}(g), \quad \forall f, g \in \mathcal F;$$
\item[($\mathcal E.3$)] $\F$ is a Hilbert space with the norm $({\E}^s_{\lambda+1})^{1/2}.$
\end{description}
Let $D \subset L^2(X,m)$ be a dense linear subspace. A bilinear form $(\E, D)$ satisfying $(\E.1)$ and $(\E.2)$ is {\it closable} if, for any $f_n \in  D$ with $\lim_{n,m \to \infty}\E(f_n-f_m)=0$ and $\lim_{n \to \infty}\|f_n\|_{2}=0$, it holds that $\lim_{n \to \infty}\E(f_n)=0.$
The form $(\E, \F)$ is said to be {\it symmetric} if $\E(f,g)=\E(g,f)$ for all $f,g\in \mathcal F$. {\it The dual form} ${\E^*}$ is defined to be ${\E^*}(f,g)=\E(g,f)$ for $f, g \in \mathcal F$. 
If $(\E, \F)$ is a coercive closed form, then the corresponding semigroups $\{T_t\}$ and $\{T^*_t\}$ exist on $L^2(X,m)$ so that $(T_tf, g)=(g, T^*_tf)$ for any $t \ge 0$ and $f, g \in L^2(X,m)$. Furthermore, the corresponding resolvents $G_\alpha$ and $G^*_\alpha$ with $\alpha \ge \lambda$ are defined as $G_\alpha f=\int_0^\infty e^{-\alpha t}T_tf dt$ and $G^*_\alpha f=\int_0^\infty e^{-\alpha t}T^*_t f dt$, and it holds that
$$\E_\alpha (G_\alpha f, g)=(f,g)=\E_\alpha(g, G^*_\alpha f), \quad \forall f \in L^2(X,m),\ g \in \mathcal F,\ \alpha \ge \lambda.$$
Concerning the Markovian property, the following statements are equivalent (e.g., \cite[Theorem 1.1.5.]{O13}):
\begin{description}
\item[{\bf (${\E}.4$)}] for any $f \in \F$ and non-negative constant $a$, the following holds:
$$\E(f\wedge a, f-f\wedge a)\ge 0; $$ 
\item[(M)] $\{T_t\}$ is Markovian: if $f \in L^2(X,m)$ satisfies $0 \le f \le 1$ $m$-a.e., then $0 \le T_tf \le 1$ $m$-a.e.
\end{description}
Note that the property (M) is called {\it sub-Markovian} in \cite[Theorem 1.1.5.]{O13}, but we call it {Markovian} in this paper. 
The bilinear form $(\E, \F)$ is said to be a {\it Dirichlet form} if ($\mathcal E.1$)--($\mathcal E.4)$ hold. The form $(\E, \F)$ is said to satisfy {\it the dual Markov property} if 
\begin{description}
\item[{\bf ($\hat{\E}.4$)}] for any $f \in \F$ and non-negative constant $a$, the following hold
$$\E(f-f\wedge a, f\wedge a)\ge 0.$$ 
\end{description}

Now the regularity/quasi-regularity for Dirichlet forms is presented, which is a sufficient condition for the existence of Hunt processes/$m$-tight special standard processes associated with Dirichlet forms (see \cite[Theorem 3.5 in Chapter IV]{MR92}). An increasing sequence $\{E_n\}_{n \in \N}$ of closed subsets of $X$ is said to be an {\it$\E$-nest} if 
$$\bigcup_{n \in \N }\mathcal F_{E_n} \text{ is dense in $\mathcal F$ with respect to $({\E}^s_{\lambda+1})^{1/2}$},$$ 
whereby $\mathcal F_{A}=\{u \in \mathcal F: u=0\ m\text{-a.e.\ on $A^c$}\}$.
 A subset $N \subset X$ is said to be {\it $\E$-exceptional} if 
 $$N \subset \bigcap_{n \in \N}E_n^c \text{ for some $\E$-nest $\{E_n\}_{n \in \N}$.}$$
A property of points in $X$ holds {\it $\E$-quasi-everywhere ($\E$-q.e.)} if the property
 holds outside some $\E$-exceptional set. A function $f$ $\E$-q.e.\ defined on $X$ is called {\it $\E$-quasi-continuous} if there exists an $\E$-nest $\{E_n\}_{n \in \N}$ so that $f \in C(\{E_n\})$, whereby 
 $$C(\{E_n\})=\{f: A \to \R: \bigcup_{n \in \N} E_n\ \subset A \subset X,\ f|_{E_n}\text{ is continuous for all}\ n \in \N\}.$$
The Dirichlet form $(\E, \F)$ on $L^2(X,m)$ is said to be  {\it quasi-regular} if the following three conditions hold:
\begin{enumerate}
\item[(i)] there exists an $\E$-nest $\{E_n\}_{n \in \N}$ consisting of compact sets;
\item[(ii)] there exists an $({\E}^s_{\lambda+1})^{1/2}$-dense subset of $\F$ whose elements have $\E$-quasi-continuous $m$-versions;
\item[(iii)] there exist $u_n \in \F$ for $n \in \N$ having $\E$-quasi-continuous $m$-versions $\tilde{u}_n$ and an $\E$-exceptional set $N \subset X$ so that $\{\tilde{u}_n\}_{n \in \N}$ separates points of $X \setminus N$.
\end{enumerate}
Assume that $(X,d)$ is locally compact. The Dirichlet form $(\E, \F)$ is said to be {\it regular} with a core $\mathcal C_1$ if $\mathcal C_1 \subset C_0(X)\cap \mathcal F$ is dense both in $C_0(X)$ with the uniform norm $\|\cdot\|_\infty$ and in $\mathcal F$ with $({\E}^s_{\lambda+1})^{1/2}$ respectively. The form $(\E, \F)$ is quasi-regular if it is regular (\cite[Chapter IV Section 4a)]{MR92}.

Let $\{T_{t}\}$ be the semigroup corresponding to $(\E, \F)$. If the Dirichlet form $(\E, \F)$ is quasi-regular, then there exists an $m$-tight special standard process $(\Omega, \mathcal M, \mathcal M_t, S_t, \theta_t, \mathbb P_x)$ (\cite[Definition 1.5, Definition 1.13, Theorem 3.5 in Chapter IV]{MR92}) so that, for all $t \ge 0$ and $f \in \mathcal B_b(X) \cap L^2(X,m)$,
$$T_t f(x)=\mathbb E_x(f(S_t)), \quad \text{$\E$-q.e.\ $x$}.$$
Here $\mathbb E_x(f(S_t))=\int_\Omega f(S_t(\omega))\mathbb P_x(d\omega)$.

Assume that $(\E, \F)$ is quasi-regular. 
The form $(\E, \F)$ is said to be {\it conservative} if $T_t \1 =\1$ $m$-a.e.\ for any $t \ge 0$.
The form $(\E, \F)$ is said to be {\it local} if $\E(f,g)=0$ whenever $f, g \in \mathcal F$ with ${\rm supp}[f] \cap {\rm supp}[g]=\emptyset$. Note that the compactness of the supports of $f$ and $g$ are not necessary (see \cite[Proposition 1.2]{MR92}).
The form $(\E, \F)$ is said to be {\it strongly local} if, for any $f,g \in \mathcal F$, the following holds: if $g$ is constant on a neighborhood  of ${\rm supp}[f]$, then  $\E(f,g)=0.$

\subsection{RCD Spaces} \label{subsec: RCDLV}
In this subsection, RCD$(K,\infty)$/RCD$(K,N)$ spaces are briefly explained. Let $(X,d,m)$ be a metric measure space. 
Recall that ${\sf Ch}$ denotes the Cheeger energy and the property of infinitesimal Hilbertianity (IH) was defined in \eqref{defn: Cheeger}.
Under Condition (IH), considering \eqref{limit Ch}, ${\sf Ch}$ is a strongly local symmetric Dirichlet form (\cite{AGS14, AGS14b}). The corresponding semigroup $\{H_t\}$ is called {\it heat semigroup}, and the infinitesimal generator is denoted by $\Delta$. 
The following condition is considered:
\begin{align} \label{text: REP}
&\text{Every function $f \in W^{1,2}(X,d,m)$ with $|\nabla f| \le 1$ $m$-a.e. admits a continuous}  \notag
\\
&\text{$1$-Lipschitz representative $\tilde{f}$.} 
\end{align}
The gradient estimates of the heat semigroup is the following: for every $f \in W^{1,2}(X,d,m)$ with $|\nabla f| \le 1$ $m$-a.e., and every $t>0$, it holds that 
\begin{align} \label{BEI}
H_t f \in {\rm Lip}_b(X), \quad |\nabla H_tf|^2 \le e^{-Kt}H_t(|\nabla f|^2), \quad m\text{-a.e.\ in $X$}.
\end{align}
The gradient estimate with dimensional upper bounds is the following:
\begin{align}\label{BEK}
H_t f \in {\rm Lip}_b(X), \quad |\nabla H_tf|^2 +\frac{2Kt^2}{N(e^{Kt}-1)}|\Delta H_t f|^2\le e^{-Kt}H_t(|\nabla f|^2), \quad m\text{-a.e.\ in $X$}.
\end{align}
According to \cite{AGS15, AGS14, EKS15, CM16}, RCD$(K,\infty)$/RCD$(K,N)$ conditions can be identified with  \eqref{BEI}/\eqref{BEK} under (IH) and \eqref{text: REP}. Namely, 
the metric measure space $(X,d,m)$ is said to be {\it an RCD$(K,\infty)$} ({\it resp.\ RCD$(K,N)$}) space if \eqref{BEI} (resp.\ \eqref{BEK}) holds under (IH) and \eqref{text: REP}.
\begin{rem} \normalfont
Note that RCD$(K,\infty)$/RCD$(K,N)$ conditions can be defined also in terms of $K$-convexity of the relative entropy in the $L^2$-Wasserstein space with Condition (IH) (\cite{AGMR15, AGS14, EKS15, CM16})). 
\end{rem}
The class of RCD spaces contains various  (finite- and infinite-dimensional) singular spaces such as Ricci limit spaces (Sturm \cite{Sturm06, Sturm06-2}, Lott--Villani \cite{LV09}), Alexandrov spaces (Petrunin, Zhang--Zhu \cite{Pet11, ZZ10}), warped products and cones (Ketterer \cite{Ket14, Ket14a}), quotient spaces (Galaz-Garc\'ia--Kell--Mondino--Sosa \cite{GKMS17}), stratified spaces (Bertrand--Ketterer--Mondello--Richard \cite{BKMR18}) and infinite-dimensional spaces such as Hilbert spaces with log-concave measures (Ambrosio--Savar\'e--Zambotti \cite{ASZ09}), which are related to various stochastic partial differential equations. 

Under the RCD$(K,\infty)$ condition, the following gradient estimates hold (see \cite[Theorem 6.2]{AGS14b} and \cite[Theorem 5.5.2]{BGL14}):
\begin{align}
&|\nabla H_t f| \le e^{-Kt/2}H_t(|\nabla f|), \ m\text{-a.e.,}\quad \forall f \in W^{1,2}(X,d,m), \label{ineq: GE}
\\
& |\nabla H_t f|^2 \le e^{-Kt}H_t(|\nabla f|^2), \ m\text{-a.e.,}\quad \forall f \in W^{1,2}(X,d,m). \label{ineq: GE2}
\end{align}
The stability of RCD$(K,\infty)$/RCD$(K,N)$ spaces under the pmG convergence is presented.
\begin{thm} [\cite{ AGMR15, AGS14b, EKS15, GMS13, Sturm06, Sturm06-2}] {\bf (Stability of the RCD$(K,\infty)$/RCD$(K,N)$)}
\\
Let $K \in \R$ and $1 < N <\infty$.
Let $(X_n, d_n, m_n, \x_n)$ be a p.m.m.\ space satisfying RCD$(K,\infty)$/RCD$(K,N)$ for $n \in \N$.
If $(X_n, d_n, m_n, \x_n)$ converges to $(X_\infty, d_\infty, m_\infty, \x_\infty)$ in the pmG sense, then the limit space $(X_\infty, d_\infty, m_\infty, \x_\infty)$ also satisfies the RCD$(K,\infty)$/RCD$(K,N)$ condition. 
\end{thm}
It is known that the following Poincar\'e inequality, the log-Sobolev inequality and the Sobolev inequality hold (see e.g., \cite{V09}):
\begin{thm} \label{thm: PSm}
Let $(X,d,m)$ be an RCD$(K,\infty)$ with $K>0$ (called the first case hereinafter), or with ${\rm diam}(X)<D$ (called the second case hereinafter).
Then the following two statements hold:
\begin{enumerate} 
\item[(i)] there exists a constant $C_P>0$ depending only on $K$ in the first case, and on $K,D$ in the second case so that the $L^2$-Poincar\'e inequality holds:
\begin{align} \label{ineq: PI}
\int_X|f-m(f)|^2dm \le C_P \int_X |\nabla f|^2dm, \quad f \in W^{1,2}(X,d,m);
\end{align}
\item[(ii)] there exists a constant $C_{LS}>0$ depending only on $K$ in the first case, and on $K,D$ in the second case so that 
\begin{align} \label{ineq: LS}
\int_X f^2 \log f^2 dm \le \frac{C_{LS}}{2}\int_X |\nabla f|^2 dm, \quad  f \in W^{1,2}(X,d,m).
\end{align}
\end{enumerate}
Let $(X,d,m)$ be an RCD$(K,N)$ with $K>0$ and $N>2$. Then the following holds:
\begin{enumerate} 
\item[(iii)] there exists a constant $C_{S}>0$ depending only on $K,N$  so that 
\begin{align} \label{ineq: SL}
\|f\|^2_{2N/(N-2)}\le \|f\|_2^2+C_S \|\nabla f\|_2^2.
\end{align}
\end{enumerate}
\end{thm}
By the log-Sobolev inequality, the heat semigroup $\{H_t\}$ satisfies the hyper-contractive property (\cite[Theorem 6.1.14]{DS89}): 
\begin{align} \label{ineq: HC}
\|H_t\|_{p \to q} \le 1, \quad t>0,\ 1<p<q<\infty, \ \frac{q-1}{p-1} \le e^{4t/C_{LS}}.
\end{align}

\subsection{Test Functions and Test Vector Fields} \label{subsec: TFTV}
In this subsection, test objects on RCD spaces are presented: one is the class of test functions, the other is the class of test vector fields, both of which play fundamental roles as test objects. Readers are referred to \cite[Section 3]{G16} for more details. 

Let $(X,d,m)$ be an RCD$(K,\infty)$ space. The class ${\rm TestF}(X)$ of test functions is defined in the following manner:
$${\rm TestF}(X)=\Bigl\{f \in \mathcal D(\Delta) \cap L^\infty(X,m): |\nabla f| \in L^\infty(X,m) \quad \text{and} \quad \Delta f \in W^{1,2}(X,d,m)\Bigr\}.$$
By the Sobolev--to--Lipschitz property (\cite[Theorem 6.2]{AGS14b}), the following holds:
$$\text{any $f \in {\rm TestF}(X)$ has a Lipschitz representative $\tilde{f}: X \to \R$ with ${\rm Lip}(\tilde{f}) \le \| |\nabla {f}| \|_\infty$}. $$
By the $L^\infty-{\rm Lip}$ regularization property of the heat semigroup $\{H_t\}$ (\cite[Theorem 6.2]{AGS14b}), it can be seen that (see also \cite[\S 3.1]{G16})
$$f \in L^2(X,m) \cap L^\infty(X,m) \implies H_tf \in {\rm TestF}(X) \quad \forall t>0.$$
Therefore,  
$${\rm TestF}(X)\text{ is dense in}\ W^{1,2}(X,d,m).$$
Now the class of test vector fields is recalled. The class ${\rm TestV}(X) \subset L^2(TX)$ of test vector fields is defined in the following manner:
$${\rm TestV}(X)=\Bigl\{ \sum_{i=1}^n g_i \nabla f_i: n \in \N, f_i, g_i \in {\rm TestF}(X), i=1,2,..., n \Bigr\}.$$
By the approximating argument of the heat semigroup discussed in \cite[Section 3.2]{G16}, the following holds:
$${\rm TestV}(X)\text{ is dense in}\ L^2(TX).$$


\subsection{Brownian Motion and {Stochastic Analysis}} \label{subsec: B}
In this subsection, stochastic analysis concerning Brownian motions on RCD spaces is briefly explained. See \cite[Chapter 5]{FOT11} for a comprehensive account of stochastic analysis with respect to Dirichlet forms.
Let $(X,d,m)$ be an RCD$(K,\infty)$ space and ${\sf Ch}$ denote the Cheeger energy with the domain $W^{1,2}(X,d,m)$. For an open set $O \subset X$, define $\mathcal L_O=\{u \in W^{1,2}(X,d,m): u \ge 1\ m\text{-a.e.\ on $O$}\} $ and define the capacity of $O$ with respect to ${\sf Ch}$ as follows:
\begin{align*}
{\rm Cap}(O)=
\begin{cases}
\inf_{u \in \mathcal L_O} {\sf Ch}_1(u), & \text{if} \ \mathcal L_O \neq \emptyset,
\\
\infty, & \text{if}\ \mathcal L_O = \emptyset.
\end{cases}
\end{align*}
For any set $A \subset X$, define ${\rm Cap}(A)=\inf \{ {\rm Cap}(O): A \subset O, \ \text{$O$ is open}\}$. 
It is known that $({\sf Ch}, W^{1,2}(X,d,m))$ is a symmetric strongly local conservative quasi-regular Dirichlet form (\cite[Lemma 6.7]{AGS14b}). Let $(1/2)\Delta$, $\{H_t\}, \{R_\alpha\}$ denote the infinitesimal generator, the heat semigroup and the resolvent on $L^2$ corresponding to $({\sf Ch}, W^{1,2}(X,d,m))$ respectively.  The $L^2$-semigroup $\{H_t\}$ can be extended uniquely to $L^p$ for $p \in [1,\infty]$ and the extension is strongly continuous for any $1 \le p < \infty$ (e.g., \cite[Theorem 1.4.1]{D89}). Let $((1/2)\Delta_p, \mathcal D(\Delta_p))$ denote the infinitesimal generator associated with $\{H_t\}$ on $L^p(X,m)$ for $1 \le  p < \infty$. The subscript $p$ is used only if a particular $p$ is relevant. It holds that $\mathcal D(\Delta_p) \subset \mathcal D(\Delta_q)$ if $q \le p$ with $p,q \in [1,\infty)$.
By the quasi-regularity of $({\sf Ch}, W^{1,2}(X,d,m))$, there exists the corresponding conservative diffusion process $(\Omega, \mathcal M, \mathcal M_t, B_t, \theta_t, \mathbb P_x)$, called {\it Brownian motion}, so that (see Section \ref{subsec: NSD})
$$H_t f(x)=\mathbb E_x(f(B_t)), \quad \forall f \in \mathcal B_b(X), \ t\ge 0,\ \text{${\sf Ch}$-q.e.}\  x.$$ 

 Let $A_t(\omega)$ be an extended real valued function for $t \ge 0$ and $\omega \in \Omega$. The function $A$ is said to be {\it an additive functional (AF for brevity)} with respect to the Brownian motion $(\Omega, \mathcal M, \mathcal M_t, B_t, \theta_t, \mathbb P_x)$ if the following two conditions hold (\cite[\S 5.1]{FOT11}):
\begin{enumerate}
\item[(A.1)]
$A_t(\cdot)$ is $\mathcal M_t$-measurable for any $t \ge 0$;  
\item[(A.2)]
there exists a set $\Lambda \in \mathcal M_\infty$ and an exceptional set $N \subset X$ with respect to the capacity associated with the Cheeger energy ${\sf Ch}$ so that $\mathbb P_x (\Lambda)=1$ for all $x \in X \setminus N$, $\theta_t \Lambda \subset \Lambda$ for all $t >0$, and, moreover, for each $\omega \in \Lambda$, the path $A_\cdot (\omega)$ is right continuous and has the left limit on $[0,\infty)$, $A_0(\omega)=0$, $|A_t(\omega)| <\infty$ for $t>0$ and $A_{t+s}(\omega)=A_s(\omega)+A_t(\theta_s \omega)$ for $t,s \ge 0$. 
\end{enumerate}
The sets $\Lambda$ and $N$ referred to in (A.2) are called {\it a defining set} and {\it an exceptional set of the AF $A$} respectively. 
Two AFs $A^1$ and $A^2$ are said to be {\it equivalent} if, for each $t>0$, $\mathbb P_x(A^1_t=A_t^2)=1$ ${\sf Ch}$-q.e.\ $x$. These two AF $A^1$ and $A^2$ are identified with respect to this equivalence relation. An AF $A$ is said to be {\it continuous} if $A_t(\omega)$ is continuous in $t \in [0,\infty)$ for any $\omega$ in a defining set. A $[0,\infty]$-valued continuous AF is said to be {\it a positive continuous AF (PCAF)}. Let $\mathbf A_c^+$ denote the family of all equivalence classes of PCAFs. For any $f \in \mathcal B_b(X)$ with $f \ge 0$ and $A \in \mathbf A_c^+$, the PCAF $(fA)_t$ is defined as follows:
$$(fA)_t=\int_0^tf(B_s)dA_s.$$

A Borel measure $\mu$ is said to be {\it smooth} if $\mu$ charges no set of zero-capacity and there is an increasing sequence $\{F_k\}_k$ of  closed sets so that 
\begin{align*}
\mu(F_k)<\infty \quad (k=1,2,...),
\\
\lim_{k \to \infty}{\rm Cap}(K \setminus F_k)=0 \quad \text{for any compact set $K$}.
\end{align*}
Let $\mathbf S$ denote the family of all smooth measures. A PCAF $A \in \mathbf A_c^+$ and $\mu_A \in \mathbf S$ are said to be {\it in the Revuz correspondence} if 
\begin{align} \label{cor: Rev}
\lim_{t \downarrow 0}\frac{1}{t}\mathbb E_{h\cdot m}[(fA)_t]=\int_X h f\mu_A,
\end{align}
for any $\gamma$-excessive function $h$ with $\gamma \ge 0$ and $f \in \mathcal B(X)$ with $f \ge 0$. Here a function $f$ is {\it $\gamma$-excessive} if $f \ge 0$ and $e^{-\gamma t}H_tf \le f$ $m$-a.e., for any $t \ge 0$.
It is known that the family $\mathbf A_c^+$ and $\mathbf S$ are in one-to-one correspondence under the Revuz correspondence (see e.g., \cite[Theorem 5.1.4]{FOT11}). For $A \in \mathbf A_c^+$, the corresponding smooth measure is denoted by $\mu_A$, which is called {\it the Revuz measure of $A$}.

The energy $e(A)$ of the AF $A$ is defined in the following manner:
$$e(A)=\lim_{t \downarrow 0}\frac{1}{2t}\mathbb E_m[A^2_t],$$
if the limit exists in $[0,\infty]$. The mutual energy $e(A^1, A^2)$ is defined by the polarization: 
$$e(A^1, A^2)=\frac{1}{4}\bigl(e(A^1+A^2)-e(A^1-A^2)\bigr).$$ 
Let us define
$${\mathcal M}=\{M: \text{$M$ is a continuous AF s.t.\ for each $t \ge 0$, $\mathbb E_x[M_t^2]<\infty$, $\mathbb E_x[M_t]=0$\ ${\sf Ch}$-q.e.\ $x$.}\},$$
and 
$$\accentset{\circ}{\mathcal M}=\{M \in \mathcal M: e(M)<\infty\}.$$
Each $M \in \accentset{\circ}{\mathcal M}$ admits a unique PCAF $\langle M \rangle_t$, called {\it quadratic variation}, satisfying the following:
$$\mathbb E_x[\langle M \rangle_t]=\mathbb E_x[M_t^2], \quad \forall t >0, \ \text{${\sf Ch}$-q.e.}\ x.$$
Its Revuz measure $\mu_{\langle M \rangle}$ is called {\it the energy measure of $M$}. By \eqref{cor: Rev}, the energy measure satisfies the following equality:
\begin{align} \label{eq: ER}
e(M)=\lim_{t \downarrow 0}\frac{1}{2t}\mathbb E_m[M^2_t]=\lim_{t \downarrow 0}\frac{1}{2t}\mathbb E_m[\langle M \rangle_t]=\frac{1}{2}\mu_{\langle M \rangle}(X).
\end{align}
 The polarizations of the quadratic variation and the energy measure are denoted by $\langle \cdot, \cdot \rangle_t$ and $\mu_{\langle \cdot, \cdot \rangle}$ respectively.
 The space $\accentset{\circ}{\mathcal M}$ is a Hilbert space with the inner product $e(\cdot, \cdot)$. For $M \in \accentset{\circ}{\mathcal M}$ and $f \in L^2(X, \mu_{\langle M \rangle})$, there exists a unique element $f \bullet M \in \accentset{\circ}{\mathcal M}$ so that (see \cite[Theorem 5.6.1]{FOT11})
$$e(f \bullet M, L)=\frac{1}{2}\int_X f d\mu_{\langle M, L \rangle}, \quad \forall L \in \accentset{\circ}{\mathcal M}.$$
The map $f \mapsto f\bullet M$ is linear and continuous from $L^2(X, d\mu_{\langle M \rangle})$ into $(\accentset{\circ}{\mathcal M}, e)$.
The martingale $f \bullet M$ is called {\it stochastic integral}. 
The space $(\accentset{\circ}{\mathcal M}, e)$ is a $\mathcal B_b(X)$-module by the stochastic integral $(f\bullet M)_t$ for $f \in \mathcal B_b(X)$. 
For $f \in \mathcal B_b(X)$ and $M^1, M^2 \in \accentset{\circ}{\mathcal M}$, the following equality holds: 
\begin{align} \label{eq: energy}
d\mu_{\langle f\bullet M^1, M^2 \rangle}=f d\mu_{\langle M^1, M^2 \rangle}.
\end{align}
For $f \in W^{1,2}(X,d,m)$, there is a quasi-continuous modification $\tilde{f}$. By the Fukushima decomposition,
the following holds (See \cite[Theorem 5.2.2]{FOT11}):
$$\tilde{f}(B_t)-\tilde{f}(B_0)=M_t^{[f]}+N_t^{[f]}, \quad \text{for any $t \ge 0$}, \mathbb P_x\text{-a.e.}, \text{${\sf Ch}$-q.e.\ }x,$$
whereby $M_t^{[f]} \in \accentset{\circ}{\mathcal M}$ and $N^{[f]}$ is a continuous AF with $e(N^{[f]})=0$. Such $M^{[f]}$ and $N^{[f]}$ are uniquely determined as AFs. By \cite[Theorem 5.2.3]{FOT11}, it holds that $\mu_{\langle M^{[f]} \rangle}=|\nabla f|^2dm$ for any $f \in W^{1,2}(X,d,m)$. Combined with \eqref{eq: ER}, the quadratic variation of $M^{[f]}$ and $M^{[g]}$ for $f, g \in W^{1,2}(X,d,m)$, and its signed measure can be written as follows:
\begin{align} \label{eq: EM}
\langle M^{[f]}, M^{[g]} \rangle_t=\int_0^t \langle \nabla f, \nabla g \rangle (B_s)ds, \quad d\mu_{\langle M^{[f]}, M^{[g]} \rangle}=\langle \nabla f, \nabla g \rangle d\mu.
\end{align}



\subsection{Convergence of Non-symmetric Forms} \label{sec: CNSF}
In this subsection, the convergence of non-symmetric forms under varying metric measure spaces is recalled. In the case of fixed underlying spaces, the convergence of non-symmetric forms was introduced by \cite{H98}. In the case of varying underlying spaces, this notion has been generalized by \cite{T06, T10} and \cite{S18}. In this section, we follow \cite{S18}. The $L^2$-convergence of functions on varying metric measure spaces is recalled first.
\begin{defn} {\rm \cite[Definition 6.1]{GMS13}} \label{defn: Weak} \normalfont
Let $(X_n, d_n, m_n, {\x_n})$ be a sequence of p.m.m.\ spaces for $n \in \EN$. Assume that $(X_n,d_n,m_n, \x_n)$ converges to $(X_\infty, d_\infty, m_\infty, \x_\infty)$ in the pmG sense.  Let $(X,d)$ be a complete separable metric space and
$\iota_n: X_n \to X$ be an isometry for each $n \in \EN$ as in Definition \ref{prop: Dconv}. Identify $(X_n, d_n, m_n)$ with $(\iota_n(X_n), d, {\iota_n}_\#m_n)$ and omit $\iota_n$. 
\begin{description}
	\item[(i)] {\it $f_n \in L^2(X, m_n)$ converges weakly to $f_\infty \in L^2(X,m_\infty)$} if the following hold:
$$\sup_{n \in \N}\int|f_n	|^2\ dm_n < \infty \quad \text{and} \quad \int \phi f_n \ dm_n \to \int \phi f_\infty \ dm_\infty \quad \forall \phi \in C_{bs}(X).$$
 \item[(ii)] {\it $f_n \in L^2(X,m_n)$ converges strongly to $f_\infty \in L^2(X, m_\infty)$} if $f_n$ converges weakly to $f_\infty$ and the following holds:
$$\limsup_{n \to \infty}\int |f_n|^2 \ dm_n \le \int|f_\infty|^2\ dm_\infty.$$
\end{description}
\end{defn}

Let $(X_n, d_n, m_n, \x_n)$ be a sequence of p.m.m.\ spaces converging to $(X_\infty, d_\infty, m_\infty, \x_\infty)$ in the pmG sense. Let $(X,d)$ be a complete separable metric space and
$\iota_n: X_n \to X$ be an isometry for each $n \in \EN$ as in Definition \ref{prop: Dconv}. Let $(\E_n,\F_n)$ be a sequence of coercive forms with a common coercivity constant $\lambda \ge 0$. 
Let $\Phi_n(f)=\sup\{\E^n_{\lambda+1}(g, f): ({\E}^{n,s}_{\lambda+1}(g))^{1/2}=1\}$ for $f \in \mathcal F_n$,
 whereby ${\E}^{n,s}_{\lambda+1}$ denotes the symmetric part of $\E^n_{\lambda+1}$.
\begin{defn} \label{thm: Mosco of Ch} \normalfont (See also \cite[Definition 4.2]{S17})
{\it $(\E_n, \F_n)$ converges to $(\E_\infty, \F_\infty)$} if the following two conditions hold:
\begin{description}
	\item[(N1)] if a sequence $f_n \in L^2(X,m_n)$ converges weakly to $f_\infty \in L^2(X,m_\infty)$ with $\liminf_{n \to \infty}\Phi_n(f_n)<\infty$, then 
	$$f_\infty \in \mathcal F_\infty;$$
	\item[(N2)] for any sequence $f_n \in \F_n$ converging weakly in $L^2$ to $f_\infty \in \F_\infty$, and any $w_\infty \in \mathcal F_\infty$, there exists a sequence $w_n \in \mathcal F_n$ converging strongly in $L^2$ to $w_\infty \in \mathcal F_\infty$ so that 
	$$\lim_{n \to \infty}\mathcal E_n(f_n,w_n)=\mathcal E_\infty(f_\infty, w_\infty).$$
\end{description}
\end{defn}
In the case of symmetric forms, see, e.g., \cite{M94}, \cite{KS03} and \cite{GMS13}.

Verifying (N2) is not always easy, so another condition is introduced below:
\begin{description}
	\item[(N2$'$)] For any sequence $\{n_k\}_k \uparrow \infty$ and any sequence $f_k \in L^2(X; m_{n_k})$ weakly convergent in $L^2$ to $f_\infty \in \mathcal F_\infty$ with $\sup_{k \in \N}\Phi_{n_k}(f_k)<\infty$, there exists a dense subset $\mathcal C \subset \mathcal F_\infty$ for the topology with respect to ${\E}^{\infty,s}_{\lambda+1}$ so that every $w \in \mathcal C$ has a sequence $\{w_k\}$ \added{with $w_k \in \mathcal F_{n_k}$} converging to $w$ strongly \added{in $L^2$} with 
	$$\liminf_{k \to \infty}\E_{n_k}(w_k,f_k) \le \E_\infty(w,f_\infty).$$
\end{description}
Define (N$1_*$) by replacing $\Phi_n(f_n)$ with $({\E}^{n,s}_{\lambda+1}(f_n))^{1/2}$ in (N1), and (N$2_*'$) by replacing $\Phi_{n_k}(f_k)$ with $({\E}^{n_k,s}_{\lambda+1}(f_k))^{1/2}$ in (N2$'$).

Let $\{T^n_t\}$ and $\{G_\alpha^n\}$ be the $L^2$-contraction semigroup and resolvent associated with $\mathcal E_n$. 
\begin{description}
	\item[(R)] For any sequence ${f_n}$ converging to ${f_\infty}$ strongly in $L^2$, the resolvent $G^n_{{\alpha}} f_n$ converges to $G_{{\alpha}}^\infty f_\infty$ strongly in $L^2$ for any ${\alpha>0}$;
	\item[(S)] For any sequence ${f_n}$ converging to ${f_\infty}$ strongly in $L^2$, $T^n_t f_n$ converges to $T_t^\infty f_\infty$ strongly in $L^2$. The convergence is uniform on any compact time interval $[0,T].$ 
\end{description}

\begin{thm} \cite[Thm.\ 4.4]{S18}   \label{thm: EQUIVHINO}
The following statements hold:
\begin{enumerate}
\item[(i)] $(N1)\ \& \ (N2) \iff (N1)\ \& \ (N2') \iff (R) \iff (S)$;
\item[(ii)]  $(N1_*) \ \& \ (N2_*') \implies (R)$.
\end{enumerate}
\end{thm}

Several convergences of derivations are presented by following \cite{AST16}.
Let $\mathscr A\subset {\rm Lip}_b(X_\infty)$ denote the smallest algebra containing the following functions:
$$\min\{d(\cdot, x), k\},\quad k \in \Q, \ x \in D \subset X_\infty, \ \text{dense subset}.$$
The algebra $\mathscr A$ becomes a vector space over $\Q$. Let $\mathscr A_{bs}$ be a subalgebra consisting of bounded support functions. 
Let $\{H^\infty_t\}$ be the heat semigroup associated with Cheeger energy ${\sf Ch}_\infty$.
Let $H_{\Q_+}\mathscr A_{bs}=\{H^\infty_sf: f \in \mathscr A_{bs}, s \in \Q_+\} \subset {\rm Lip}_{b}(X)$.
\begin{defn} \normalfont $(${\rm \cite[Definition 4.3, 5.3]{AST16}}$)$ Let $(X_n,d_n,m_n,\x_n)$ converge to $(X_\infty, d_\infty, m_\infty, \x_\infty)$ in the pmG sense. Let $(X,d)$ be a complete separable metric space and
$\iota_n: X_n \to X$ be an isometry for each $n\in \EN$ as in Definition \ref{prop: Dconv}. Identify $(X_n, d_n, m_n)$ with $(\iota_n(X_n), d, {\iota_n}_\#m_n)$ and omit $\iota_n$. 
\begin{itemize}
\item[(i)] (Weak Convergence) We say that $\b_n \in {\rm Der}^2(X,d,m_n)$ converges weakly to $\b_\infty \in  {\rm Der}^2(X,d,m_\infty)$  in duality with $H_{\Q_+}\mathscr A_{bs}$ if, for all $f \in H_{\Q_+}\mathscr A_{bs}$, 
$$\int_X\b_n(f)h dm_n \to \int_X\b_\infty(f)h dm_\infty \quad \forall h \in C_{bs}(X).$$
\item[(ii)] (Strong Convergence) We say that $\b_n \in {\rm Der}^2(X,d,m_n)$ converges strongly to $\b_\infty \in  {\rm Der}^2(X,d,m_\infty)$ if, for all $f \in H_{\Q_+}\mathscr A_{bs}$, the function $\b_n(f)$ converges in measure to $\b_\infty(f)$, i.e., 
$$\int_X\Phi(\b_n(f))h dm_n \to \int_X\Phi(\b_\infty(f))h dm_\infty \quad \forall h \in C_{bs}(X), \quad \forall \Phi \in C_b(\R).$$
\item[(iii)] ($L^p$-strong Convergence) Let $p \in [1,\infty)$. We say that $\b_n \in {\rm Der}^p(X,d,m_n)$ converges $L^p$-strongly to $\b_\infty \in  {\rm Der}^p(X,d,m_\infty)$  if $\b_n$ converges strongly to $\b_\infty$ and,  for all  $f\in H_{\Q_+}\mathscr A_{bs}$, the following holds:
$$\limsup_{n \to \infty} \int_{X}|\b_n(f)|^pdm_n \le \int_{X}|\b_\infty(f)|^pdm_\infty.$$
\end{itemize}
\end{defn}

The $W^{1,2}$-convergence of functions on varying metric measure spaces is recalled. 
\begin{defn}$(${\rm\cite[Definition 5.2]{AH16}}$)$ \label{defn: WConv} \normalfont
Let $(X_n, d_n, m_n, \x_n)$ be a sequence of p.m.m.\ spaces. Assume that $(X_n,d_n,m_n, \x_n)$ converges to $(X_\infty, d_\infty, m_\infty, \x_\infty)$ in the pmG sense.  Let $(X,d)$ be a complete separable metric space and
$\iota_n: X_n \to X$ be an isometry for each $n \in \EN$ as in Definition \ref{prop: Dconv}. Identify $(X_n, d_n, m_n)$ with $(\iota_n(X_n), d, {\iota_n}_\#m_n)$ and omit $\iota_n$. 
\begin{description}
	\item[(i)] {\it $f_n \in W^{1,2}(X, m_n)$ converges to $f_\infty \in W^{1,2}(X,m_\infty)$ weakly in $W^{1,2}$} if $f_n \to f_\infty$ weakly in $L^2$ in the sense of Definition \ref{defn: Weak} and 
$\sup_{n \in \N}{\sf Ch}_n(f_n)< \infty$;
 \item[(ii)] {\it $f_n \in W^{1,2}(X,m_n)$ converges to $f_\infty \in W^{1,2}(X, m_\infty)$ strongly in $W^{1,2}$} if $f_n$ converges to $f_\infty$ strongly in $L^2$ in the sense of Definition \ref{defn: Weak} and $\lim_{n \to \infty}{\sf Ch}_n(f_n)={\sf Ch}_\infty(f_\infty)$.
\end{description}
\end{defn}


\section{Construction of Semigroups and Diffusion Processes} \label{sec: CSD}
In this section, the construction of diffusion processes on RCD spaces is discussed.
 Throughout this section, $(X,d,m)$ is an RCD$(K,\infty)$ space for $K \in \R$ and $(\Omega, \mathcal M, \mathcal M_t, B_t, \theta_t, \mathbb P_x)$ denotes the Brownian motion on $(X,d,m)$ (see Subsection \ref{subsec: B}).
\subsection{Construction of Diffusions by Girsanov's Transformation} \label{subsec: CD}
An isometry from the tangent module $L^2(TX)$ to $\accentset{\circ}{\mathcal M}$ is constructed in order to utilize the Girsanov transformation. For a given $u \in W^{1,2}(X,d,m)$, define a map $I$ as follows:
\begin{align} \label{map: I}
I(\nabla u)=M^{[u]},
\end{align}
whereby $M^{[u]} \in \accentset{\circ}{\mathcal M}$ denotes the martingale part of the Fukushima decomposition with respect to the Brownian motion $(B_t, \mathbb P_x)$ (see Subsection \ref{subsec: B}): 
$$u(B_t)-u(B_0)=M_t^{[u]}+N_t^{[u]}, \quad \mathbb P_x\text{-a.s.} \quad \text{${\sf Ch}$-q.e. $x \in X$}. $$
By the following proposition, the map $I$ can be extended to the isometry $I: L^2(TX) \to \accentset{\circ}{\mathcal M}$. 
\begin{prop} \label{prop: isom}
The map $I$ defined in \eqref{map: I} can be uniquely extended to the isometry $I: L^2(TX) \to \accentset{\circ}{\mathcal M}$. Furthermore,  the map $I$ is an $L^\infty$-module morphism in the sense of \eqref{prop: mod} and the following equality holds:
\begin{align} \label{eq: Iso}
\langle I(V), I(W) \rangle_t = \frac{1}{2}\int_0^t \langle V, W \rangle(B_s)ds, \quad \forall V, W \in L^2(TX),\ t \ge 0.
\end{align}
\end{prop}
\proof
Let $\mathcal V=\{\sum_{i=1}^nf_i \nabla u_i: f_i \in \mathcal B_b(X), u_i \in W^{1,2}(X,d,m), n \in \N\}$. For $V=\sum_{i=1}^nf_i \nabla u_i \in\mathcal V$, define the map $I$ as follows:
$$I(V)=\sum_{i=1}^nf_i \bullet M^{[u_i]}.$$
Let $W=\sum_{i=1}^ng_i \nabla v_i \in \mathcal V$ be another vector field. Then, by \eqref{eq: ER}, \eqref{eq: energy} and \eqref{eq: EM}, the following holds:
\begin{align*}
e(I(V), I(W))&=\frac{1}{2}\int_M\sum_{1 \le i,j \le n}f_ig_jd\mu_{\langle u_i, v_j \rangle}
=\frac{1}{2}\int_M\sum_{1 \le i,j \le n}f_ig_j \langle \nabla u_i, \nabla v_j \rangle dm
=\frac{1}{2}\langle V, W\rangle_{L^2(TX)}. 
\end{align*}
Therefore, $I$ is an isometry and well-defined as a map from $\mathcal V$ to $\CC{\mathcal M}$. 
By the same argument as in \cite[Lemma 5.6.3]{FOT11}, $\mathcal V$ generates a dense subset  in $\CC{\mathcal M}$. 
Since ${\rm Test}V(X) \subset \mathcal V \subset L^2(TM)$ is dense in $L^2(TM)$ (see Subsection \ref{subsec: TFTV}), the map $I$ can be extended to the isometry $I: L^2(TX) \to \accentset{\circ}{\mathcal M}$. The uniqueness and the property of  $L^\infty$-module morphism can be easily seen by the construction of $I$. 

Finally, the equality \eqref{eq: Iso} is checked. For any Borel set $A \subset X$, by \eqref{eq: ER}, \eqref{eq: energy} and \eqref{eq: EM}, the following holds: 
\begin{align*}
\frac{1}{2}\int_A \langle V, W \rangle dm=e(I(V), I(\1_{A}W))=e(I(V), \1_{A}\bullet I(W))=\mu_{\langle I(V), I(W)\rangle}(A).
\end{align*}
Therefore, $(1/2)\langle V, W \rangle dm=d\mu_{\langle I(V), I(W)\rangle}$, which implies \eqref{eq: Iso}. 
\qed

Let $\b$ be a derivation with $|\b| \in L^2(X,m)$ and $\b^* \in L^2(TX)$ denote the dual of $\b$ (see Subsection \ref{subsec: Der}).  Set $N=I(\b^*)$ and $M_t=\exp{\{N_t-(1/2)\langle N \rangle_t\}}$.
\begin{prop} \label{prop: mart}
$\{M_t\}$ is a continuous martingale under $\mathbb P_x$ for ${\sf Ch}$-q.e.\ $x \in M$.
If, moreover, $\exp\{\theta |\b|^2\} \in L^1(X,m)$ holds for some $\theta>0$, then the following statements hold:
\begin{itemize} 
\item[(i)]   for any $p \ge 1$, there exists $T>0$ so that $\{M_t\}$ is an $L^p$-martingale on $[0,T]$ under $\mathbb P_m$ and $\mathbb P_x$ for  ${\sf Ch}$-q.e.\ $x \in M$;
\item[(ii)] for any $p>0$, $r>1$ and $t \ge 0$, it holds that
$$\mathbb E_x[M_t^p] \le \mathbb E_x \Bigl[\exp\Bigl\{\frac{pr(pr-r+1)}{p(r-1)}\int_0^t|\b(B_s)|^2ds\Bigr\}\Bigr]^{1/r} \quad \text{ ${\sf Ch}$-q.e.\ $x \in X$}.$$ 
\end{itemize}
\end{prop}
\proof
The first part can be proved by the combination of the Markov property of  $(B_t, \mathbb P_x)$ and the super-martingale property of $\{N_t\}$ under $\mathbb P_x$  ${\sf Ch}$-q.e.\ $x$.  The second part can be shown by using \eqref{eq: Iso} and the H\"older inequality. See \cite[Proposition 2.2]{H98-2} for details. 
\qed

The construction of diffusion processes is presented by making use of the Girsanov transformation. Let $N_1$ be an exceptional Borel set so that $\{M_t\}$ is a continuous martingale under $\mathbb P_x$ for $x \in M \setminus N_1$. Redefine $M_t=\exp{\{N_t-(1/2)\langle N \rangle_t \cdot \1_{B_0 \in X \setminus N_1}\}}$. Set the following probability measure:
\begin{align} \label{GTR}
\mathbb Q_x|_{\mathcal M_t} = M_t\mathbb P_x|_{\mathcal M_t} \quad \forall t>0.
\end{align}
 Then $(\Omega, \mathcal M, \mathcal M_t, B_t, \theta_t, \mathbb Q_x)$ is a diffusion process on $X$. The choice of $N_1$ does not matter to the process in the following sense: Let $N_2 \supset N_1$ be an exceptional Borel set so that $M \setminus N_2$ is $(B_t, \mathbb P_x)$-invariant. Then $M \setminus N_2$ is also $(B_t, \mathbb Q_x)$-invariant. See \cite[Proposition 2.3]{H98-2}. 
Now define $T_tf$ as follows: for any $f \in \mathcal B_b(M)$ and $t>0$, 
$$T_t f(x)=\mathbb E^\mathbb Q_x[f(B_t)]=\int_{\Omega}f (B_t)d\mathbb Q_x=\mathbb E_x[f(B_t)M_t], \quad T_0f=f.$$
The family of the operators $\{T_t\}$ can be uniquely extended to an $L^\infty$-contraction Markovian semigroup $\{T_t\}$. In the next section, we extend $\{T_t\}$ to the $L^p$-space under Condition (B), for which the log-Sobolev inequality plays a key role. 
\subsection{$L^p$-semigroups} \label{subsec: L^p-es}
In this subsection, an $L^p$-estimate of the semigroup $\{T_t\}$ is obtained.
Throughout this subsection, $(X,d,m)$ is an RCD$(K,\infty)$ space with $K>0$ or with a bounded diameter. Under these assumptions, the log-Sobolev inequality \eqref{ineq: LS} holds with the log-Sobolev constant $C_{LS}$. Set $\alpha=C_{LS}$. Furthermore, assume that  
$$\exp\{\theta |\b|^2\} \in L^1(X,m) \quad \text{for some $\theta>\frac{\alpha}{2}$}.$$
Fix $p>1$ so that
 $$\theta > \frac{p^2}{2(p-1)^2}\alpha.$$
 We name this set of assumptions as (B).
\begin{prop} \label{prop: L2}
Suppose (B). 
Then the semigroup $\{T_t\}$ can be extended to a strongly continuous semigroup on $L^p(X,m)$. Furthermore, $\{T_t\}$ is hypercontractive in the following sense: for $t>0$, $p<q<\infty$, $\alpha'>\alpha$ with $\frac{\alpha}{2}\frac{p}{(p-1)^2}\frac{\alpha'^2}{(\alpha'-\alpha)^2}<\theta$ and $\frac{q-1}{p-1} \le e^{4t/\alpha'}$, it holds that
\begin{align} \label{eq: OPNORM}
\|T_t\|_{p \to q}\le \biggl\|\exp\biggl\{ \frac{\alpha}{2}\biggl(\frac{p\alpha'}{(p-1)(\alpha'-\alpha)}\biggr)^2|\b|^2\biggr)\biggr\} \biggr\|_1^{\frac{2}{p-1}\frac{(q-1)^2}{q^2}(\frac{1}{\alpha} - \frac{1}{\alpha'})t}.
\end{align}
In particular, 
\begin{align} 
\|T_t\|_{p\to p} \le  \biggl\|\exp\biggl\{ \frac{\alpha}{2}\biggl(\frac{p}{(p-1)}\biggr)^2|\b|^2\biggr)\biggr\} \biggr\|_1^{\frac{2(p-1)t}{p^2\alpha}}.
\end{align}
\end{prop}
\proof
The proof is similar to \cite[Proposition 3.1]{H98-2}. The hyper-contractivity of the heat semigroup $\{H_t\}$ is used, which follows from the log-Sobolev inequality \eqref{ineq: LS}. 
\qed

Let $(L_p, \mathcal D(L_p))$ be the infinitesimal generator corresponding to $\{T_t\}$ on $L^p(X,m)$.
Let $\overline{L}=\frac{1}{2}\Delta+\b$ with the domain 
$$\mathcal D(\overline{L})=\{f: f \in \mathcal D(\Delta_1)\cap S^2(X,d,m),\ |\nabla f| \in L^1(X,m)\},$$ whereby $\Delta_1$ denotes the infinitesimal generator of the heat semigroup $\{H_t\}$ on $L^1(X,m)$. A criterion for functions belonging to the domain $\mathcal D(L_p)$ of $L_p$ is shown in the following proposition.
\begin{prop} \label{prop: domain} If $f \in \Bigl( \cup_{r>1}\mathcal D(\Delta_r)\Bigr) \cap W^{1,2}(X,d,m) \cap L^p(X,m)$ and $\overline{L}f \in L^p(X,m)$, then $f \in \mathcal D(L_p)$ and $L_pf=\overline{L}f.$
\end{prop}
\proof
The proof can be obtained by the same way as \cite[Proposition 3.2]{H98-2}.  
\qed
\\
As a corollary of Proposition \ref{prop: domain}, it can be easily checked that $\mathcal D(\Delta_2) = \mathcal D(L_2)$ if $|\b| \in L^\infty(X,m)$.
\subsection{Dirichlet Forms} \label{subsec: DFC}
In this subsection, several important properties of Dirichlet forms are proved. 
Let $\E$ be the following bilinear form:
\begin{align*} 
\E(f,g)=\frac{1}{2}\int_X \langle \nabla f, \nabla g \rangle dm-\int_X \b(f)g dm, \quad f,g \in {\rm Lip}_b(X).
\end{align*}
The following three types of assumptions are considered: Let ${\rm diam}(X)=\sup\{d(x,y): x,y \in X\}.$
\begin{itemize}
\item[(A)] Let $(X, d, m)$ be an RCD$(K,\infty)$ space with $m(X)<\infty$. Assume $|\b| \in L^\infty(X,m)$, ${\rm div}\b \ge 0$.
\item[(B)] Let $(X, d, m)$ be an RCD$(K,\infty)$ space with $K>0$ or ${\rm diam}(X)<D<\infty$.
Assume $\exp\{{\theta|\b|^2}\} \in L^1(X,m)$ for some $\theta > \frac{C_{LS}}{2}$, whereby $C_{LS}$ denotes the constant in the log-Sobolev inequality \eqref{ineq: LS}. Let $p>1$ with $\theta>C_{LS}\frac{p^2}{2(p-1)^2}$.
 
\item[(C)]  Let $(X, d, m)$ be an RCD$(K,N)$ space with $K>0$ and $2<N<\infty$. Assume that $|\b| \in L^N(X,m)$ and ${\rm div} \b \ge -c$ for some constant $c \ge 0$.
\end{itemize}
Note that $m(X)<\infty$ holds under each of the three assumptions. 
\begin{prop} \label{prop: form}
Under (A), (B) with $\theta >2C_{LS}$, or (C), the bilinear form $(\E, {\rm Lip}_b(X))$ is closable and its closure $(\E, \F)$ is a quasi-regular strongly local conservative coercive Markovian form. Furthermore, $\mathcal F=W^{1,2}(X,d,m)$ and 
\begin{align} \label{eq: DF}
\E(f,g)=\frac{1}{2}\int_X \langle \nabla f, \nabla g \rangle dm-\int_X \b(f)g dm, \quad f,g \in \mathcal F=W^{1,2}(X,d,m),\ g\ \text{is bounded}.
\end{align}
Moreover, the following statements hold:
\begin{itemize}
\item[(i)] Under (A), the form $(\E, \F)$ is positivity-preserving and has the dual-Markov property. The coercivity constant is $\lambda=0$. 
\item[(ii)] Under (B) with $\theta>2C_{LS}$, the infinitesimal generator $(L, \mathcal D(L)))$ associated with $(\E, \F)$ coincides with $(L_2, \mathcal D(L_2))$ defined in Subsection \ref{subsec: L^p-es}. 
\item[(iii)] Under (C), the form $(\E, \F)$ is regular. 
\end{itemize}
\end{prop}
\proof
The conservativeness is obvious since the total measure $m(X)$ is finite.
Under (A), all the properties stated in the proposition have been proved in \cite[(b) Proposition 3.2]{S18}.  It suffices to give a proof only for the case of (B) and (C). 
Note that $({\sf Ch}, W^{1,2}(X,d,m))$ is a quasi-regular symmetric Dirichlet form in the case RCD$(K,\infty)$, and a regular symmetric Dirichlet form in the case of RCD$(K,N)$ (see \cite[Lemma 6.7]{AGS14b}). Therefore, 
concerning the coercivity, the (quasi-)regularity, $\mathcal F=W^{1,2}(X,d,m)$ and \eqref{eq: DF} for $(\E, \F)$, it suffices to show that there are constants $\lambda \ge 0$, $c, c'>0$ and $C \ge 1$ so that 
\begin{align} \label{ineq: 1}
c \|f\|^2_{1,2} \le \E_\lambda (f) \le c' \|f\|^2_{1,2}, \quad \forall f \in {\rm Lip}_{b}(X),
\end{align}
and the weak sector condition:
\begin{align} \label{ineq: 2}
\E_{\lambda+1}(f,g) \le C \E_{\lambda +1}^{1/2}(f) \E_{\lambda+1}^{1/2}(g), \quad \forall f, g\in {\rm Lip}_{b}(X).
\end{align}
In order to show \eqref{ineq: 1} and \eqref{ineq: 2}, we only have to show that there are $\e\in (0,1)$ and $a, b \ge 0$ so that 
\begin{align} \label{ineq: DD}
\int_X \b(f)f dm \le  \e {\sf Ch}(f)+a\|f\|_2^2,
\end{align}
and 
\begin{align} \label{ineq: DD1}
\biggl|\int_X \b(f)g dm \biggr| \le b {\sf Ch}^{1/2}_1(f) {\sf Ch}^{1/2}_1(g),
\end{align}
for any $f, g \in {\rm Lip}_b(X)$. (See e.g., \cite[(2.12a), (2.12b)]{BRZ00}, which imply the equivalence of the norms \eqref{ineq: 1} and the weak sector condition \eqref{ineq: 2}. Note that the plus-minus sign in front of $\b$ in \eqref{eq: DF} is opposite to  the Dirichlet form in \cite[p. 96]{BRZ00}, so the plus-minus sign in \eqref{ineq: DD} is also opposite to \cite[(2.12b)]{BRZ00}.) In this case, the constant $\lambda$ can be taken as follows:
\begin{align} \label{const: lambda}
\lambda=a+1 - \e.
\end{align}

Assume (B) with $\theta>2C_{LS}$. We first check \eqref{ineq: DD1}. By making use of the Hausdorff--Young inequality $st \le e^s+t\log t^t$ for $s \in \R$ and $t > 0$ and the log-Sobolev inequality \eqref{ineq: LS}, it holds that, for $f, g \in {\rm Lip}_b(X)$,
\begin{align*} 
\Bigl| \int_X \b(f)g dm \Bigr|  \le \int_X |\b||\nabla f||g| dm &\le \|\nabla f\|_2\||\b| g\|_2 \notag
\\
& \le \|\nabla f\|_2\|g\|_2 \Bigl\||\b| \frac{g}{\|g\|_2}\Bigr\|_2 \notag
\\
&\le  \|\nabla f\|_2 \|g\|_2 \biggl\{ \int_X \theta^{-1} \biggl( e^{\theta |\b|^2}+\frac{g^2}{\|g\|_2^2} \log{\frac{g^2}{\|g\|_2^2}} \biggr)dm \biggr\}^{1/2} \notag
\\
&\le  \theta^{-1/2} \|\nabla f\|_2 \|g\|_2 \biggl\{ \Bigl\|e^{\theta |\b|^2}\Bigr\|_1+\frac{C_{LS}}{2}\frac{\|\nabla g\|_2^2}{\|g\|_2^2} \biggr\}^{1/2} \notag
\\
&\le \theta^{-1/2}\|\nabla f\|_2 \biggl(  \Bigl\|e^{\theta |\b|^2}\Bigr\|_1^{1/2}\|g\|_2 +\sqrt{\frac{C_{LS}}{2}}\|\nabla g\|_2 \biggr)
\\
&\le \theta^{-1/2}\Bigl\|e^{\theta |\b|^2}\Bigr\|_1^{1/2} {\sf Ch}^{1/2}_1(f) {\sf Ch}^{1/2}_1(g) +\theta^{-1/2}\sqrt{\frac{C_{LS}}{2}}{\sf Ch}^{1/2}_1(f){\sf Ch}^{1/2}_1(g) \notag
\\
&= b{\sf Ch}^{1/2}_1(f) {\sf Ch}^{1/2}_1(g). \notag
\end{align*}
Thus \eqref{ineq: DD1} has been checked. Now we check \eqref{ineq: DD}.
Take a positive constant $\alpha$ so that 
\begin{align} \label{alpha}
0<\alpha<1-\sqrt{2C_{LS}}\theta^{-1/2}.
\end{align}
Then we have that 
\begin{align*}
\int_X \b(f)f dm &\le \theta^{-1/2}\sqrt{\frac{C_{LS}}{2}}\|\nabla f\|_2^2+\theta^{-1/2}\Bigl\|e^{\theta |\b|^2}\Bigr\|_1^{1/2}\|\nabla f\|_2 \|f\|_2
\\
&=\Biggl(\sqrt{\frac{C_{LS}}{2}}\theta^{-1/2}+\frac{\alpha}{2}\Biggr) \|\nabla f\|^2_2 -\Biggl(\frac{\sqrt{\alpha}}{\sqrt{2}}\|\nabla f\|_2 - \frac{\Bigl\|e^{\theta |\b|^2}\Bigr\|_1^{1/2} \theta^{-1/2}}{\sqrt{2\alpha}} \|f\|_2\Biggr)^2
+\frac{\Bigl\|e^{\theta |\b|^2}\Bigr\|_1 }{2\alpha \theta}\|f\|^2_2
\\
& \le (\sqrt{2C_{LS}}\theta^{-1/2}+\alpha) {\sf Ch}(f)
+\frac{\|e^{\theta |\b|^2}\|_1 }{2\alpha \theta}\|f\|^2_2
\\
&=  \e {\sf Ch}(f)+a\|f\|_2^2.
\end{align*}
By \eqref{alpha} and $\theta>2C_{LS}$, it holds that $\e \in (0,1).$ Therefore, \eqref{ineq: DD} has been proved. 

Assume (C). 
Since ${\rm div} \b \ge -c$ for some $c \ge 0$, it holds that 
\begin{align*}
\int_X \b(f)f dm = \frac{1}{2}\int_X \b(f^2) dm =-\frac{1}{2}\int_X {\rm div}\b |f|^2 dm 
\le \frac{c}{2}\int_X  |f|^2 dm.
\end{align*}
Therefore, \eqref{ineq: DD} has been proved. 
We now check \eqref{ineq: DD1}. By the Sobolev inequality \eqref{ineq: SL}, it holds that
\begin{align*} 
\Bigl| \int_X \b(f)g dm \Bigr|  \le\int_X |\b||\nabla f||g| dm  \le \|\nabla f\|_2\||\b|^2 g^2\|^{1/2}_1 
& \le  \|\nabla f\|_2 \||\b|^2\|^{1/2}_{N/2}\|g^2\|^{1/2}_{N/(N-2)} \notag
\\
& =  \|\nabla f\|_2 \|\b\|_{N}\|g\|_{2N/(N-2)} \notag
\\
& \le  \|\nabla f\|_2 \|\b\|_{N} \biggl(\|g\|_2^2+C_S\|\nabla g\|^2_{2}\biggr)^{1/2} \notag
\\
& \le  \|\nabla f\|_2 \|\b\|_{N} \bigl(\|g\|_2+\sqrt{C_S}\|\nabla g\|_{2}\bigr) \notag
\\
&= b{\sf Ch}^{1/2}_1(f) {\sf Ch}^{1/2}_1(g).
\end{align*}
We have proved \eqref{ineq: DD1}.

Now it is proved that $(\E, \F)$ is Markovian. It suffices to check that, for any non-negative constant  $a$ and $f \in {\rm Lip}_b(X)$, it holds that $\E(f \wedge a, f-f\wedge a) \ge 0$ (see \cite[Theorem 1.2.2]{O13}). 
Let $\{f<a\}=\{x \in X: f(x)<a\}$ for $f \in \mathcal D(\E)$ and a non-negative constant $a$. Noting the Markovian property of ${\sf Ch}$ and $\b(a)=0$, the following holds:
\begin{align*}
\E(f \wedge a, f \wedge a)& = {\sf Ch}(f \wedge a)-\int_{\{f<a\}} \b(f) f dm
\le {\sf Ch}(f \wedge a, f)-\int_{X} \b(f \wedge a) f dm = \E(f\wedge a, f).
\end{align*}
Thus $(\E, \F)$ is Markovian and $0 \le T_t f \le 1$ $m$-a.e.\ whenever $0 \le f \le 1$ $m$-a.e.\ with $f \in L^2(X,m)$.
The strong locality follows easily from the local property of ${\sf Ch}$ and the derivation $\b$. 

{It is now proved that the infinitesimal generator corresponding to $(\E, \F)$ coincides with $(L_2, \mathcal D(L_2))$ under (B).
Let $L$ be the infinitesimal generator associated with $(\E, \F)$. Let $f \in \mathcal D(L) \subset W^{1,2}(X,d,m)$. For $g \in \mathcal D(\Delta_q) \cap L^\infty(X,m)$ with any $q \ge  2$,
\begin{align*}
\int_X f (-\Delta_q g) dm=\int_X \langle \nabla f, \nabla g \rangle dm=\E(f,g)+\int_X \b(f)g dm
 = \int_X (-Lf+\b(f))gdm.
\end{align*}
Hence, $f \in \mathcal D(\Delta_p)$ and $\Delta_p f=Lf-\b(f)$ for any $1<p \le 2$. By Proposition \ref{prop: domain}, $f$ belongs to $\mathcal D(L_2)$ and $L_2f=Lf$. Therefore,  $L \subset L_2$. Thus, $(\lambda -L)^{-1} \supset (\lambda - L_2)^{-1}$ for sufficiently large $\lambda$.  Since both operators $(\lambda -L)^{-1}$ and $(\lambda - L_2)^{-1}$ are everywhere defined on $L^2(X,m)$, which implies the $L^2$-continuity of these operators,  and coincide on a dense subspace $\mathcal D(L_2) \subset L^2(X,m)$, it can be concluded that $L=L_2$.}
\qed

\section{Existence and Uniqueness of Invariant Measures} \label{sec: EU}
In this section, the existence and the uniqueness of invariant measures for the infinitesimal generator $(L, \mathcal D(L))$ are proved. 
Now we state the existence and the uniqueness of invariant measures.
\begin{thm} \label{thm: mthm2}
The following hold:
\begin{itemize}
\item[(i)]
Suppose one of the two assumptions (A) or (C). Then $(L, \mathcal D(L))$ has a unique 
(up to positive constant multiplication) non-zero invariant measure $\mu=\rho m$ in the sense of \eqref{eq: invL0} with $\rho \in L^2(X,m)$. In the case of (C), $\rho$ is strictly positive. 
\item[(ii)]
Suppose (B). Then $(L_p, \mathcal D(L_p))$ has a unique (up to positive constant multiplication) non-zero invariant measure $\mu=\rho m$ in the sense of \eqref{eq: invL0} with $\rho \in L^{\frac{p}{p-1}}(X,m)$. The density $\rho$ is strictly positive. 
\end{itemize}
\end{thm}
We first prove the case of (A). 
\\
{\it Proof for Theorem \ref{thm: mthm2} under Condition (A)}.
By Proposition \ref{prop: form}, $(\E, \F)$ is a positivity preserving coercive form under Condition (A). Therefore, \cite[Corollary 3.20]{BCR15} can be applied to show the existence of non-zero invariant measures for $(L, \mathcal D(L))$. 

The uniqueness of invariant measures can be proved by a similar way to \cite[Corollary 2.13]{BRZ00}, for which (2.12a) and (2.12b) in p. 96 in \cite{BRZ00} have been checked in Proposition \ref{prop: form}. Note that \cite{BRZ00} required the underlying space $X$ to be a locally convex topological vector space, but the proof of Corollary 2.13 is available also for the setting in Theorem \ref{thm: mthm2}. To apply \cite[Corollary 2.13]{BRZ00}, we only have to check the following Liouville-type property.
\begin{lem} \label{lem: LIU}
Let $(X,d,m)$ be an RCD$(K,\infty)$ space with $K \in \R$ and $m(X)=1$. 
If $u \in \mathcal D(\Delta)$ and $\Delta u=0$, then $u$ has an $m$-modification $\tilde{u}$ so that $\tilde{u}$ is a constant function.
\end{lem}
\proof
Take a function $u$ with $u \in \mathcal D(\Delta)$ and $\Delta u=0$. Then, it holds that
\begin{align*}
0=\int_Xu \Delta u dm=-\int_X\langle \nabla u, \nabla u \rangle dm.
\end{align*}
Therefore, $|\nabla u|=0$ $m$-a.e. By the Sobolev-to-Lipschitz property in \cite[Theorem 6.2]{AGS14b}, a Lipschitz modification $\tilde{u}$ of $u$ can be taken with the Lipschitz constant ${\rm Lip}(\tilde{u})\le  \||\nabla u|\|_\infty=0 $. Thus, $\tilde{u}$ is a constant function. 
\qed

We resume the proof for Theorem \ref{thm: mthm2} under Condition (A).
\\
{\it Proof for Theorem \ref{thm: mthm2} under Condition (A)}. By applying \cite[Corollary 2.13]{BRZ00}, we finished the proof. 
\qed

In the remainder of this section, Condition (B) or (C) is assumed.
We utilize \cite[Theorem 3.6]{H00}, which gives sufficient conditions for the existence and the uniqueness of strictly positive invariant measures. To apply \cite[Theorem 3.6]{H00}, we have to check two criteria: 
\begin{itemize}
\item[(i)] (\cite[Definition 2.1]{H00}) The semigroup $\{T_t\}$ is said to satisfy {\it Condition (I)$_p$} if,
for every $t>0$, there is a constant $J>0$, so that 
\begin{align} \label{ineq: I}
\psi_{T_{t}}(J)=\sup_{f \in L^p_+(X,m),\ \|f\|_p \le 1}\| (T_{t}f -J)_+\|_p<1.
\end{align}
\item[(ii)] (\cite[Definition 3.1]{H00})
Set 
$$\chi_{T_t}(\e)=\inf\Bigl\{ \int_X T_t \1_{B_1} \1_{B_2} dm: m(B_1) \ge \e, m(B_2) \ge \e\Bigr\}.$$
The semigroup $\{T_t\}$ is said to satisfy {\it Condition (E)} if, for each $\e>0$, there exists some $t >0$ so that $\chi_{T_{t}}(\e)>0$. 
\end{itemize}
\begin{lem} \label{lem: (I)} 
The following hold:
\begin{itemize}
\item[(i)] Under (B), 
the semigroup $\{T_t\}$ satisfies Condition (I)$_p$.
\item[(ii)] Under (C), 
the semigroup $\{T_t\}$ satisfies Condition (I)$_2$.
\end{itemize}
\end{lem}
\proof
We first check Condition (I)$_p$ under (B). By Proposition \ref{prop: L2}, the semigroup $\{T_t\}$ is a bounded operator from $L^p(X,m)$ to $L^q(X,m)$ for $1<p<q<\infty$ satisfying the condition in Proposition \ref{prop: L2}. Thus, by \cite[(i) in Lemma 2.2]{H00}, we can conclude that $\{T_t\}$ satisfies (I)$_p$.  

Now we check Condition (I)$_2$ under (C). By the Sobolev inequality \eqref{ineq: SL} and \eqref{ineq: 1}, we have that, for any $\alpha>\lambda$, there is a constant $C>0$, so that 
\begin{align} \label{ineq: BDDR}
\|f\|^2_{\frac{2N}{N-2}} \le C_S\| \nabla f\|_2^2 +\|f\|_2^2 \le C(\E_\alpha (f) + \|f\|_2^2), \quad \forall f \in \mathcal F.
\end{align}
Let $\{R_\alpha\}$ denote the resolvent of $(\E, \F)$. Then we have that $(f,g)=\E_\alpha (R_\alpha f, g)$ for any $\alpha>\lambda$ and any $f, g \in W^{1,2}(X,d,m)$. Thus, by \eqref{ineq: BDDR}, the $L^{2N/(N-2)}$-norm of $R_\alpha f$ can be estimated as follows:
\begin{align*}
\|R_\alpha f\|^2_{\frac{2N}{N-2}} \le C\E_\alpha (R_\alpha f) + C\|R_\alpha f\|_2^2 
\le  C\|f\|_2^2 + C\|R_\alpha f\|_2^2
 \le C\|f\|_2^2 + C\|R_\alpha\|^2_{2 \to 2}\|f\|_2^2.
\end{align*}
Since $\|R_\alpha\|_{2 \to 2} \le \frac{1}{\alpha-\lambda}$ (e.g., \cite[Theorem 1.1.2]{O13}), the resolvent $R_\alpha$ is a bounded operator from $L^2(X,m)$ to $L^{2N/(N-2)}(X,m)$ for each $\alpha >\lambda$. Therefore, by \cite[(i) in Lemma 2.2]{H00}, we can conclude that $\{T_t\}$ satisfies (I)$_2$.
\qed 

Now we show that $\{T_t\}$ satisfies Condition (E). 
It suffices to check the conditions (F1)--(F8) in \cite[Section 4]{H00}:
\begin{itemize}
\item[(F1)] there exist a positive constant $\lambda \ge 0$ so that $\E_\lambda (f) \ge 0$ for any $f \in \F$, and the domain $\F$ is dense in $L^2$ and closed under the norm $(\E_{\lambda+1})^{1/2}$. Furthermore, the weak sector condition holds: there exists a constant $C \ge 1$ so that
$$\E(f,g) \le C \E_{\lambda +1}^{1/2}(f) \E_{\lambda +1}^{1/2}(g);$$
\item[(F2)] the corresponding semigroup $\{T_t\}$ is positivity-preserving, that is, 
$$0 \le f \quad  m\text{-a.e.} \implies 0 \le T_t f \quad m\text{-a.e.\ for any $t > 0$;}$$
\item[(F3)] there exists a bilinear map $\Gamma: \F \times \F \to L^1(X,m)$ (not necessarily symmetric) so that 
$$\E(f,g)=\int_X \Gamma(f,g) dm, \quad f,g \in \F;$$
\item[(F4)] $\Gamma$ has a derivation property with respect to the first component: for any $f_i \in \F$ with $i=1,...,n$ and any $\phi \in C^\infty(\R^n)$, $\phi(f_1,...,f_n)$ belongs to $\F$ and 
$$\Gamma(\phi(f_1,...,f_n),g)=\sum_{i=1}^n \frac{\partial \phi}{\partial x_i}(f_1, ..., f_n) \Gamma(f_i, g), \quad \forall g \in \F.$$
\end{itemize}
Let us set the following set:
\begin{align} \label{eq: ExtS}
U=\{p>1: \{T_t\} \text{\ is extended or restricted to a strongly continuous semigroup on $L^p$}\}.
\end{align}
By the Riesz--Thorin interpolation theorem, it holds that $[2,\infty) \subset U$. 
Furthermore, the following conditions are considered:
there exist another bilinear map $\Gamma_0: \F \times \F \to L^1(X,m)$, functions $\eta \in \F$, $\chi \in L^2(X,m)$, $\sigma>0$, $\xi \ge 0$ $m$-a.e., and constants $r>0$, $\kappa>0$ so that 
\begin{itemize}
\item[(F5)] $\Gamma_0$ has a derivation property as in $(F4)$ with $n=1$ with respect to both components;
\item[(F6)] $\|\eta\|_2=1, \eta \in \bigcup_{s \in U}L^{2s}(X,m)$ and $\eta^2 \xi \in L^1(X,m)$;
\item[(F7)] the following Poincar\'e type inequality holds: for every $f \in \F \cap L^\infty(X,m)$,
\begin{align} \label{ineq: POINCARE}
\int_X \Gamma_0(f,f) \eta^2 dm \ge \Bigl( \int_X \Bigl| f-\int_X f\eta^2 dm\Bigr|^r \sigma dm\Bigr)^{2/r};
\end{align}
\item[(F8)] for every $f \in \F$, it holds that $\Gamma_0(f,f)+\xi f^2 \ge 0$ $m$-a.e., 
\begin{align*} 
&\Gamma(\eta, f) \ge \chi(\Gamma_0(f,f)+\xi f^2)^{1/2} \quad m\text{-a.e.,}
\end{align*}
and
\begin{align}\label{ineq: POINC3}
&\Gamma(f,f) \ge \kappa \Gamma_0(f,f) - \xi f^2\quad m\text{-a.e.}
\end{align}
\end{itemize}

\begin{lem} \label{prop: F1-F8}
Under (B) or (C), 
the form $(\E, \mathcal F)$ satisfies all conditions (F1)--(F8) with 
$$\Gamma(f)=\frac{1}{2}|\nabla f|^2-\b(f)f, \quad \Gamma_0(f)=\frac{1}{2}|\nabla f|^2,$$ 
and $\eta \equiv 1$, $\chi \equiv 0$, $r=2$, $\sigma=C_P$, $\kappa=1/2$ and $\xi=|\b|^2$.
Here $C_P$ denotes the constant appearing in the Poincar\'e inequality in \eqref{ineq: PI}.
\end{lem}
\proof
(F1) has been already checked in Proposition \ref{prop: form}.
(F2) is now checked. It suffices to show that, for any $f \in \F$, it holds that $f_+ \in \F$ and $\E_{\lambda}(f_+, f) \ge 0$, whereby $f_+=f \vee 0$.
Noting that $\mathcal F=W^{1,2}(X,d,m)$, it is clear that $f_+ \in \F$ whenever $f \in \F$. By \eqref{ineq: 1}, the following holds:
\begin{align*}
\E_{\lambda}(f_+, f)& =\E_{\lambda}(f\1_{\{f \ge 0\}}) 
\ge c\|f\1_{\{f \ge 0\}}\|^2_{1,2} 
 \ge 0.
\end{align*}
Thus (F2) has been checked. 

(F3) can be easily obtained by taking $\Gamma(f,g)=(1/2)\langle \nabla f, \nabla g \rangle - \b(f)g$. The derivation property (F4) holds by making use of the derivation property of $\b$ and the inner product $\langle \cdot, \cdot \rangle$ in $L^2(TX)$. By taking $\Gamma_0(f,g)=(1/2)\langle \nabla f, \nabla g \rangle$, Condition (F5) and (F6) can be checked. Since $({\sf Ch}, W^{1,2}(X,d,m))$ satisfies the Poincar\'e inequality owing to the RCD$(K,\infty)$ condition with $K>0$, or with a bounded diameter as in \eqref{ineq: PI}, Condition (F7) holds. Concerning (F8), it suffices to check \eqref{ineq: POINC3}. It can be seen by the following argument with $\kappa=1/2$ and $\xi=|\b|^2$ that 
\begin{align*}
\Gamma(f) -\frac{1}{2} \Gamma_0(f) +|\b|^2 f^2 =\frac{1}{4}|\nabla f|^2-\b(f)f  + |\b|^2 f^2
 \ge \frac{1}{4}|\nabla f|^2-|\b||\nabla f||f|  + |\b|^2 f^2
& = \frac{1}{4} \Bigl( |\nabla f|-2|\b|f\Bigr)^2
\\
& \ge 0.
\end{align*}
Thus, Condition (F8) has been checked. 
\qed

{\it Proof for Theorem \ref{thm: mthm2} under Condition (B) or (C)}.
By Lemma \ref{lem: (I)} and \ref{prop: F1-F8}, we can apply \cite[Theorem 3.6, Proposition 4.3]{H00}, which concludes the existence and the uniqueness of strictly positive invariant measures. 
\qed

\section{Regularity of Invariant Measures} \label{sec: NSDF}
Throughout this section, we always assume without explicitly mentioning that 
\begin{align} \label{con: AssReg}
\text{$(X,d,m)$ is an RCD$(K,\infty)$ space with $K >0$, or ${\rm diam}(X)<\infty$, and $\b$ satisfies $|\b| \in L^1(X,\rho m)$.}
\end{align}
Let $\rho \in L^1_+(X,m)$ satisfy the following equality:
\begin{align} \label{eq: invL5-1}
\int_X \Bigl(\frac{1}{2}\Delta+\b\Bigr)(\phi) \rho dm =0, \quad \forall \phi \in {\rm TestF}(X).
\end{align} 
Note that \eqref{eq: invL5-1} is weaker than \eqref{eq: invL0} since ${\rm TestF}(X)$ is not necessarily dense in the domain of $(L, \mathcal D(L))$ with respect to the graph norm, but we also call $\mu=\rho m$ invariant measure. 
In this section, the Sobolev regularity of the density $\rho$ is proved. 
Recall that $\{H_t\}$ denotes the heat semigroup and $\rho_\e:=H_\e \rho$.
\begin{lem} \label{lem: 1}
If one of the following conditions holds:
\begin{itemize}
\item[(i)] $\rho \in L^2(X,m)$ and $|\b|\rho \in L^2(X,m)$;
\item[(ii)] $H_t: L^1(X,m) \to L^\infty(X,m)$ for any $t>0$ (i.e., $\{H_t\}$ is ultra-contractive);
\end{itemize}
then it holds that
\begin{align} \label{eq: Inv}
\int_{X} \langle \nabla \phi, \nabla \rho_{\e} \rangle dm= \int_{X}\b(H_\e \phi) d\mu, \quad \forall \phi \in W^{1,2}(X,d,m).
\end{align}
\end{lem}
\proof
For any $\phi \in {\rm TestF}(X)$, the following holds:
\begin{align} \label{eq: Inv0}
\int_{X} \langle \nabla \phi, \nabla \rho_{\e} \rangle dm &= -\int_{X} \Delta \phi \rho_{\e} dm 
= -\int_{X} \Delta \phi H_\e \rho dm  
= -\int_{X} H_\e \Delta \phi \rho dm  
= -\int_{X} \Delta H_\e \phi d\mu.
\end{align}
By \eqref{eq: invL5-1} and $H_\e \phi \in {\rm TestF}(X)$, it can be seen that 
\begin{align}
-\int_{X} \Delta H_\e \phi d\mu = \int_{X}\b(H_\e \phi) d\mu. \notag
\end{align}
Therefore, by \eqref{eq: Inv0}, the following equality is obtained:
$$\int_{X} \langle \nabla \phi, \nabla \rho_{\e} \rangle dm =\int_{X}\b(H_\e \phi) d\mu, \quad \forall \phi \in {\rm TestF}(X).$$
We now prove that the above equality holds for $\phi \in W^{1,2}(X,d,m)$.
Let $\{\phi_i\}_{i \in \N} \subset {\rm TestF}(X)$ with $\phi_i \to \phi$ in $W^{1,2}(X,d,m)$ as $i \to \infty$. 
Then, the following holds: 
\begin{align}
&\biggl| \int_{X} \langle \nabla \phi, \nabla \rho_{\e} \rangle dm- \int_{X}\b(H_\e \phi) d\mu \biggr|  \notag
\\
&\le \biggl| \int_{X} \langle \nabla \phi, \nabla \rho_{\e} \rangle dm- \int_{X} \langle \nabla \phi_i, \nabla \rho_{\e} \rangle dm \biggr| \notag
+ \biggl|\int_{X} \langle \nabla \phi_i, \nabla \rho_{\e} \rangle dm - \int_{X}\b(H_\e \phi_i) d\mu\biggr| \notag
\\
& \quad + \biggl|\int_{X}\b(H_\e \phi_i) d\mu -\int_{X}\b(H_\e \phi) d\mu \biggr| \notag
\\
&=\biggl| \int_{X} \langle \nabla \phi, \nabla \rho_{\e} \rangle dm- \int_{X} \langle \nabla \phi_i, \nabla \rho_{\e} \rangle dm \biggr|\notag
+ \biggl|\int_{X}\b(H_\e \phi_i) d\mu -\int_{X}\b(H_\e \phi) d\mu \biggr| 
\\
&\le\int_{X} \bigl|  \langle \nabla (\phi-\phi_i), \nabla \rho_{\e} \rangle \bigr| dm  \label{ineq: DAI}
+ \int_{X}|\b|\bigl| \nabla H_\e (\phi_i-\phi)\bigr| d\mu.
\end{align}
Since the first term in the last line goes to zero as $i \to \infty$ by $\phi_i \to \phi$ in $W^{1,2}(X,d,m)$, it suffices to show that the second term goes to zero as $i \to \infty$. 
If $|\b|\rho \in L^2(X,m)$ or $\{H_t\}$ is ultra-contractive, then, by the gradient estimate \eqref{ineq: GE}, the second term in \eqref{ineq: DAI} can be estimated as follows:
\begin{align*}
\int_{X}|\b|\bigl| \nabla \bigl(H_\e \phi_i-H_\e \phi \bigr)\bigr| d\mu 
\le e^{-K\e/2}\int_{X} H_\e \bigl| \nabla (\phi_i- \phi) \bigr| |\b|\rho dm 
& = e^{-K\e/2}\int_{X} \bigl| \nabla (\phi_i- \phi) \bigr|  H_\e (|\b|\rho)dm.
\\
& \to 0 \quad (i \to \infty).
\end{align*}
Thus, the second term in \eqref{ineq: DAI} goes to zero as $i \to \infty$, which implies the equality \eqref{eq: Inv}.
\qed

The regularity of the density $\rho$ is now presented. 
\begin{thm} \label{thm: mthm1}
Suppose \eqref{con: AssReg} and \eqref{eq: invL5-1}.
Assume $\rho, |\b|\rho \in L^2(X,m)$. Then $\rho \in W^{1,2}(X,d,m)$ and the following inequality holds:
   \begin{align} \label{ineq: SNE}
   \|\nabla \rho\|_2  \le \||\b|\rho\|_2. 
   \end{align}
   \end{thm}

\proof
By $|\b|\rho \in L^2(X,m)$, the equality \eqref{eq: Inv} in Lemma \ref{lem: 1} is available. Taking $\phi=\rho_{\e}$, we see the following estimate:
\begin{align}
\int_{X} \langle \nabla \phi, \nabla \rho_{\e} \rangle dm= \int_{X}|\nabla \rho_{\e}|^2 dm
&=\int_{X} \b(H_\e \rho_{\e}) d\mu \notag
\\
&\le \int_{X} |\b|| \nabla H_\e \rho_{\e}| d\mu \notag
\\
&= \int_{X}| \nabla H_\e \rho_{\e}| |\b|\rho dm \notag
\\
&\le e^{-K\e/2}\int_{X} H_\e \bigl( | \nabla \rho_{\e}|\bigr) |\b|\rho dm \notag
\\
&= e^{-K\e/2}\int_{X}| \nabla \rho_{\e}| H_\e (|\b|\rho) dm \notag 
\\
&\le e^{-K\e/2}\biggl(\int_{X}| \nabla \rho_{\e}|^2dm\biggr)^{1/2} \biggl(\int_{X}|H_\e(|\b|\rho)|^2 dm\biggr)^{1/2} \notag
\\
& \le e^{-K\e/2}\biggl(\int_{X}| \nabla \rho_{\e}|^2dm\biggr)^{1/2} \biggl(\int_{X}|\b|^2|\rho|^2 dm\biggr)^{1/2}. \notag
\end{align} 
Thus, it has been obtained that
\begin{align} \label{ineq: SES}
\biggl(\int_{X}|\nabla \rho_{\e}|^2 dm\biggr)^{1/2} \le e^{-K\e/2}\biggl(\int_{X}|\b|^2|\rho|^2 dm\biggr)^{1/2}<\infty.
\end{align}
Therefore, $\{\rho_\e\}_{\e>0}$ is a bounded sequence in $W^{1,2}(X,d,m)$. By \cite[Theorem 6.3, Proposition 6.7]{GMS13}, 
under \eqref{con: AssReg}, there exists a subsequence $\{\rho_{\e'}\}$ and $\rho' \in W^{1,2}(X,d,m)$ so that $\rho_{\e'}$ converges to $\rho'$ strongly in $L^2(X,m)$ and weakly in $W^{1,2}(X,d,m)$ as $\e'$ goes to zero. Since $\{H_t\}$ is a strongly continuous $L^2$-semigroup, we have that $\rho_\e$ converges to $\rho$ in $L^2(X,m)$ strongly and $W^{1,2}(X,d,m)$ weakly as $\e$ goes to zero. Thus, $\rho=\rho'$ and $\rho \in W^{1,2}(X,d,m).$ Furthermore, by \eqref{ineq: SES} and the weak $W^{1,2}$-convergence of $\{\rho_\e\}$, we obtain that  
$$\biggl(\int_{X}|\nabla \rho|^2 dm\biggr)^{1/2} \le \liminf_{\e \to 0}\biggl(\int_{X}|\nabla \rho_{\e}|^2 dm\biggr)^{1/2} \le \limsup_{\e \to 0}\biggl(\int_{X}|\nabla \rho_{\e}|^2 dm\biggr)^{1/2} \le \biggl(\int_{X}|\b|^2|\rho|^2 dm\biggr)^{1/2}.$$
Thus, we have that $\rho \in W^{1,2}(X,d,m)$ and 
$$\|\nabla \rho\|_2 \le \||\b|\rho\|_2.$$  
We have obtained the desired result. \qed

\begin{rem} \normalfont 
The Sobolev regularity of $\rho$ in the case of the Wiener space with a bounded drift has been shown by Shigekawa \cite{Shi87}, which relies on the finite-dimensional approximation of the Wiener space.
Theorem \ref{thm: mthm1} relies only on the regularity of the heat semigroup $\{H_t\}$ induced by the Bakry--\'Emery gradient estimate, which gives a different proof for the Wiener space. 
\end{rem}

In order to show the Sobolev regularity $\rho \in W^{1,2}(X,d,m)$ by applying Theorem \ref{thm: mthm1}, it is necessary to know when the assumption $\rho, |\b|\rho \in L^2(X,m)$ in Theorem \ref{thm: mthm1} holds.
For instance, 
\begin{itemize}
\item[(a)] Let $|\b| \in L^\infty(X,m)$ and $\rho \in L^2(X,m)$. In this case, the following inequality holds:
\begin{align} \label{ineq: GER}
\|\nabla \rho\|_2 \le  \||\b|\rho\|_2 \le \|\b\|_\infty\|\rho\|_2<\infty.
 \end{align}
 \item[(b)] Let $1<p<2$. Let $|\b| \in L^{2p/(2-p)}(X,m)$ and $\rho \in L^{p/(p-1)}(X,m)$. In this case, the following inequality holds:
 \begin{align} \label{ineq: GER5}
 \|\nabla \rho \|_2 \le \||\b|\rho\|_2 \le \|\b\|_{\frac{2p}{2-p}}\| \rho\|_{\frac{p}{p-1}}<\infty.
 \end{align}
 \end{itemize}
In the case of finite-dimensional spaces, it is proved in the next theorem that the condition of the integrability of $|\b|$ and $\rho$ can be relaxed: 
\begin{thm} \label{thm: mthm1-2}
Assume \eqref{eq: invL5-1} and that $(X,d,m)$ is an RCD$(K,N)$ space with $2<N<\infty$ and $K>0$, and $|\b| \in L^N(X,\mu)$.  
Then, $\rho \in W^{1,2}(X,d,m) \cap L^p(X,m)$ for any $1 \le p<\infty$ and the following inequality holds:
   \begin{align} \label{ineq: GER6}
\|\nabla \rho\|_2 \le  \||\b|\rho\|_2 \le \|\b\|_{L^N(X,d\mu)} \|\rho\|_{\frac{N}{N-2}+1}^{\frac{N-1}{N}}<\infty.
  \end{align}
 If, moreover, $|\b| \in L^k(X,\mu)$ for some $k>N$, then $\rho \in L^\infty(X,m)$. 
 \end{thm}
  \begin{rem} \normalfont \label{rem: 2}
  In the case of the Euclidean space, the Sobolev regularity in Theorem \ref{thm: mthm1-2} was studied by Metafune--Pallara--Rhandi \cite{MPR05}. Their proof relies on the Friedrichs mollifier, which is not available in our non-smooth setting. Alternatively, the heat semigroup $\{H_t\}$ plays a role as a mollifier in our setting, whose regularity comes from the lower Ricci curvature bounds. Compared to the Euclidean case, our proof needs more careful arguments on the integrability estimate concerning the heat semigroup action. We also note that the Sobolev inequality \eqref{ineq: SL} is different from the Sobolev inequality in the Euclidean space, which will be utilized in order to do Moser's iteration in Lemma \ref{thm: Lp}. 
  \end{rem}

In what follows in this section, we focus on the proof of Theorem \ref{thm: mthm1-2}. Note that under the assumption of Theorem \ref{thm: mthm1-2}, ${\rm diam}(X)<\infty$ and ${\sf m}(X)<\infty$ by the generalised Bonnet--Myers theorem \cite[Corollary 2.6]{Sturm06-2}.
We may assume that $m(X)=1$ by normalizing $m/m(X)$. 

We prove several lemmas. We first give an gradient estimate of $\sqrt{\rho}$.
\begin{lem} \label{lem: W12-1}
Assume \eqref{eq: invL5-1} and that $(X,d,m)$ is an RCD$(K,N)$ space with $2<N<\infty$, $K \in \R$, and ${\rm diam}(X)<D<\infty$.  Assume $|\b| \in L^2(X,\mu)$. Then $\sqrt{\rho} \in W^{1,2}(X,d,m)$ and 
\begin{align} \label{ineq: GER2}
\|\nabla \sqrt{\rho}\|_{2} 
\le \|\b\|_{L^2(X,\mu)}.
\end{align}
\end{lem}
\proof
Since $\rho \in L^1(X,m)$ and $\{H_t\}$ is ultra-contractive by \cite{JLZ15} and the Bishop--Gromov inequality,   by applying \cite[Propisition 6.4]{AGS14b}, we have that $\rho_\e:=H_\e \rho \in {\rm Lip}_b(X) \subset W^{1,2}(X,d,m)$. Let $\phi_\e$ be defined as follows: 
$$\phi_\e=\log(\rho_\e+\delta)-\log \delta \in W^{1,2}(X,d,m) \quad \delta>0.$$
Applying \eqref{eq: Inv} with $\phi=\phi_\e$, the gradient estimate \eqref{ineq: GE} and the chain rule for $\phi_\e$ (see \cite[Theorem 2.2.6]{G16}), we have that
\begin{align*}
\int_X \frac{|\nabla \rho_\e|^2}{(\rho_\e+\delta)} dm 
= \int_{X}\b(H_\e \phi_\e ) d\mu
 & \le  \int_{X}|\b| | \nabla H_\e \phi_\e | d\mu
 \\
 & \le  \biggl( \int_X |\b|^2 d\mu \biggr)^{1/2} \biggl( \int_{X}| \nabla H_\e \phi_\e |^2 d\mu \biggr)^{1/2}
  \\
 & \le  e^{-K \e}\biggl( \int_X |\b|^2 d\mu \biggr)^{1/2} \biggl( \int_{X}|\nabla \phi_\e |^2 H_\e (\rho) dm \biggr)^{1/2}
   \\
 & =  e^{-K \e}\biggl( \int_X |\b|^2 d\mu \biggr)^{1/2} \biggl( \int_{X}\frac{|\nabla \rho_\e|^2}{(\rho_\e+\delta)}\frac{\rho_\e}{(\rho_\e+\delta)}  dm \biggr)^{1/2}
    \\
 & \le  e^{-K \e}\biggl( \int_X |\b|^2 d\mu \biggr)^{1/2} \biggl( \int_{X}\frac{|\nabla \rho_\e|^2}{(\rho_\e+\delta)}dm \biggr)^{1/2}.
\end{align*}
Thus, by the Fatou lemma, taking $\delta \to 0$, we have that  
\begin{align*}
\Bigl( \int_X \frac{|\nabla \rho_\e|^2}{\rho_\e} dm \Bigr)^{1/2} \le \liminf_{\delta \to 0}\Bigl( \int_X \frac{|\nabla \rho_\e|^2}{(\rho_\e+\delta)} dm \Bigr)^{1/2} \le e^{-K\e}\|\b\|_{L^2(X,\mu)}.
\end{align*}
Therefore, it holds that
\begin{align*}
\limsup_{\e \to 0}\Bigl( \int_X |\nabla \sqrt{\rho_\e}|^2 dm \Bigr)^{1/2} 
=\limsup_{\e \to 0} \Bigl( \int_X \frac{|\nabla \rho_\e|^2}{\rho_\e} dm \Bigr)^{1/2} \le \|\b\|_{L^2(X,\mu)}.
\end{align*}
Hence, $\{\sqrt{\rho_\e}\}_\e$ is a bounded sequence in $W^{1,2}(X,d, m)$.
Therefore, there is a converging subsequence $\{\sqrt{\rho_{\e'}}\}_{\e'}$ in the $L^2$-strong sense and the $W^{1,2}$-weak sense. Since the measure $\rho_\e m$ converges weakly to $\rho m$ as $\e \to 0$, we have that $\sqrt{\rho_{\e}} \to \sqrt{\rho}$ weakly in $W^{1,2}(X,d,m)$ and strongly in $L^2$ as $\e \to 0$. Thus, we have that 
\begin{align*}
\Bigl( \int_X |\nabla \sqrt{\rho}|^2 dm \Bigr)^{1/2} 
\le \liminf_{\e \to 0}\Bigl( \int_X |\nabla \sqrt{\rho_\e}|^2 dm \Bigr)^{1/2}
 \le \limsup_{\e \to 0}\Bigl( \int_X |\nabla \sqrt{\rho_\e}|^2 dm \Bigr)^{1/2}
&  \le \|\b\|_{L^2(X,\mu)}.
\end{align*}
We finished the proof. \qed

We prove the following key lemma:
\begin{lem} \label{lem: Sob-M}
Assume \eqref{eq: invL5-1} and that $(X,d,m)$ is an RCD$(K,N)$ space with $2<N<\infty$ and $K>0$ and $|\b| \in L^{k}(X,\mu)$ for some $k >2$. Fix $\alpha>0$.  If $\rho \in L^{\alpha\frac{k}{k-2}+1}(X,m)$, then 
\begin{align} \label{ineq: Sob-M}
\int_X \rho^{\alpha-1}|\nabla \rho|^2 dm \le \int_X |\b|^2\rho^{\alpha+1} dm<\infty.
\end{align}
\end{lem}
\proof 
It can be checked easily that 
\begin{align} \label{ineq: Sob-M1}
\int_X |\b|^2\rho^{\alpha+1} dm =\int_X |\b|^2\rho^{2/k}\rho^{\alpha+1-2/k} dm  
\le\biggl( \int_X |\b|^k \rho dm \biggr)^{2/k} \biggl( \int_X \rho^{\alpha\frac{k}{k-2}+1} dm\biggr)^{(k-2)/k} 
<\infty.
\end{align}
Since $(X,d,m)$ is an RCD$(K,N)$ space with $K>0$, the heat semigroup $\{H_t\}$ is ultra-contractive. 
Thus, \eqref{eq: Inv} can be applied. 

We first prove the case $\alpha \ge 1$. 
Take $\phi=(\rho_\e)^\alpha$ (just simply write $\rho_\e^\alpha=(\rho_\e)^\alpha$ hereinafter). Note that $\phi \in W^{1,2}(X,d,m)\cap L^\infty(X,m)$ owing to the ultra-contrativity of $\{H_t\}$. Then we see that (see \cite[Theorem 2.2.6]{G16})
\begin{align} \label{eq: DFo1}
\nabla \phi= \alpha \rho_\e ^{\alpha-1}\nabla \rho_\e.
\end{align}
Plugging $\phi$ into \eqref{eq: Inv}, we have that 
\begin{align*}
\alpha \int_X \rho_\e^{\alpha-1}|\nabla \rho_\e|^2 dm  
=\int_X \langle \nabla \phi, \nabla \rho_\e \rangle dm 
&= \int_{X}\b(H_\e \phi ) d\mu \notag
\\
& \le  \int_{X}|\b|| \nabla H_\e \phi | d\mu \notag
\\
& \le  e^{-K\e/2} \int_{X} |\nabla \phi | H_\e (|\b| \rho)dm \notag
\\
& =  \alpha e^{-K\e/2} \int_{X}\rho_\e^{\alpha-1} |\nabla \rho_\e | H_e(|\b| \rho) dm \notag
\\
& \le  \alpha e^{-K\e/2} \biggl( \int_{X}\bigl( H_\e (|\b| \rho) \bigr)^{2} \rho_\e^{\alpha-1} dm\biggl)^{1/2} \biggl( \int_X \rho_\e^{\alpha-1} |\nabla \rho_\e |^2 dm\biggr)^{1/2}.
\end{align*}
Therefore, we have that 
\begin{align} \label{ineq: IIII}
\biggl( \int_X  \rho_\e^{\alpha-1}|\nabla \rho_\e|^2 dm\biggr)^{1/2} \le e^{-K\e/2} \biggl( \int_{X}\bigl( H_\e (|\b| \rho) \bigr)^{2} \rho_\e^{\alpha-1} dm\biggl)^{1/2}.
\end{align}
We now show that the right-hand side in \eqref{ineq: IIII} converges as $\e \to 0$. By the assumption $\rho \in L^{\alpha\frac{k}{k-2}+1}(X,m)$, it can be checked easily that $\rho^{\alpha-1} \in L^\beta(X,m)$ so that $\beta$ is determined by the following equality:
$$(\alpha-1)\beta=\alpha\frac{k}{k-2}+1.$$
By a simple calculation, we have that 
$$ \beta=\frac{\alpha k +k -2}{(\alpha-1)(k-2)}.$$ 
Therefore, $\rho_\e^{\alpha-1}$ converges to $\rho^{\alpha-1}$ strongly in $L^\beta(X,m)$ by the strong continuity of $\{H_t\}$ on $L^p(X,m)$, i.e., $H_t f \to f$ in $L^p(X,m)$ as $t \to 0$ for every $f \in L^p(X,m)$ ($1 \le p < \infty$). 
Let $\beta_*$ denote the conjugate exponent of $\beta$, that is, 
$$\beta_*=\frac{\beta}{\beta-1}=\frac{\alpha k +k -2}{2k +2\alpha-4}.$$
In order to show that the right-hand side in \eqref{ineq: IIII} converges, it suffices to show that $( H_\e (|\b| \rho) )^{2}$ converges to $|\b|^2 \rho^2$ strongly in $L^{\beta_*}(X,m)$. To show this, we only have to show that $|\b|^2\rho^2 \in L^{\beta_*}(X,m)$ because of the strong continuity of $\{H_t\}$ on $L^p(X,m)$ for any $1 \le p<\infty$. We have that 
\begin{align} \label{ineq: IVI}
\int_X |\b|^{2\beta_*}\rho^{2\beta_*} dm =\int_X |\b|^{2\beta_*}\rho^{2\beta_*/k} \rho^{2\beta_*-2\beta_*/k} dm
\le\biggl( \int_X |\b|^k \rho dm \biggr)^{2\beta_*/k} \biggl( \int_X \rho^{\frac{4\beta^2_*(k-1)}{k(k-2\beta_*)}} dm\biggr)^{\frac{k-2\beta_*}{2\beta_*}}.
\end{align}
By the assumption $|\b| \in L^{k}(X,\mu)$, it holds that $\int_X |\b|^k \rho dm<\infty$. 
By a simple calculation, 
$$\frac{4\beta^2_*(k-1)}{k(k-2\beta_*)}= \frac{4\beta^2_*(k+\alpha-2)}{k(k-2)}.$$ 
Since $k>2$ and $\alpha \ge 1$, we see that 
\begin{align*}
\alpha\frac{k}{k-2}+1-\frac{4\beta^2_*(k-1)}{k(k-2\beta_*)} 
&= \alpha\frac{k}{k-2}+1-\frac{4\beta^2_*(k+\alpha-2)}{k(k-2)}
\\
&=\frac{1}{k(k-2)}(\alpha k^2+k(k-2)-4\beta^2_*(k+\alpha-2) )
\\
&=\frac{1}{k(k-2)}\biggl(\alpha k^2+k(k-2)-4\frac{(\alpha k +k -2)^2}{(2k +2\alpha-4)^2}(k+\alpha-2) \biggr)
\\
&=\frac{1}{k(k-2)}\biggl(k\bigl(\alpha k+(k-2)\bigr)-\frac{(\alpha k +k -2)^2}{(k +\alpha-2)} \biggr)
\\
&=\frac{(\alpha k +k -2)}{k(k-2)}\biggl(k-\frac{(\alpha k +k -2)}{(k +\alpha-2)} \biggr)
\\
&=\frac{(\alpha k +k -2)}{k(k-2)}\frac{(k-1)(k-2)}{(k +\alpha-2)}
\\
&=\frac{(\alpha k +k -2)}{k}\frac{(k-1)}{(k +\alpha-2)}
\\
& > 0.
\end{align*} 
Therefore, we have that 
$$\frac{4\beta^2_*(k-1)}{k(k-2\beta_*)}  < \alpha\frac{k}{k-2}+1.$$
By the assumption $\rho \in L^{\alpha\frac{k}{k-2}+1}(X,m)$ and the fact that $\|f\|_p \le \|f\|_q$ whenever $p \le q$ since $m(X)=1$, we obtain that 
$$\biggl( \int_X \rho^{\frac{4\beta^2_*(k-1)}{k(k-2\beta_*)}} dm\biggr)^{\frac{k-2\beta_*}{2\beta_*}} =\biggl( \|\rho\|_{\frac{4\beta^2_*(k-1)}{k(k-2\beta_*)}}^{\frac{k-2\beta_*}{2\beta_*}} \biggr)^{\frac{k(k-2\beta_*)}{4\beta^2_*(k-1)}}  \le \biggl( \|\rho\|_{\alpha\frac{k}{k-2}+1}^{\frac{k-2\beta_*}{2\beta_*}}\biggr)^{\frac{k(k-2\beta_*)}{4\beta^2_*(k-1)}}  <\infty.$$
Thus, \eqref{ineq: IVI} is finite, which concludes that $|\b|^2\rho^2 \in L^{\beta_*}(X,m)$. 
Therefore,  the right-hand side of \eqref{ineq: IIII} converges as $\e \to 0$ and 
\begin{align*} 
\limsup_{\e \to 0}\biggl( \int_X  \rho_\e^{\alpha-1}|\nabla \rho_\e|^2 dm\biggr)^{1/2} 
\le \biggl( \int_{X}|\b|^2  \rho^{\alpha+1} dm\biggl)^{1/2}.
\end{align*}
 Note that, by \eqref{eq: DFo1}, it holds that 
 $$\frac{2}{\alpha+1}\biggl(\int_X \Bigl| \nabla \rho^{\frac{\alpha+1}{2}}_\e\Bigr|^2 dm\biggr)^{1/2}=\biggl( \int_X \rho_\e^{\alpha-1} |\nabla \rho_\e|^2 dm\biggr)^{1/2}.$$
Thus, we see that
 \begin{align*}
  \limsup_{\e \to 0}\frac{2}{\alpha+1}\biggl(\int_X \Bigl| \nabla \rho^{\frac{\alpha+1}{2}}_\e\Bigr|^2 dm\biggr)^{1/2}
  = \limsup_{\e \to 0}\biggl( \int_X \rho_\e^{\alpha-1} |\nabla \rho_\e|^2 dm\biggr)^{1/2}
  \le \biggl( \int_{X}|\b|^2  \rho^{\alpha+1} dm\biggl)^{1/2}.
 \end{align*}
Note that $\rho^{(\alpha+1)/2}_\e \in L^2(X,m)$ since $\alpha+1 \le \alpha\frac{k}{k-2}+1$ and $\rho \in L^{\alpha\frac{k}{k-2}+1}(X,m)$ by the assumption.  Therefore, $\{ \rho^{(\alpha+1)/2}_\e\}_\e$ is bounded in $W^{1,2}(X,d,m)$, which implies that there is a converging subsequence weakly in $W^{1,2}(X,d,m)$ and strongly in $L^2(X,m)$. 
Since $\rho_\e$ converges to $\rho$ strongly in $L^{\alpha\frac{k}{k-2}+1}(X,m)$ by the strong continuity of $\{H_t\}$, we have that the full sequence $\{ \rho^{\frac{\alpha+1}{2}}_\e\}_\e$ converges to $ \rho^{\frac{\alpha+1}{2}}$ weakly in $W^{1,2}(X,d,m)$ and strongly in $L^2(X,m)$. 
Thus, we obtain that  
\begin{align*} 
\biggl( \int_X \rho^{\alpha-1} |\nabla \rho|^2 dm\biggr)^{1/2} = \frac{2}{\alpha+1}\biggl(\int_X \Bigl| \nabla \rho^{\frac{\alpha+1}{2}} \Bigr|^2 dm\biggr)^{1/2}
 &\le  \frac{2}{\alpha+1}\liminf_{\e \to 0}\biggl(\int_X \Bigl| \nabla \rho^{\frac{\alpha+1}{2}}_\e\Bigr|^2 dm\biggr)^{1/2}
 \notag
 \\
  &\le \frac{2}{\alpha+1}\limsup_{\e \to 0}\biggl(\int_X \Bigl| \nabla \rho^{\frac{\alpha+1}{2}}_\e\Bigr|^2 dm\biggr)^{1/2}
\\
 & =\limsup_{\e \to 0}\biggl( \int_X \rho_\e^{\alpha-1} |\nabla \rho_\e|^2 dm\biggr)^{1/2}
 \\
 & \le \biggl( \int_{X}|\b|^2  \rho^{\alpha+1} dm\biggl)^{1/2}.
 \end{align*}
 The first equality holds by applying the chain rule in \cite[Theorem 2.2.6]{G16}. In fact, by \cite[Theorem 2.2.6]{G16}, we have that $|\nabla (\rho \wedge M)^{(\alpha+1)/2}|^2=\rho^{\alpha-1} |\nabla \rho|^2\1_{\{\rho \le M\}}$ $m$-a.e.\ for $0<M<\infty$. By $M \to \infty$ and the monotone convergence theorem, we have the first equality. 
We have finished the proof for $\alpha \ge 1$. 

We now prove the case $0<\alpha<1.$ Let $\psi_{\alpha, \eta}(x)=(x+\eta)^{\alpha}$ for $\eta>0$. We use the second subscription $\eta$ only where necessary, and otherwise write $\psi_\alpha$. 
Take $\phi=\psi_\alpha \circ \rho_\e$. Note that $\phi \in W^{1,2}(X,d,m)\cap L^\infty(X,m)$ owing to the ultra-contrativity of $\{H_t\}$. Then we see that 
\begin{align} \label{eq: DFo}
\nabla \phi = \alpha \psi_{\alpha-1}\circ \rho_\e\nabla \rho_\e.
\end{align}
Plugging $\phi$ into \eqref{eq: Inv}, we have that 
\begin{align} \label{ineq: InvM-1}
\alpha \int_X \psi_{\alpha-1}\circ \rho_\e|\nabla \rho_\e|^2 dm  \notag
& =\int_X \langle \nabla \phi, \nabla \rho_\e \rangle dm \notag
\\
& = \int_{X}\b(H_\e \phi ) d\mu \notag
\\
& \le  \int_{X}|\b|| \nabla H_\e \phi | d\mu \notag
\\
& \le  e^{-K\e/2} \int_{X} |\nabla \phi | H_\e (|\b| \rho)dm \notag
\\
& =  \alpha e^{-K\e/2} \int_{X}\psi_{\alpha-1}\circ \rho_\e |\nabla \rho_\e | H_\e(|\b| \rho) dm \notag
\\
& \le  \alpha e^{-K\e/2} \biggl( \int_{X}\bigl( H_\e (|\b| \rho) \bigr)^{2} \psi_{\alpha-1}\circ \rho_\e dm\biggl)^{1/2} \biggl( \int_X \psi_{\alpha-1}\circ \rho_\e |\nabla \rho_\e |^2 dm\biggr)^{1/2}.
\end{align}
 By the assumption $|\b| \in L^{k}(X, \mu)$, we have that 
\begin{align} \label{ineq: InvM-2}
H_\e(|\b|\rho) 
= H_\e\Bigl(|\b|\rho^{\frac{1}{k}}\rho^{1-\frac{1}{k}}\Bigr)  
&\le \Bigl(H_\e\bigl(|\b|^k\rho\bigr)\Bigr)^{\frac{1}{k}} \Bigl( H_\e \bigl( \rho^{\frac{k}{k-1}-\frac{1}{k-1}}\bigr)\Bigr)^{\frac{k-1}{k}} 
\le \Bigl(H_\e\bigl(|\b|^k\rho\bigr)\Bigr)^{\frac{1}{k}} \rho_\e^{\frac{k-1}{k}}. 
\end{align}
Here the first inequality above follows from the H\"older inequality for the heat kernel measure $h_t(x,dy)$ with $H_tf(x)=\int_X f h_t(x, dy)$.
By \eqref{ineq: InvM-1} and \eqref{ineq: InvM-2}, it holds that 
\begin{align} \label{ineq: II}
\biggl( \int_X \psi_{\alpha-1}\circ \rho_\e |\nabla \rho_\e|^2 dm\biggr)^{1/2}
& \le e^{-K\e/2} \biggl( \int_{X}\bigl( H_\e (|\b| \rho) \bigr)^{2} \psi_{\alpha-1}\circ \rho_\e   dm \biggl)^{1/2} \notag
\\
& \le e^{-K\e/2} \biggl( \int_{X}\Bigl(H_\e(|\b|^k\rho\bigr)\Bigr)^{\frac{2}{k}} \rho_\e^{\frac{2(k-1)}{k}} \psi_{\alpha-1}\circ \rho_\e   dm \biggl)^{1/2}.
 \end{align}
 To show the convergence of the right-hand side in \eqref{ineq: II}, we prove the following statements:
 \begin{itemize}
 \item[(i)] $\bigl(H_\e(|\b|^k\rho\bigr)\bigr)^{\frac{2}{k}}$ converges to $|\b|^2\rho^{\frac{2}{k}}$ strongly in $L^{\frac{k}{2}}(X,m)$ as $\e \to 0$;
 \item[(ii)] $f_\e:=\rho_\e^{\frac{2(k-1)}{k}} \psi_{\alpha-1}\circ \rho_\e$ converges to $f:=\rho^{\frac{2(k-1)}{k}}\psi_{\alpha-1}\circ \rho$ strongly in $L^{\frac{k}{k-2}}(X,m)$ as $\e \to 0$, whereby ${\frac{k}{k-2}}$ is the conjugate exponent of $\frac{k}{2}$.
 \end{itemize}
Statement (i) can be proved by the assumption $|\b| \in L^k(X,\mu)$ and the strong continuity of $\{H_t\}$  on $L^p(X,m)$ for any $1 \le p < \infty$. 
 Thus, it suffices to show (ii).  By  $0<\alpha<1$, 
\begin{align} \label{eq: LLL}
f_\e=\rho_\e^{\frac{2(k-1)}{k}} \psi_{\alpha-1}\circ \rho_\e 
= \frac{\rho_\e^{\frac{2(k-1)}{k}}}{(\rho_\e+\eta)^{1-\alpha}}
=\frac{\rho_\e^{1-\alpha}}{(\rho_\e+\eta)^{1-\alpha}}\rho_\e^{\frac{2(k-1)}{k}-(1-\alpha)} \le \rho_\e^{\frac{2(k-1)}{k}-(1-\alpha)}.
\end{align}
Thus, we obtain that 
\begin{align} \label{eq: LLL1}
\biggl( \frac{2(k-1)}{k}-(1-\alpha)\biggr)\frac{k}{k-2} = \frac{k-2+k\alpha}{k-2} =\alpha\frac{k}{k-2}+1.
\end{align}
Noting $\rho \in L^{\alpha\frac{k}{k-2}+1}(X,m)$, we conclude that 
$$f_\e \in L^{\frac{k}{k-2}}(X,m). $$
Considering \eqref{eq: LLL}, \eqref{eq: LLL1} and the contraction of $\{H_t\}$, we have that 
\begin{align*}
\sup_{\e}\|f_\e\|_{\frac{k}{k-2}} & \le \sup_{\e}\biggl\|\rho_\e^{\frac{2(k-1)}{k}-(1-\alpha)}\biggr\|_{\frac{k}{k-2}} 
 \le \sup_{\e}\|\rho_\e\|_{\alpha\frac{k}{k-2}+1}^{\alpha\frac{k}{k-2}+1} 
 <\infty.
\end{align*}
Thus, the sequence $\{f_\e\}_\e$ is uniformly bounded in $L^{\frac{k}{k-2}}(X,m)$. Combining $\rho_\e \to \rho$ strongly in $L^{\alpha\frac{k}{k-2}+1}(X,m)$, we have proved that $f_\e \to f$ strongly in $L^{\frac{k}{k-2}}(X,m)$, which is Statement (ii). 

  Therefore, by \eqref{eq: DFo}, \eqref{ineq: II}, (i) and (ii), we have that 
   \begin{align*}
  \limsup_{\e \to 0}\frac{2}{\alpha+1}\biggl(\int_X \Bigl| \nabla \psi_{\frac{\alpha+1}{2}}\circ \rho_\e\Bigr|^2 dm\biggr)^{1/2}
  & = \limsup_{\e \to 0}\biggl( \int_X \psi_{\alpha-1}\circ \rho_\e |\nabla \rho_\e|^2 dm\biggr)^{1/2}
  \\
 & \le \biggl( \int_{X}|\b|^2\rho^{2}  \psi_{\alpha-1}\circ \rho   dm \biggl)^{1/2}
 \\
 & \le   \biggl( \int_{X}|\b|^2  \psi_{\alpha+1}\circ \rho   dm \biggl)^{1/2}
 <\infty.
 \end{align*}
Here the last two inequalities above can be seen by similar arguments as in \eqref{ineq: Sob-M1} and \eqref{eq: LLL}.
  Note that $\rho^{(\alpha+1)/2}_\e \in L^2(X,m)$ since $\alpha+1 \le \alpha\frac{k}{k-2}+1$ and $\rho \in L^{\alpha\frac{k}{k-2}+1}(X,m)$ by the assumption.  
Therefore, $\{\psi_{\frac{\alpha+1}{2}}\circ \rho_\e\}_\e$ is bounded in $W^{1,2}(X,d,m)$, which implies that there is a converging subsequence weakly in $W^{1,2}(X,d,m)$ and strongly in $L^2(X,m)$. 
Since $\rho_\e$ converges to $\rho$ strongly in $L^{\alpha\frac{k}{k-2}+1}(X,m)$ by the strong continuity of $\{H_t\}$, we have that the full sequence $\{\psi_{\frac{\alpha+1}{2}}\circ \rho_\e\}_\e$ converges to $\psi_{\frac{\alpha+1}{2}}\circ \rho$ weakly in $W^{1,2}(X,d,m)$ and strongly in $L^2(X,m)$. 
Thus, we obtain that  
\begin{align} \label{ineq: II-6-1}
 \frac{2}{\alpha+1}\biggl(\int_X \Bigl| \nabla \psi_{\frac{\alpha+1}{2}, \eta}\circ \rho\Bigr|^2 dm\biggr)^{1/2} 
 & \le \frac{2}{\alpha+1}\liminf_{\e \to 0}\biggl(\int_X \Bigl| \nabla \psi_{\frac{\alpha+1}{2}, \eta}\circ \rho_\e \Bigr|^2 dm\biggr)^{1/2} \notag
 \\
  & \le \frac{2}{\alpha+1} \limsup_{\e \to 0}\biggl(\int_X \Bigl| \nabla \psi_{\frac{\alpha+1}{2}, \eta}\circ \rho_\e \Bigr|^2 dm\biggr)^{1/2} \notag
  \\
& =  \limsup_{\e \to 0}\biggl( \int_X \psi_{\alpha-1, \eta}\circ \rho_\e |\nabla \rho_\e|^2 dm\biggr)^{1/2} \notag
\\
 & \le \biggl( \int_{X}|\b|^2\rho^{2}  \psi_{\alpha-1, \eta}\circ \rho   dm \biggl)^{1/2} \notag
 \\
 & \le \biggl( \int_{X}|\b|^2\psi_{\alpha+1, \eta}\circ \rho   dm \biggl)^{1/2} \notag
 \\
 & \to \biggl( \int_{X}|\b|^2 \rho^{\alpha+1}   dm \biggl)^{1/2} \quad (\eta \to 0).
 \end{align}
%
Since $\psi_{\alpha-1, \eta}\circ \rho \to \rho^{\alpha-1}$ $m$-a.e.\ as $\eta\to 0$, it holds that, by the Fatou lemma,
\begin{align} \label{ineq: II-6-2}
\frac{2}{\alpha+1}\biggl( \int_X \Bigl| \nabla \rho^{\frac{\alpha+1}{2}}\Bigr|^2  dm\biggr)^{1/2} 
=\biggl( \int_X \rho^{\alpha-1} |\nabla \rho|^2 dm\biggr)^{1/2} & \le \liminf_{\eta \to 0} \biggl( \int_X \psi_{\alpha-1, \eta}\circ \rho |\nabla \rho|^2 dm\biggr)^{1/2} \notag
\\
& =  \frac{2}{\alpha+1} \liminf_{\eta \to 0} \biggl(\int_X \Bigl| \nabla \psi_{\frac{\alpha+1}{2}, \eta}\circ \rho\Bigr|^2 dm\biggr)^{1/2} .
\end{align}
The two equalities above hold by applying the $S^p(X)$-version ($1<p<\infty$) of \cite[Theorem 2.2.6]{G16} (see also \cite[the third paragraph in (2) in \S 2.5]{G16}) combined with the result \cite{GH16}. In fact, by applying this chain rule, we have that $|\nabla (l\vee\rho \wedge M)^{(\alpha+1)/2}|^2=\rho^{\alpha-1} |\nabla \rho|^2\1_{\{l \le \rho \le M\}}$ $m$-a.e.\ for $0<l<M<\infty$. By $M \to \infty$, $l \to 0$ and the monotone convergence theorem, we have the first equality.
By \eqref{ineq: II-6-1} and \eqref{ineq: II-6-2}, we obtain the desired result:
 \begin{align*} 
\frac{2}{\alpha+1}\biggl( \int_X \Bigl| \nabla \rho^{\frac{\alpha+1}{2}}\Bigr|^2  dm\biggr)^{1/2} 
=\biggl( \int_X \rho^{\alpha-1} |\nabla \rho|^2 dm\biggr)^{1/2} 
\le \biggl( \int_{X}|\b|^2 \rho^{\alpha+1} dm \biggl)^{1/2}.
 \end{align*}
We finished the proof. \qed
 
By \eqref{ineq: Sob-M} and \eqref{ineq: Sob-M1}, we have the following corollary:
\begin{cor} \label{thm: Sob-M}
Assume  \eqref{eq: invL5-1} and that $(X,d,m)$ is an RCD$(K,N)$ space with $2<N<\infty$ and $K>0$. 
Let $|\b| \in L^{k}(X,\mu)$ for some $k>2$. If $\rho \in L^{\alpha\frac{k}{k-2}+1}(X,m)$ for some $\alpha>0$, then $\rho^{(\alpha+1)/2} \in W^{1,2}(X,d,m)$ and  
\begin{align} \label{ineq: Sob-M2}
\int_X |\nabla \rho^{\frac{\alpha+1}{2}}|^2 dm \le \biggl(\frac{\alpha+1}{2}\biggr)^2\biggl( \int_X |\b|^k \rho dm \biggr)^{2/k} \biggl( \int_X \rho^{\alpha\frac{k}{k-2}+1} dm\biggr)^{(k-2)/k}.
\end{align}
\end{cor}

We now show the integrability of $\rho$ in Theorem \ref{thm: mthm1-2}. 
\begin{lem} \label{thm: Lp}
Assume \eqref{eq: invL5-1} and that $(X,d,m)$ is an RCD$(K,N)$ space with $2<N<\infty$ and $K>0$. Then the following statements hold:
\begin{itemize}
\item[(i)] 
if $|\b| \in L^{k}(X,\mu)$ with $k>N$, then $\rho \in L^\infty(X,m)$; 
\item[(ii)] 
if $|\b| \in L^{N}(X,\mu)$, then $\rho \in L^p(X,m)$ for any $1 \le p <\infty$.
\end{itemize}
\end{lem}
\proof
Assume that $\rho \in L^{\alpha\frac{k}{k-2}+1}(X,m)$ for some $\alpha>0$.
By \eqref{ineq: Sob-M2} and the Sobolev inequality \eqref{ineq: SL}, we have that $\rho \in L^{(\alpha+1)\frac{N}{N-2}}(X,m)$ and 
\begin{align} \label{ineq: imp}
\biggl( \int_X \rho^{(\alpha+1)\frac{N}{N-2}} dm \biggr)^{(N-2)/N}
& \le C_S\|\nabla \rho^{(\alpha+1)/2}\|^2_2+\|\rho^{(\alpha+1)/2}\|^2_2 \notag
\\
& \le C_S \|\b\|^2_{L^k(X,\mu)}\biggl(\frac{\alpha+1}{2}\biggr)^2\biggl( \int_X \rho^{\frac{\alpha k}{k-2}+1} dm\biggr)^{(k-2)/k} +\|\rho^{(\alpha+1)/2}\|^2_2. 
\end{align}
Let $C=C_S \|\b\|^2_{L^k(X,\mu)}$. Noting $\frac{\alpha k}{k-2}+1 \ge \alpha+1$, we see that the right-hand side of \eqref{ineq: imp} is dominated by the following quantity:
\begin{align*}
C\biggl(\frac{\alpha+1}{2}\biggr)^2\|\rho\|_{\frac{\alpha k}{k-2}+1}^{\alpha+\frac{k-2}{k}}+\|\rho\|^{\alpha+1}_{\alpha+1}
\le C\biggl(\frac{\alpha+1}{2}\biggr)^2\|\rho\|_{\frac{\alpha k}{k-2}+1}^{\alpha+\frac{k-2}{k}}+\|\rho\|^{\alpha+1}_{\frac{\alpha k}{k-2}+1}.
\end{align*}
Combining this with \eqref{ineq: imp}, we obtain that 
\begin{align} \label{ineq: imp2}
\|\rho\|_{\frac{(\alpha+1)N}{N-2}}^{{\alpha+1}}
\le 
C\biggl(\frac{\alpha+1}{2}\biggr)^2\|\rho\|_{\frac{\alpha k}{k-2}+1}^{\alpha+\frac{k-2}{k}}+\|\rho\|^{\alpha+1}_{\frac{\alpha k}{k-2}+1}.
\end{align}
The inequality \eqref{ineq: imp2} improves the integrability of $\rho$ inductively. Let $\eta_n$ be the integrability of $\rho$ in the $n$th step, i.e., $\rho \in L^{\eta_n}(X,m)$. Then, we have the following relation:
\begin{align} \label{ineq: imp1}
\eta_{n+1}=r\Bigl(\eta_n+\frac{2}{k-2}\Bigr), \quad r=\frac{N}{N-2}\frac{k-2}{k}.
\end{align}
By the assumption $k>N$, it holds that $r>1$. Since $\sqrt{\rho} \in W^{1,2}(X, d, m)$ by Lemma \ref{lem: W12-1}, by using the Sobolev inequality \eqref{ineq: SL}, it holds that $\rho \in L^{N/(N-2)}(X,m)$. Thus, we can take $\eta_0=N/(N-2)$ and $\alpha_n>0$ so that 
$$\alpha_n+1=\frac{N-2}{N} \eta_{n+1}=\frac{k-2}{k}\Bigl( {\eta_n+\frac{2}{k-2}}\Bigr).$$
We note that $\eta_n \to \infty$ and $\alpha_n \to \infty$ as $n \to \infty$. 
We may assume that $\|\rho\|_{\eta_n} \ge 1$ for sufficiently large $n$ since if not, $\|\rho\|_{\eta_n} \le 1$ for any $n \in \N$, which implies $\|\rho\|_\infty \le 1$ and concludes the desired result. We may also assume that 
$$C\biggl(\frac{\alpha_n+1}{2}\biggr)^2 \ge 1,$$
for sufficiently large $n$.
Noting that, for any $\alpha>0$ and $\|\rho\|_{\frac{\alpha k}{k-2}+1} \ge 1$, 
$$\|\rho\|_{\frac{\alpha k}{k-2}+1}^{\alpha+\frac{k-2}{k}} \le \|\rho\|^{\alpha+1}_{\frac{\alpha k}{k-2}+1},$$
it holds that, for sufficiently large $n$, 
\begin{align} \label{ineq: imp3}
\|\rho\|_{\eta_{n+1}}^{{\alpha_n+1}}
 \le 
C\biggl(\frac{\alpha_n+1}{2}\biggr)^2\|\rho\|_{\eta_n}^{\alpha_n+\frac{k-2}{k}}+\|\rho\|^{\alpha_n+1}_{\eta_n} \notag
& \le 
C\biggl(\frac{\alpha_n+1}{2}\biggr)^2\|\rho\|_{\eta_n}^{\alpha_n+1}+C\biggl(\frac{\alpha_n+1}{2}\biggr)^2\|\rho\|^{\alpha_n+1}_{\eta_n} \notag
\\
& \le 
2C\biggl(\frac{\alpha_n+1}{2}\biggr)^2\|\rho\|_{\eta_n}^{\alpha_n+1}.
\end{align}
Therefore, we have that, with some constant $C_1>0$,  
$$\|\rho\|_{\eta_{n+1}} \le \biggl(C_1\frac{\alpha_n+1}{2}\biggr)^{\frac{2}{\alpha_n+1}}\|\rho\|_{\eta_n}.$$
Taking the logarithm in both sides of the above inequality with $\theta_n=\log \|\rho\|_{\eta_n}$, it can be seen that 
$$\theta_{n+1} \le \frac{2}{\alpha_n+1} \log \biggl\{C_1\frac{\alpha_n+1}{2}\biggr\}+\theta_{n}.$$
We have that 
\begin{align} \label{ineq: imp6}
\theta_{n+1}-\theta_n \le \frac{2}{\alpha_n+1} \log \biggl\{C_1\frac{\alpha_n+1}{2}\biggr\}  
\le \frac{2N}{\eta_{n+1}(N-2)} \log \biggl\{\frac{C_1(N-2)}{2N}\eta_{n+1}\biggr\}  
 \le \frac{C_2}{\eta_{n+1}^{1-\e}},
\end{align}
for any $0<\e<1$ with some constant $C_2>0$. Since the right-hand side of \eqref{ineq: imp6} is summarable in $n$, it concludes that $\theta_n \to \theta=\log \|\rho\|_\infty<\infty$ as $n \to \infty.$ 
We finished the proof of (i). 
%

In the case of $k=N$ in (ii), by \eqref{ineq: imp1}, it holds that 
$$\eta_{n+1}=\eta_n+\frac{2}{N-2}.$$ Thus, $\eta_n \uparrow \infty$ and we have finished the proof.  

 \qed
 
By Lemma \ref{thm: Lp}, we finished the proof of the integrability of $\rho$ in Theorem \ref{thm: mthm1-2}.

We now give the proof of \eqref{ineq: GER6} in Theorem \ref{thm: mthm1-2}.
\\
{\it Proof of \eqref{ineq: GER6} in Theorem \ref{thm: mthm1-2}}. 
Let $k=N$, and $\alpha=1$ in \eqref{ineq: Sob-M2}, and  $\alpha=(N-2)/N$ in \eqref{ineq: imp}. Then, we have that 
\begin{align} \label{ineq: iii}
\int_X |\nabla \rho|^2 dm & \le \int_X |\b|^2\rho^2 dm  
\le \biggl( \int_X |\b|^N \rho dm \biggr)^{2/N} \biggl( \int_X \rho^{\frac{N}{N-2}+1} dm\biggr)^{(N-2)/N} 
 <\infty.
\end{align}
Thus, $\rho \in W^{1,2}(X,d,m)$ and \eqref{ineq: GER6} has been proved.  \qed

The following corollary will be used in Section \ref{sec: SIM} in order to show the stability of invariant measures. 
\begin{cor} \label{cor: grad}
Assume \eqref{eq: invL5-1} and that $(X,d,m)$ is an RCD$(K,N)$ space with $2<N<\infty$ and $K>0$ with $|\b| \in L^N(X,\mu) \cap L^{2N}(X,m)$. Then it holds that 
$$\|\nabla \rho\|_2^2 \le  C_S \biggl(\frac{N-1}{N}\biggr)^2\|\b\|_{2N}^{4}\|\rho\|_2^{2} +\|\b\|_{2N}^{2} \|\rho\|_2^{2}.$$
\end{cor}
\proof
Since $2(N-1)/N \le 2$ and $\|f\|_p \le \|f\|_q$ for any $p \le q$ owing to $m(X)=1$, we have that 
$$\|\rho\|_{\frac{2(N-1)}{N}}^{\frac{2(N-1)}{N}} \le \|\rho\|_{2}^{\frac{2(N-1)}{N}}.$$
Therefore, by the Sobolev inequality \eqref{ineq: SL},  the right-hand side of \eqref{ineq: iii} can be estimated from above as follows:
\begin{align*}
& \biggl( \int_X |\b|^{N} \rho dm \biggr)^{2/N} \Biggl( C_S  \biggl( \int_X |\b|^N \rho dm \biggr)^{2/N}\biggl(\frac{N-1}{N}\biggr)^2\biggl( \int_X \rho^2 dm\biggr)^{(N-2)/N} +\|\rho^{\frac{N-1}{N}}\|^2_2 \Biggr) \notag
\\
&\le  C_S  \biggl(\frac{N-1}{N}\biggr)^2\biggl( \int_X |\b|^N \rho dm \biggr)^{4/N}\biggl( \int_X \rho^2 dm\biggr)^{(N-2)/N} +\biggl( \int_X |\b|^N \rho dm \biggr)^{2/N}\|\rho^{\frac{N-1}{N}}\|^2_2 \notag
 \\
 & \le  C_S \biggl(\frac{N-1}{N}\biggr)^2 \biggl( \int_X |\b|^N \rho dm \biggr)^{4/N}\biggl( \int_X \rho^2 dm\biggr)^{(N-2)/N} +\biggl( \int_X |\b|^N \rho dm \biggr)^{2/N}\|\rho\|_{2}^{\frac{2(N-1)}{N}}
  \\
 & \le  C_S \biggl(\frac{N-1}{N}\biggr)^2\|\b\|_{2N}^{4}\|\rho\|_2^{\frac{4}{N}}\|\rho\|_2^{\frac{2(N-2)}{N}} +\|\b\|_{2N}^{2} \|\rho\|_2^{\frac{2}{N}}\|\rho\|_{2}^{\frac{2(N-1)}{N}}
   \\
 & =  C_S \biggl(\frac{N-1}{N}\biggr)^2\|\b\|_{2N}^{4}\|\rho\|_2^{2} +\|\b\|_{2N}^{2} \|\rho\|_2^{2}.
\end{align*}
Therefore, we conclude the desired result.
\qed

 \section{Symmetry of Semigroups} \label{sec: SS}
 In this section, we characterize the symmetrizability of the semigroup $\{T_t\}$. 
 Throughout this section (except Proposition \ref{lem: Dom1}), Condition (B) in Section \ref{sec: EU}, $m(X)=1$ and $|\b| \in L^\infty(X,m)$ are assumed.
Let $(\E, \F)$ be the Dirichlet form defined in \eqref{eq: DF}.
Let $(L, \mathcal D(L))$ and $\{T_t\}$ be the infinitesimal generator and the semigroup on $L^2(X,m)$ corresponding to $(\E, \F)$ respectively. By Theorem \ref{thm: mthm2} and Theorem \ref{thm: mthm1}, there exists a unique (up to constant multiplication) invariant measure $\mu=\rho m$ for $(L, \mathcal D(L))$ so that $\rho$ is strictly positive and $\rho \in W^{1,2}(X,d,m)$.

Before the main result in this section, we prove the essential self-adjointness of $(\Delta, {\rm TestF}(X))$ in $L^2(X, m)$.
\begin{prop} \label{lem: Dom1}
Let $(X,d,m)$ be an RCD$(K,\infty)$ space with $K \in \R$.  Then $(\Delta, {\rm TestF}(X))$ is essentially self-adjoint in $L^2(X, m)$, i.e., 
${\rm TestF}(X)$ is dense in $\mathcal D(\Delta)$ with respect to the graph norm $\|\cdot\|_\Delta^2=\|\cdot\|^2_2+\|\Delta \cdot\|^2_2$.
\end{prop}
\proof
By using the Bakry-\'Emery gradient estimate \eqref{ineq: GE}, it is straightforward to check that the heat semigroup action preserves {\rm TestF}(X), i.e., $H_tf \in {\rm TestF}(X)$ whenever $f \in {\rm TestF}(X)$ for any $t \ge 0$. Thus, by \cite[Theorem X.49]{RS75}, we can conclude the statement. 
\qed

The following result is the main theorem in this section. 
\begin{thm} \label{thm: mthm3}
Assume (B), $m(X)=1$ and
$|\b| \in L^\infty(X,m)$.  
Then $\{T_t\}$ is symmetrizable if, and only if, there exists $f \in {\rm Lip}(X) \cap L^2(X,m)$ so that $$\b=\langle \nabla f, \nabla \cdot\rangle.$$ The symmetrizing measure $\mu$ is unique up to constant multiplication, and $\mu=\rho dm$ with $\rho=ce^{2f}$ $m$-a.e.\ with some constant $c>0$. 
\end{thm}
\begin{rem} \normalfont  
 In the case of the Wiener space, the set of cylinder functions  can be taken as a dense subset in $\mathcal D(\Delta)$ (see Shigekawa \cite[Theorem 2.1]{Shi87}).
 By replacing Lipschitz functions with $H$-Lipschitz functions in the statement of Theorem \ref{thm: mthm3}, we can apply the same proof of Theorem \ref{thm: mthm3} to the case of the Wiener space.
\end{rem}

Hereinafter, we prove Theorem \ref{thm: mthm3}.
\begin{lem} \label{lem: Dom}
${\rm TestF}(X)$ is dense in $\mathcal D(L)$ with respect to the graph norm $\|\cdot\|^2_L=\|\cdot\|^2_2+\|L \cdot\|_2^2$.
\end{lem}
\proof  
By Proposition \ref{lem: Dom1}, for any $f \in \mathcal D(L)=\mathcal D(\Delta)$ (by Proposition \ref{prop: domain}), we can take $\{f_i\} \subset {\rm TestF}(X)$ so that $f_i$ converges to $f$ in the graph norm $\|\cdot\|_{\Delta}$. Since $|\b| \in L^\infty(X,m)$, in order to prove that $f_i$ converges to $f$ with respect to $\|\cdot\|_{L}$, it suffices to show that $f_i$ converges to $f$ in $W^{1,2}(X,d,m)$. Since $\Delta f_i \to \Delta f$ and $f_i \to f$ both in $L^2$, we have that 
 \begin{align*}
 \int_X |\nabla (f_i-f)|^2 dm= -\int_X \Delta (f_i-f) (f_i-f) dm   \to 0.
 \end{align*}
We finished the proof. \qed

\begin{lem} \label{lem: symE}
Let $f \in W^{1,2}(X,d,m)$ satisfy $|\nabla f| \in L^\infty(X,m)$. Then $e^{f} \in W^{1,2}(X,d,m)$ and $\nabla e^{f}=e^f\nabla f$.
\end{lem}
\proof
Take a smooth function $\phi: \R \to \R$ so that $\phi(z)=z$ for $z \le 0$, and $\phi(z)=1$ for $z \ge 1$. Let $\phi_n(z)=n+\phi(z-n)$.
Set $f_n(x)=\phi_n(f(x))$. Then $f_n \in W^{1,2}(X,d,m)$ and 
$$|\nabla f_n|=(\phi_n'\circ f)|\nabla f|, \quad \text{$m$-a.e.}$$
Thus, $|\nabla f_n| \in L^\infty(X,m)$ and 
\begin{align} \label{ineq: sym}
\|\nabla f_n\|_\infty \le \|\phi'_n\|_\infty\|\nabla f\|_\infty \le \|\phi'\|_\infty\|\nabla f\|_\infty.
\end{align}
Let $g_n=c_ne^{f_n}$ with the constant $c_n$ determined by $\|g_n\|_2=1$. Then $g_n \in W^{1,2}(X,d,m)$ and $\nabla g_n=g_n \nabla f_n$. 
By \eqref{ineq: sym}, it holds that 
$$\|\nabla g_n\|_2 \le \|g_n\|_2\|\nabla f_n\|_\infty \le \|\phi'\|_\infty\|\nabla f\|_\infty.$$
Hence, $\{g_n\}_{n \in \N}$ is a bounded sequence in $W^{1,2}(X,d,m)$. By \cite[Theorem 6.3, Proposition 6.7]{GMS13}, a converging subsequence $\{g_{n'}\}_{n' \in \N}$ can be taken with respect to the $L^2$-strong convergence and the $W^{1,2}$-weak convergence. 
Noting that $c_n$ is non-increasing, let $\lim_{n \to \infty}c_n=c$. Then it can be seen that $g_n$ converges to $g=ce^f$ $m$-a.e. Hence, $g_n$ converges to $g$ strongly in $L^2$ and weakly in $W^{1,2}$. Therefore, $g \in W^{1,2}(X,d,m)$, and noting $\|g_n\|_2=1$ for any $n \in \N$, it can be seen that $c>0$. It is easy to check that $\nabla e^{f}=e^f\nabla f$.
\qed
\begin{lem} \label{lem: W12}
$\log \rho \in W^{1,2}(X,d,m)$ and 
\begin{align} \label{ineq: GER2}
\|\nabla \log \rho\|_{2} \le \|\b\|_\infty.
\end{align}
\end{lem}
\proof
Let $\rho_\e=H_\e \rho$ and $\phi_\e$ be defined as follows: 
$$\phi_\e=\frac{1}{\rho_\e+\delta}-\frac{1}{\delta} \in W^{1,2}(X,d,m) \quad \delta>0.$$
Then, by \eqref{eq: Inv} with the above $\phi=\phi_\e$ and the gradient estimate \eqref{ineq: GE}, it can be seen that
\begin{align*}
\int_X \frac{|\nabla \rho_\e|^2}{(\rho_\e+\delta)^2} dm 
= \int_{X}\b(H_\e \phi_\e ) d\mu
 & \le  \int_{X}|\b| | \nabla H_\e \phi_\e | d\mu
 \\
& \le  e^{-(K/2) \e}\|\b\|_\infty \int_{X} \frac{|\nabla \rho_\e|}{(\rho_\e+\delta)} \frac{H_\e (\rho)}{(\rho_\e+\delta)} dm
\\
& \le  e^{-(K/2) \e}\|\b\|_\infty \Bigl( \int_{X} \frac{|\nabla \rho_\e|^2}{(\rho_\e+\delta)^2} dm \Bigr)^{1/2}\Bigl(\int_X \frac{\rho_\e^2}{(\rho_\e+\delta)^2} dm\Bigr)^{1/2}.
\end{align*}
Therefore, the following holds:
\begin{align*}
\Bigl( \int_X \frac{|\nabla \rho_\e|^2}{(\rho_\e+\delta)^2} dm \Bigr)^{1/2} \le e^{-(K/2) \e}\|\b\|_\infty \Bigl(\int_X \frac{\rho_\e^2}{(\rho_\e+\delta)^2} dm\Bigr)^{1/2}
\le e^{-(K/2) \e}\|\b\|_\infty.
\end{align*}
Thus, it holds that 
\begin{align*}
\limsup_{\e \to 0}\Bigl( \int_X \frac{|\nabla \rho_\e|^2}{(\rho_\e+\delta)^2} dm \Bigr)^{1/2} \le \|\b\|_\infty.
\end{align*}
Hence, $\{\phi_\e\}_\e$ is a bounded sequence in $W^{1,2}(X,d, m)$. 
Therefore, there is a converging subsequence $\{\phi_{\e'}\}_{\e'}$ and its limit $\phi' \in W^{1,2}(X,d,m)$ in the $L^2$-strong sense and the $W^{1,2}$-weak sense. Since $\phi_\e$ converges to $\phi=1/(\rho+\delta)+1/\delta$ in the $L^2$-strong sense, it holds that $\phi=\phi'$.
Thus, $\phi_\e$ converges to $\phi$ in the $L^2$-strong sense and the $W^{1,2}$-weak sense. Therefore, the following inequality is obtained:
\begin{align*}
\Bigl( \int_X \frac{|\nabla \rho|^2}{(\rho+\delta)^2} dm \Bigr)^{1/2} \le \liminf_{\e \to 0}\Bigl( \int_X \frac{|\nabla \rho_\e|^2}{(\rho_\e+\delta)^2} dm \Bigr)^{1/2} 
 \le \limsup_{\e \to 0}\Bigl( \int_X \frac{|\nabla \rho_\e|^2}{(\rho_\e+\delta)^2} dm \Bigr)^{1/2} 
 \le \|\b\|_\infty.
\end{align*}
As $\delta$ tends to zero with a similar argument to \eqref{ineq: II-6-1} and \eqref{ineq: II-6-2}, the desired inequality \eqref{ineq: GER2} is obtained. 


Now it is proved that $\log \rho \in W^{1,2}(X,d,m)$, for which it suffices to show that $\log \rho \in L^2(X,m)$ according to \eqref{ineq: GER2}. Under (B), the Poincar\'e inequality \eqref{ineq: PI} holds:  
\begin{align*}
\int_X|f-m(f)|^2dm \le C_P \int_X |\nabla f|^2dm, \quad f \in W^{1,2}(X,d,m).
\end{align*}
Let $f_n=| \log (\min\{\rho, n\}+n^{-1})|.$ According to \eqref{ineq: GER2}, it can be seen that $|\nabla f_n|$ is uniformly bounded in $L^2(X,m)$. By combining this with the Poincar\'e inequality, it holds that $f_n-m(f_n)$ is uniformly bounded in $L^2(X,m)$. Noting that $\rho$ is strictly positive by Theorem \ref{thm: mthm2}, it holds that $f_n \to |\log \rho|$ $m$-a.e. 
By using Fatou's lemma, the sequence $\{m(f_n)\}_{n}$ is bounded. In fact, if $\{m(f_n)\}_{n}$ is unbounded, a subsequence $\{m(f_{n'})\}_{n'}$ can be taken so that 
$|m(f_{n'})| \to \infty. $
By Fatou's lemma, it holds that
\begin{align*}
\infty=\int_X\lim_{n' \to \infty}|f_{n'}-m(f_{n'})|^2dm 
\le \liminf_{n' \to \infty}\int_X|f_{n'}-m(f_{n'})|^2dm.
\end{align*}
This contradicts the Poincar\'e inequality and the $L^2$-uniform boundedness of $|\nabla f_n|$. Thus $\{m(f_n)\}_{n}$ is bounded. Hence, $\|f_n\|_2$ is uniformly bounded, which implies $|\log \rho| \in L^2(X,m)$. The proof has been completed. 
\qed

Now the proof of Theorem \ref{thm: mthm3} is presented.
\\
{\it Proof of Theorem \ref{thm: mthm3}.}  
The ``only if" part is proved first. Assume that $\{T_t\}$ is symmetric with respect to a measure $\mu$. Then, $\mu$ is an invariant measure and unique up to constant multiplication by Theorem \ref{thm: mthm2}. By Theorem $\ref{thm: mthm1}$, we can write $\mu=\rho m$ with $\rho \in W^{1,2}(X,d,m)$. Let $\phi, \psi \in {\rm TestF}(X)$. By Proposition \ref{prop: domain}, it holds that ${\rm TestF}(X) \subset \mathcal D(L_2)$ and 
$$\int_X(T_t \phi-\phi)\psi d\mu=\int_X \phi (T_t \psi-\psi) d\mu.$$
Thus, by taking the time derivative $d/dt$ on both sides,  it can be obtained that 
$$\int_X \frac{d}{dt}(T_t\phi-\phi) \psi d\mu=\int_X \phi \frac{d}{dt}(T_t\psi-\psi)d\mu.$$
Letting $t \to 0$,  the following equality holds:
\begin{align} \label{eq: inv1}
\int_X L(\phi) \psi d\mu=\int_X \phi L(\psi) d\mu.
\end{align}
By the invariance of $\mu$, it holds that 
$$\int_X L(\phi \psi) d\mu=0.$$
By the Leibniz rule of $\Delta$ and $\b$, it can be seen that
$L(\phi \psi)=L(\phi)\psi+\phi L(\psi)+\langle \nabla \phi, \nabla \psi \rangle$, which implies that
\begin{align}
\int_X L(\phi \psi) d\mu=\int_X \Bigl(L(\phi) \psi + \phi L\psi +\langle \nabla \phi, \nabla \psi \rangle \Bigr)d\mu=0. \label{eq: inv-sym1}
\end{align}
Therefore, by \eqref{eq: inv1} and the Sobolev regularity of $\rho$ obtained in Theorem \ref{thm: mthm1}, it holds that 
\begin{align}
0=\int_X \Bigl(L(\phi) \psi + \phi L(\psi) +\langle \nabla \phi, \nabla \psi \rangle \Bigr)d\mu 
&=\int_X \Bigl(\phi L(\psi) + \phi L(\psi) +\langle \nabla \phi, \nabla \psi \rangle \Bigr)d\mu \notag
\\
&=\int_X \Bigl(\phi \Delta(\psi) + 2\phi \b(\psi)  +\langle \nabla \phi, \nabla \psi \rangle \Bigr)d\mu  \notag
\\
&=-\int_X \langle \nabla (\phi \rho), \nabla \psi\rangle dm + \int_X2\phi \b(\psi) \rho m  + \int_X\langle \nabla \phi, \nabla \psi \rangle \rho dm  \notag
\\
&=-\int_X\Bigl( \phi \langle \nabla \rho, \nabla \psi \rangle +\rho \langle \nabla \phi, \nabla \psi \rangle \Bigr)dm + \int_X2\phi \b(\psi) \rho m  \notag
\\
& \quad +\int_X \langle \nabla \phi, \nabla \psi \rangle \rho dm  \notag
\\
&=\int_X  \Bigl(2\rho \b(\psi) - \langle \nabla \rho, \nabla \psi\rangle \Bigr)\phi dm. \label{eq: inv-sym2}
\end{align}
Since ${\rm TestF}(X)$ is dense in $W^{1,2}(X,d, m)$, the following equality is obtained:
$$\langle \nabla \rho, \nabla v \rangle=2\rho \b(v), \quad \forall v \in W^{1,2}(X,d,m).$$
There is a unique vector field $\b^* \in L^2(TX)$ so that $\b(v)=\langle \b^*, \nabla v \rangle$ for any $v \in W^{1,2}(X,d,m)$ (see Subsection \ref{subsec: Der}). 
Since $\{\nabla v: v \in W^{1,2}(X,d,m)\}$ is dense in $L^2(TX)$ (see Subsection \ref{subsec: TFTV}), it holds that
\begin{align} \label{eq: GRI}
\nabla \rho=2\rho\b^*.
\end{align}
Take $f=(1/2)\log \rho$. Then, by Lemma \ref{lem: W12}, it holds that $f \in W^{1,2}(X,d,m)$ and, by \eqref{eq: GRI},
$$\nabla f=\frac{\nabla \rho}{2\rho}=\b^*.$$
Since $|\nabla f|=|\b| \in L^\infty(X,m)$ by the assumption, using the Sobolev-to-Lipschitz property (\cite[Theorem 6.2]{AGS14b}), we have that $f$ can be taken as a Lipschitz continuous function.
The proof of the ``only if" part has been completed.

The ``if" part is now shown. Assume that there is a function $f \in {\rm Lip}(X) \cap L^2(X,m)$ so that $\b=\langle \nabla f, \nabla \rangle$, in other words, $\b^*=\nabla f$. 
Then the measure $\rho'=e^{2f}m$ is an invariant measure for $(L, \mathcal D(L))$. In fact, 
by Lemma \ref{lem: symE}, it holds that $\rho' \in W^{1,2}(X,d,m)$ and 
$$\nabla \rho'=2e^{2f}\nabla f=2\rho' \b^*.$$
Noting $ \mathcal D(L)=\mathcal D(\Delta)$ by Proposition \ref{prop: domain}, the following holds: for $\phi \in \mathcal D(L)$,
\begin{align*}
\int_X L (\phi) \rho' dm=\frac{1}{2}\int_X \Delta(\phi) \rho' dm+\int_X \b(\phi) \rho' dm
&=-\frac{1}{2}\int_X \langle \nabla \phi, \nabla \rho' \rangle dm+\int_X \b(\phi) \rho' dm
\\
&=-\int_X \b(\phi) \rho' dm+\int_X \b(\phi) \rho' dm=0.
\end{align*}
Therefore, $\rho'$ is an invariant measure. Since $(L, \mathcal D(L))$ has a unique invariant measure up to constant multiplication by Theorem \ref{thm: mthm2}, it is obtained that $\rho'=c\rho$ for some constant $c>0$. We may assume that $c=1$ for simplicity of the following argument. 

Now it is shown that $\{T_t\}$ is symmetric with respect to $\mu=\rho m$. 
Since $\rho=e^{2f}$ and $\b^*=\nabla f$, the right-hand side of \eqref{eq: inv-sym2} is equal to zero. 
By following the arguments in \eqref{eq: inv-sym2} conversely and combining these with \eqref{eq: inv-sym1}, for any $\phi,\psi \in {\rm TestF}(X)$, it holds that 
\begin{align} \label{eq: inv-sym3}
\int_X L (\phi) \psi d\mu=\int_X \phi L(\psi) d\mu.
\end{align}

Since ${\rm TestF}(X)$ is dense in $\mathcal D(L)$ with respect to the graph norm $\|\cdot\|_L$ by Lemma \ref{lem: Dom}, it holds that, 
for any $\phi, \psi \in{\rm TestF}(X)$,
$$\int_X LT_t(\phi) T_s(\psi) d\mu =\int_X T_t(\phi) LT_s(\psi)d\mu.$$
Define
$$g(s)=\int_X T_{t-s}(\phi)T_s(\psi) d\mu, \quad 0 \le s \le t.$$
Then, by differentiating the both sides in $s$, it can be seen that
$$g'(s)=-\int_X LT_{t-s}(\phi)T_s(\psi) d\mu+\int_X T_{t-s}(\phi)LT_s (\psi) d\mu=0.$$
Thus, $g(t)=g(0)$, and 
\begin{align} \label{eq: symm}
\int_X T_t (\phi)\psi d\mu=\int_X \phi T_t (\psi)d\mu.
\end{align}
By the density of ${\rm TestF}(X)$ in $L^2(X,m)$, the equality \eqref{eq: symm} holds for any $\phi, \psi \in \mathcal B_b(X)$. The proof has been completed.
\qed

\section{Stability of Invariant Measures} \label{sec: SIM}
In this section, the stability of invariant measures with respect to the perturbation of $\b$ and $X$ is proved.
Let $(X_n, d_n, m_n, \x_n)$ be a metric measure space with a fixed point $\x_n \in X_n$ for $n \in \EN=\N\cup\{\infty\}$. The following three conditions are considered: 
\begin{itemize}
\item[(A')] {Set $K \in \R$ and $0<D<\infty$}. Let $(X_n, d_n, m_n)$ be an RCD$(K,\infty)$ space with $K>0$, or $\sup_{n \in \N}{\rm diam}(X_n)<D$.  Let $(X_n, d_n, m_n)$ be an RCD$(K,\infty)$ space and $\b_n$ be an $L^\infty$-derivation, 
and 
$$\sup_{n \in \N}\|\b_n\|_\infty <\infty.$$
\item[(B')] {Set $K \in \R$ and $0<D<\infty$}. Let $(X_n, d_n, m_n)$ be an RCD$(K,\infty)$ space with $K>0$, or $\sup_{n \in \N}{\rm diam}(X_n)<D$.
Let $\b_n$ be a derivation so that there is a constant $\theta >2C_{LS}$ independent of $n$ satisfying the following condition:
$$\sup_{n \in \N}\Bigl\| e^{\theta |\b_n|^2}\Bigr\|_1<\infty, 
$$
whereby $C_{LS}$ denotes the constant in the log-Sobolev inequality \eqref{ineq: LS}. 
\item[(C')]  Set $K>0$, $2<N<\infty$ and $c \ge 0$.  Let $(X_n, d_n, m_n)$ be an RCD$(K,N)$ space for $n \in \N$. Assume that $\sup_{n \in \N}\|\b_n\|_{2N} <\infty$ and ${\rm div} \b_n \ge -c$.
\end{itemize}
Under each of these assumptions, there is a unique non-zero invariant measure $\mu_n=\rho_nm_n$ with $\rho_n \in L^2(X_n,m_n)$ by Theorem \ref{thm: mthm2}.
For $n \in \EN$, let $\mu_n$ be the invariant measure for $(L_n, \mathcal D(L_n))$ with $\|\rho_n\|_2=1$ (by the normalization $\rho_n/\|\rho_n\|_2$, we can always take such $\rho_n$).  
Let ${\rm INV}_{>0}(L)$ and ${\rm Inv}_{>0}(L)$ denote the set of non-zero invariant measures and the set of densities of non-zero invariant measures for $(L, \mathcal D(L))$ respectively. The set of all finite Borel measures on $X$ is denoted by $\mathcal M(X)$. 
The following theorem shows the stability of the invariant measure $\mu_n$ and the density $\rho_n$ with respect to the pmG convergence of $X_n$ and the $L^2$-convergence of the derivation $\b_n$. 
\begin{thm}\label{thm: mthm4}
Suppose that $(X_n, d_n, m_n, \x_n)$ converges to $(X_\infty, d_\infty, m_\infty, \x_\infty)$ in the sense of the pmG, and $\b_n \to \b_\infty$ and ${\rm div}\b_n \to {\rm div}\b_\infty$ strongly in $L^2$ respectively.
Let $(X,d)$ be a common ambient space for $\{X_n\}_{n \in \N}$ associated with the pmG convergence (see Subsection \ref{subsec: D}).
Then the following statements hold:
\begin{itemize}
\item[(i)] Under (A') or (C'), $\rho_n$ converges to $\rho_\infty$ weakly in $W^{1,2}$ and strongly in $L^2$;
\item[(ii)]  under (B'), let $A$ denote the set of limit points of $\{\rho_n\}_{n}$ with respect to the weak $L^2$-sense.  Then $A \subset {\rm Inv}_{>0}(L_\infty)$.


%
\end{itemize} 
In particular, $\mu_n \to \mu_\infty$ weakly in $\mathcal M(X)$ in the cases of (i). In the case of (ii), let $\mathcal A$ denote the set of limit points of $\{\mu_n\}$ with respect to the weak convergence. Then $\mathcal A \subset {\rm INV}_{>0}(L_\infty)$.
\end{thm}
\begin{rem} \normalfont  
We give two remarks on Theorem \ref{thm: mthm4}. 
\label{rem: mthm4}
\begin{itemize}
\item[(i)]
If the underlying spaces are fixed, i.e., $X_n=X_\infty$ for any $n \in \N$, then the result (i) in Theorem \ref{thm: mthm4} gives the convergence of $\mu_n$ in the total variation distance $d_{TV}$: 
\begin{align*} 
d_{TV}(\mu_n, \mu_\infty)= \sup_{|f| \le 1}|\mu_n(f)-\mu_\infty(f)|.
\end{align*}
\item[(ii)] Theorem \ref{thm: mthm4} tells us that the weak convergence of invariant measures can be derived from weaker assumptions compared to the convergence in law of diffusion processes. Indeed, in order to show the convergence in law of the corresponding diffusion processes, we need the boundedness of $|\b|$ for the tightness of the laws (see \cite[Theorem 1.3]{S18}, and also \cite[Theorem 11.1.4]{SV79} in the case of the Euclidean space). However, for the weak convergence of invariant measures, the $L^{2N}$-integrability is sufficient in the case of (C'). 
\end{itemize} 
\end{rem}

We now start to prove Theorem \ref{thm: mthm4}. Under Condition (A'), (B') or (C'), the bilinear form $(\E_n, \F_n)$ defined in \eqref{eq: DF} is a closed form for any $n \in \EN$ and let $\{T_t^n\}$ be the corresponding semigroup. The $L^2$-convergence of the semigroups $\{T_t^n\}$ is proved first.  
\begin{prop} \label{lem: Hino_Conv}
Assume that (A'), (B') or (C') hold. If 
\begin{itemize}
\item[(i)] $X_n$ converges to $X_\infty$ in the pmG sense;
\item[(ii)] $\b_n \to \b_\infty$ and ${\rm div}\b_n \to {\rm div}\b_\infty$ strongly in $L^2$ respectivly.
\end{itemize}
Then the semigroup $T_t^n$ converges to $T^\infty_t$ in $L^2$ for all $t \ge 0$. Namely, $T^n_tf_n$ converges to $T^\infty_tf_\infty$ strongly in $L^2$ for any $t \ge 0$ whenever $f_n$ converges to $f_\infty$ strongly in $L^2$.
\end{prop}
\proof
The coercivity constant $\lambda$ can be taken uniformly in $n$ owing to \eqref{const: lambda} under (A'), (B') or (C'). 
By Theorem \ref{thm: EQUIVHINO}, it suffices to show the convergence of $(\E_n, \F_n)$ to $(\E_\infty, \F_\infty)$ in the sense of Definition \ref{thm: Mosco of Ch}. Conditions (N$1_*$) and (N$2_*'$) defined in Section \ref{sec: CNSF} are checked in the following arguments: 

 (N$1_*$): Let $f_n \to f_\infty$ weakly in $L^2$ with $\liminf_{n\to \infty}\mathcal E^n_{\lambda}(f_n)<\infty$. Since $({\sf Ch}_n, \F_n)$ converges to $({\sf Ch}_\infty, \F_\infty)$ in the Mosco sense \cite[Theorem 6.8]{GMS13}, and \eqref{ineq: 1} holds, it is obtained that   
\begin{align} \label{ineq: lambdaud}
{\sf Ch}_\infty(f_\infty) \le \liminf_{n \to \infty}{\sf Ch}_n(f_n)\le c\liminf_{n\to \infty}\mathcal E^n_{\lambda}(f_n)<\infty,
\end{align}
for some positive constant $c$ independent of $f_n$. For such $c$, see \eqref{ineq: 1} under Condition (A), (B) or (C). 
 This implies $f_\infty \in W^{1,2}(X,d, m_\infty)$.

(N$2_*'$): Let $n_k \uparrow \infty$ and $f_k \to f_\infty$ weakly in $L^2$ with $\sup_{k \in \N}\mathcal E^{n_k}_{\lambda+1}(f_{k})<\infty$ and $f_\infty \in W^{1,2}(X,d, m_\infty)$. Then, $f_k \to f_\infty$ weakly in $W^{1,2}$. Take $\mathcal C=H_{\Q_+}\mathscr A_{bs}$ (defined in Subsection \ref{sec: CNSF}) and $w \in \mathcal C$. 
By \eqref{eq: DIV}, the following holds:
\begin{align*}
|\mathcal E_{n_k}(f_k, w)-\mathcal E_{\infty}(f_\infty,w)|
&=\Bigl|\frac{1}{2}{\sf Ch}_{n_k}(f_k,w)-{\sf Ch}_\infty(f_\infty,w) \Bigr|+\Bigl|\int_{X_{n_k}}\b_{n_k}(f_k)wdm_{n_k}-\int_{X_\infty}\b_{\infty}(f_\infty)wdm_\infty\Bigr|
\\
&=\Bigl|\frac{1}{2}{\sf Ch}_{n_k}(f_k,w)-{\sf Ch}_\infty(f_\infty,w) \Bigr|+\Bigl|\int_{X_{n_k}}f_k\b_{n_k}(w)dm_{n_k}-\int_{X_\infty}f_\infty\b_{\infty}(w)dm_\infty\Bigr|\\
&\quad +
\Bigl|\int_{X_{n_k}}f_kw{\rm div} \b_{n_k}dm_{n_k}-\int_{X_\infty}f_\infty w {\rm div} \b_{\infty} dm_\infty\Bigr|
\\
&={\rm (I)}_k+{\rm (II)}_k+{\rm (III)}_k.
\end{align*} 
By the Mosco convergence of ${\sf Ch}_n$ to ${\sf Ch}_\infty$ proved by \cite{GMS13}, it can be seen that ${\rm (I)}_k \to 0$ as $k \to \infty.$
By $\b_n \to \b_\infty$ strongly in $L^2$, it holds that  (II)$_k \to 0$ as $k \to \infty$. 
By ${\rm div}\b_n \to {\rm div} \b_\infty$ strongly in $L^2$, 
the quantity (III)$_k \to 0$ also goes to zero. 
\qed



Theorem \ref{thm: mthm4} is now proved.
\\
{\it Proof of Theorem \ref{thm: mthm4}}. 
We prove (i). 
Let $\rho_n$ be the non-zero density of the invariant measure $\mu_n$ with $\|\rho_n\|_2=1$ for any $n \in \N$. By  Theorem \ref{thm: mthm1} for (A') and, by Corollary \ref{cor: grad} for (C'), it holds that  
\begin{align}
\sup_{n \in \N}\|{\rho}_n\|_{1,2} <\infty.
\end{align} 
By \cite[Theorem 6.3, Proposition 6.7]{GMS13},
 there is a converging subsequence $\{\rho_{n'}\}$ with its limit $\rho_\infty' \in W^{1,2}(X,d,m_\infty)$  so that
$$\rho_{n'} \to \rho_\infty', \quad \text{strongly in $L^2$ and weakly in $W^{1,2}$}.$$
For simplicity of notation, we write $\rho_{n}$ shortly for $\rho_{n'}$.
We first show that  $\rho'_\infty$ is an invariant measure for $(L_\infty, \mathcal D(L_\infty))$.
For any $\phi \in H_{\Q_+}\mathscr A_{bs}$, the following holds:
\begin{align} \label{eq: INVM}
0=\mathcal \E_n(\rho_n,\phi)
=\int_{X_n} \Bigl( \langle \nabla \rho_{n}, \nabla \phi \rangle +\b_n (\rho_{n})\phi \Bigr)dm_n \notag
& = \int_{X_n} \Bigl( \langle \nabla \rho_n, \nabla \phi \rangle -\rho_n\b_n(\phi)-{\rm div}\b_n\rho_n \phi\Bigr)dm_n \notag
\\
&\to  \int_{X} \Bigl( \langle \nabla \rho_\infty', \nabla \phi \rangle -\rho_\infty'\b_\infty(\phi)-{\rm div}\b_\infty\rho'_\infty\phi\Bigr)dm_\infty  \notag
\\
&=\mathcal E_\infty (\rho'_\infty, \phi).
\end{align}
Since $ H_{\Q_+}\mathscr A_{bs}$ is dense in $W^{1,2}(X,d,m)$, it holds that 
$$\mathcal E_\infty (\rho'_\infty, \phi)=0, \quad \forall \phi \in W^{1,2}(X,d,m).$$
This implies that 
$$\int_X L_\infty (\phi) \rho'_\infty dm_\infty=0, \quad \forall \phi \in \mathcal D(L).$$
Thus, $\rho_\infty'$ is an invariant measure for $(L_\infty, \mathcal D(L_\infty))$. Since $(L_\infty, \mathcal D(L_\infty))$ has a unique invariant measure by Theorem \ref{thm: mthm2}, it holds that $\rho'_\infty$ must coincide with $c\rho_\infty$ with some constant $c>0$. Since $\|\rho_n\|_2=1$ for any $n \in \EN$ and $\rho_n$ converges to $\rho_\infty'$ strongly in $L^2$ up to subsequences, we conclude that $c=1$ and the limit $\rho'_\infty$ is independent of choice of subsequences. Consequently, $\rho_n$ converges to $\rho_\infty$ strongly in $L^2$ and weakly in $W^{1,2}$. 

We prove (ii). Let $\rho_n$ be the density of the non-zero invariant measure $\mu_n$ with $\|\rho_n\|_2=1$ for any $n \in \N$.  Therefore, there is a converging subsequence $\{\rho_{n'}\}$ and its limit $\rho' \in L^2(X,m_\infty)$ so that 
 \begin{align} \label{prop: wconv}
 \rho_{n'} \to \rho' \ \text{weakly in $L^2$ as $n' \to \infty$}.
 \end{align}
By Proposition \ref{lem: Hino_Conv}, for any $f_n \to f_\infty$ strongly in $L^2$, it holds that $T_t^n f_n \to T_t^\infty f_\infty$ strongly in $L^2$. Note that, if $f_n=f|_{\iota_n(X_n)}$ for $f \in {\rm Lip}_{b}(X)$ for each $n \in \EN$, then $f_n \to f_\infty$ strongly in $L^2$. Note also that $m_n(f_ng_n) \to m_\infty(f_\infty g_\infty)$ whenever $f_n \to f_\infty$ strongly in $L^2$ and $g_n \to g_\infty$ weakly in $L^2$. Therefore, by \eqref{prop: wconv} and the invariance of $\rho_{n'}dm_{n'}$ with respect to $\{T_t^{n'}\}$, it holds that, for any $t>0,$ and $f \in {\rm Lip}_b(X)$, 
\begin{align*}
\int_{X_\infty} T_t^\infty (f) \rho'_\infty dm_\infty = \lim_{n' \to \infty} \int_{X_{n'}} T_t^{n'} (f)\rho_{n'} dm_{n'} 
= \lim_{n' \to \infty} \int_{X_{n'}}f \rho_{n'} dm_{n'} 
= \int_{X_\infty} f \rho_\infty' dm_\infty. 
\end{align*}
Noting that ${\rm Lip}_b(X)$ is dense in $L^2(X,m_\infty)$, if $\rho_\infty'$ is non-zero, then $\rho_\infty'$ must coincide with $c\rho_\infty$ with some constant $c>0$ owing to the uniqueness of the invariant measure $\mu_\infty$. 
Now it is shown that $\rho_\infty'$ is non-zero, which is proved by contradiction. Assume $\rho_\infty'=0$. Then $\rho_n$ has to converge to zero in $L^2$.
Thus, for any $t>0$, the following holds:
\begin{align} \label{eq: 1=0}
1=\int_{X_{n'}}\rho_{n'}^2 dm_{n'} =\int_{X_{n'}}T^{n'}_{t}(\rho_{n'})\rho_{n'} dm_{n'}  
= \|T^{n'}_{t} \rho_{n'}\|_2\|\rho_{n'}\|_2 
\le \|T^{n'}_{t}\|_{2 \to 2}\|\rho_{n'}\|^2_2 
& \to 0 \quad (n' \to \infty),
\end{align} 
whereby the last equality follows from the fact that $\sup_{n \in \N}\|T_{t_1}^n\|_{2 \to 2}<\infty$ by Proposition \ref{prop: L2}. Thus, \eqref{eq: 1=0} contradicts $\rho_\infty'=0$.
Consequently, $\rho_{n'} \to \rho'_\infty$ weakly in $L^2$ with $\rho'_\infty$ non-zero. Thus, $\rho'_\infty=c\rho_\infty$ for some constant $c>0$. This implies $A \subset {\rm Inv}_{> 0}(L_\infty)$, whereby $A$ denotes the set of limit points of $\{\rho_n\}$ with respect to the $L^2$-weak convergence.

Finally, it is proved that $\mu_n$ converges weakly to $\mu_\infty$ in the case of (i), and $\mathcal A \subset {\rm INV}_{> 0}(L_\infty)$, whereby $\mathcal A$ is the set of limit points of $\{\mu_n\}$ with respect to the weak convergence of measures in the case of (ii). 
It has been proved that $\rho_n$ converges to $\rho_\infty$ strongly in $L^2$ in the case of (i). 
Since $X_n$ converges to $X_\infty$ in the pmG sense, $m_n$ converges to $m_\infty$ weakly in $\mathcal P(X)$. By noting that $m_n(f_ng_n) \to m_\infty(f_\infty g_\infty)$ whenever $f_n \to f_\infty$ strongly in $L^2$ and $g_n \to g_\infty$ weakly in $L^2$,  it can be seen that, for any bounded continuous function $f \in C_b(X)$, 
\begin{align*}
\mu_n(f)&=\int_{X_n} f \rho_n dm_n  \to \int_{X_\infty} f \rho_\infty dm_\infty= \mu_\infty(f).
\end{align*}
In the case of (ii), taking a $L^2$-weak converging subsequence from $\{\rho_n\}$ and doing the same argument above, we can prove $\mathcal A \subset {\rm INV}_{> 0}(L_\infty)$.
Hence, the proof has been completed.  
\qed
 \section{Examples} \label{sec: exa} 
In this section, several examples for the main theorems are given. 
There are various concrete examples of non-smooth metric measure spaces satisfying RCD conditions. See Ricci limit spaces (Sturm \cite{Sturm06, Sturm06-2}, Lott--Villani \cite{LV09}, Cheeger--Colding \cite[Example 8]{CC97}), Alexandrov spaces (Petrunin, Zhang--Zhu \cite{Pet11, ZZ10}), warped products and cones (Ketterer \cite{Ket14, Ket14a}), quotient spaces (Galaz-Garc\'ia--Kell--Mondino--Sosa \cite{GKMS17}), stratified spaces (Bertrand--Ketterer--Mondello--Richard \cite{BKMR18}) and infinite-dimensional spaces such as Hilbert spaces with log-concave measures (Ambrosio--Savar\'e--Zambotti \cite{ASZ09}). We refer readers who are interested in non-smooth examples to these references. 
 
In this section, we give two examples of derivations $\b^n$ for the stability of invariant measures in Theorem \ref{thm: mthm4}. 
For any $f \in W^{1,2}(X,d,m)$, the derivation operator of the gradient type $\b_f$ is defined in the following manner:
$$\b_f(g)=\langle \nabla f, \nabla g\rangle, \quad g \in W^{1,2}(X,d,m).$$
It can be checked that $\b_f$ is an $L^2$-derivation with the point-wise norm $|\b_f|=|\nabla f|$ (see e.g., Gigli \cite{G16} for detail).
The strong $W^{1,2}$-convergence of $f_n \in W^{1,2}(X_n, d_n, m_n)$ is a sufficient condition for the $L^2$-strong convergence of $\b_{f_n}$ (recall the definition of the $W^{1,2}$-convergence in Definition \ref{defn: WConv}):
\begin{thm}$(${\rm \cite[Theorem 6.4]{AST16}}$)$ \label{defn: CONVD} \normalfont
Let $(X_n, d_n, m_n, \x_n)$ be a sequence of p.m.m.\ spaces with RCD$(K,\infty)$ condition. Assume that $(X_n,d_n,m_n, \x_n)$ converges to $(X_\infty, d_\infty, m_\infty, \x_\infty)$ in the pmG sense.  If $f_n \in W^{1,2}(X_n, d_n, m_n)$ converges strongly in $W^{1,2}$ to $f_\infty \in W^{1,2}(X_\infty, d_\infty, m_\infty)$, then $\b_{f_n}$  converges strongly in $L^2$ to $\b_{f_\infty}$.
\end{thm}
Using Theorem \ref{defn: CONVD}, derivations converging strongly in $L^2$ are constructed. See \cite[Example 6.6]{AST16} and \cite[Example 8.2]{S17}.
\begin{exa} {\bf (Derivation associated with resolvents)} \normalfont \label{thm: EXNS}
Let $(X_n, d_n, m_n, \x_n)$ be a sequence of p.m.m.\ spaces with RCD$(K,\infty)$ condition and $m_n(X_n)=1$ for $n \in \N$. Assume that $(X_n,d_n,m_n, \x_n)$ converges to $(X_\infty, d_\infty, m_\infty, \x_\infty)$ in the pmG sense. Let $\{G^n_\lambda\}$ and $\{H_t^n\}$ 
be the resolvent and the semigroup associated with Cheeger energy ${\sf Ch}_n$ respectively. 
Let $g_n \in L^\infty(m_n)$ with $\sup_{n \in \N}\|g_n\|_\infty<\infty$ and $g_n$ converge to $g_\infty \in L^\infty(m_\infty)$ strongly in $L^2$. 
 Let $f_n=G_\lambda g_n$ and $\b_n=\b_{f_n}$ for $n \in \EN$. Then $\b_n$ converges to $\b_\infty$ strongly in $L^2$, and $\sup_{n \in \N}(\|\b_n\|_\infty+\|{\rm div}\b_n\|_\infty)<\infty$.
\end{exa}

Another example is presented, which is given in terms of eigenfunctions of Laplacian according to \cite[Example 6.7]{AST16}. See also \cite[Example 8.3]{S17}. 
\begin{exa} {\bf (Derivation associated with eigenfunctions of Laplacian)} \normalfont \label{exa: EXNS}\\
Let $K>0$.
Let $\mathcal X_n=(X_n,d_n,m_n,\x_n)$ be an RCD$(K,\infty)$ space for all $n \in \N$ converging to $(X_\infty, d_\infty, m_\infty, \x_\infty)$ in the pmG sense. Let $u_n$ be a normalized eigenfunction $\int_{X_n}u_n^2dm_n=1$ of the generator $-\Delta_n$ associated with ${\sf Ch}_n$ with $-\Delta_nu_n=\lambda u_n$ for some $\lambda \in \R_{\ge 0}$. By \cite[Proposition 6.7]{GMS13}, $-\Delta_n$ has discrete spectra $\{\lambda_n^k\}_{k=1}^\infty$ (non-decreasing order) with the eigenfunctions $\{u_n^k\}_{k \ge 0}$. 
By \cite[Theorem 7.8]{GMS13}, $\lambda_n^k$ converges to $\lambda_\infty^k$, and $u_n^k$ converges to $u_\infty^k$ strongly in $L^2$ if the limit eigenvalue is simple (if not simple, a convergence subsequence can be extracted). By ${\sf Ch}_n(u_n^k)=\lambda^k_n \to \lambda_\infty^k={\sf Ch}_\infty(u_\infty^k)$, it holds that $u_n^k$ converges to $u_\infty^k$ strongly in $W^{1,2}$, which implies $\b_{u_n^k}$ converges to $\b_{u_\infty^k}$ strongly in $L^2$ by Theorem \ref{defn: CONVD}. If, furthermore,  $X_n$ is RCD$^*(K,N)$ with $\sup_{n \in \N}{\rm diam}(X_n)<\infty$, or with $\inf_{n \in \N}\inf_{x \in X_n}m_n(B_r(x))>0$ for any fixed $r>0$, then $\sup_{n \in \N}(\|\b_{u_n^k}\|_\infty+\|{\rm div}\b_{u_n^k}\|_\infty)<\infty$ for any $k \in \N$.  

\end{exa}
\vspace{5mm}
\noindent
Kohei Suzuki
\vspace{2mm}
\\ 
Scuola Normale Superiore di Pisa, 
\vspace{1mm}
\\ 
Piazza dei Cavalieri 7, 
\vspace{1mm}
\\ 
56126 Pisa PI, Italy.
\vspace{1mm} 
\\ 
E-mail: kohei.suzuki@sns.it
\end{document}